\documentclass[12pt]{article}%
\usepackage{amsfonts}
\usepackage{fancyhdr}
\usepackage{comment}
\usepackage[a4paper, top=2.5cm, bottom=2.5cm, left=2.2cm, right=2.2cm]%
{geometry}
\usepackage{times}
\usepackage{color}
\usepackage{amsmath,amsmath,amssymb}
\usepackage{changepage}
\usepackage{amssymb}
\usepackage{graphicx}%
\usepackage{graphicx,amssymb,enumerate,hyperref,subcaption,graphicx,amsthm,amsmath,amssymb}
\usepackage{bbm}
\usepackage{epstopdf}
\usepackage{algorithm}
\usepackage{rotating}
\usepackage {tikz}
\usetikzlibrary {positioning}

\newtheorem{theorem}{Theorem}

\newtheorem{corollary}[theorem]{Corollary}

\newtheorem{lemma}[theorem]{Lemma}

\newtheorem{proposition}[theorem]{Proposition}
\newtheorem{remark}[theorem]{Remark}

\usepackage{accents}
\newlength{\dhatheight}

\newcommand{\bfm}[1]{\ensuremath{\mathbf{#1}}}

     \def\CC{\mathbb{C}}
     
\def\be{e}     \def\EE{\mathbb{E}}

     \def\OO{\mathbb{O}}
     \def\PP{\mathbb{P}}
     
     \def\RR{\mathbb{R}}
\def\bs{\bfm s}

\def\bx{\bfm x}

\def\calC{{\cal  C}} 
 
\def\calE{{\cal  E}}

\def\calK{{\cal  K}} 
\def\calL{{\cal  L}} 
\def\calM{{\cal  M}} 
\def\calN{{\cal  N}} 
 
\def\calP{{\cal  P}}

\def\calS{{\cal  S}} 
\def\calT{{\cal  T}}

\def\calX{{\cal  X}}

\newcommand{\bfsym}[1]{\ensuremath{\boldsymbol{#1}}}

            \def\bDelta { {\Delta}}

              \def\bSigma{\bfsym \Sigma}

 \def\brho   {\bfsym {\rho}}
 
 \def\bxi{\bfsym {\xi}}
 
\DeclareMathOperator{\rank}{rank}

\DeclareMathOperator{\tr}{tr}

\newdimen\biblioindent    \biblioindent=30pt

 at 8truept

\numberwithin{equation}{section}

\begin{document}

%\begin{frontmatter}
\title{Confidence Region of Singular Subspaces for Low-rank Matrix Regression}
%\runtitle{Confidence interval of singular vectors}
%\runauthor{Dong Xia}

%\author{\fnms{Vladimir} \snm{Koltchinskii} \thanksref{t1}\ead[label=e1]{vlad@math.gatech.edu}}
%\address{686 Cherry Street, Atlanta, GA30332, USA.\\
%\printead{e1}}

\author{Dong Xia\footnote{Most of the manuscript was written when the author was affiliated with Columbia University.}\\ Hong Kong University of Science and Technology}

%\address{Columbia University, New York, NY10027, USA.}

%\printead{e2}

%\thankstext{t1}{Supported in part by NSF Grants DMS-1509739, DMS-1207808, CCF-1523768 and CCF-1415498}
%\thankstext{t2}{Supported in part by NSF Grants DMS-1207808, CCF-1523768.}

%\affiliation{Columbia University}
\date{(\today)}

\maketitle

\begin{abstract}
Low-rank matrix regression refers to the instances of recovering a low-rank matrix based on specially designed measurements and the corresponding noisy outcomes. In the last decade, numerous statistical methodologies have been developed for efficiently recovering the unknown low-rank matrices. However, in some applications, %such as quantum state tomography, phase retrieval and blind deconvolution, 
the unknown singular subspace is scientifically more important than the low-rank matrix itself. In this article, we revisit the low-rank matrix regression model and introduce a two-step procedure to construct confidence regions of the singular subspace. The procedure involves the de-biasing for  the typical low-rank estimators after which we calculate the empirical singular vectors. We investigate the distribution of the joint projection distance between the empirical singular subspace and the unknown true singular subspace. We specifically prove the asymptotical normality of the joint projection distance with data-dependent centering and normalization when $r^{3/2}(m_1+m_2)^{3/2}=o(n/\log n)$ where $m_1, m_2$ denote the matrix row and column sizes, $r$ is the rank and $n$ is the number of independent random measurements. Consequently, we propose data-dependent confidence regions of the true singular subspace which attains any pre-determined confidence level asymptotically. In addition, non-asymptotical convergence rates are also established. Numerical results are presented to demonstrate the merits of our methods. 
\end{abstract}

\section{Background and motivation}\label{sec:intro}
%Low-rank matrix regression aims to estimating a low-rank matrix based on specially designed measurements and the corresponding outcomes. 
Let ${M}\in\RR^{m_1\times m_2}$ be an unknown low-rank matrix with $r=\rank({M})\ll \min(m_1,m_2)$
whose singular value decomposition is written as ${M}={U}{\Lambda} {V}^\top$ with ${U}\in\RR^{m_1\times r}, {V}\in\mathbb{R}^{m_2\times r}$ being ${M}$'s left and right singular vectors. The diagonal matrix ${\Lambda}={\rm diag}(\lambda_1,\ldots,\lambda_r)$ with $\lambda_1\geq\cdots\geq \lambda_r>0$ consists of its corresponding singular values. The goal of matrix regression is to recover ${M}$ from a set of measurements and noisy outcomes. It has been intensively studied in the last decade. See, e.g., \cite{candes2011tight}, \cite{koltchinskii2011nuclear}, \cite{koltchinskii2015optimal}, \cite{negahban2011estimation}, \cite{rohde2011estimation}, \cite{gross2011recovering} and references therein. In general, the most popular framework for studying this problem is the so-called trace regression model
which can be described by a random pair $({X},y)$ with $X\in\RR^{m_1\times m_2}$ and $y\in\RR$ satisfying
\begin{gather}\label{eq:tr_reg}
y=\tr({M}^{\top}{X})+\xi
\end{gather}
where the noise $\xi$ is independent with ${X}$ and $\xi\sim\calN(0,\sigma_\xi^2)$. Given i.i.d. copies $\{({X}_i, y_i)\}_{i=1}^n$, the goal is to recover ${M}$ with both computational and statistical efficiency. 

Many applications can be formulated by the trace regression model (\ref{eq:tr_reg}) with $\rank({M})\ll \min(m_1,m_2)$. Among these applications, the following three examples are the most related with the motivation of this article. 
\smallskip

{\it Quantum state tomography.} In {\it quantum computing} and {\it quantum communication}, it is often necessary to recover the state (called {\it quantum state}) of quantum systems.  The {\it pure quantum state} is usually represented by a state vector ${v}\in\CC^{m}$ where $m=2^b$ with $b$ being the number of qubits in the quantum system. See \cite{koltchinskii2015optimal}, \cite{xia2016estimation}, \cite{xia2017estimation}, \cite{gross2011recovering}, \cite{carpentier2015uncertainty} and references therein. In {\it quantum mechanics}, quantum systems are allowed to have {\it mixed state} which is a statistical ensemble of pure states. Basically, it says that the quantum system is in pure state ${v}_k\in\CC^m$ with probability $p_k>0$ for $k=1,\cdots,r$ so that $\sum_{k=1}^r p_k=1$. A mixed state is commonly studied, for simplicity, as a {\it density matrix} which is defined as 
$$
\brho=\sum_{k=1}^r p_k ({v}_k {v}_k^{\dagger})
$$ 
where $v^{\dagger}$ denotes the conjugate transpose of $v$. The density matrix is therefore  self-adjoint and positively semi-definite. Meanwhile, its trace $\tr(\brho)=1$ and its rank $\rank(\brho)\leq r$. Then, {\it quantum state tomography} refers to the recovery of $\brho$ by i.i.d. pairs of special measurements and the noisy outcomes $\{({X}_i, y_i)\}_{i=1}^n$ satisfying model (\ref{eq:tr_reg}). However, recovering the density matrix $\brho$ is generally not the ultimate goal in quantum state tomography. Of course, determining the pure states $\{{v}_k\}_{k=1}^r$ is physically more important.

{\it Phase retrieval.} The goal of phase retrieval is to recover an unknown vector $a\in\RR^m$ from noisy outcomes of the squared magnitudes of $a$'s linear measurements. Formally,  the outcome $y_i$ is written as $y_i=(a^{\top}x_i)^2+\xi_i$ where the measurement vector $x_i\in\RR^m$ can be specially designed or randomly chosen and the noise $\xi_i\sim\calN(0,\sigma_\xi^2)$ has a Gaussian distribution. Given $\{(x_i,y_i)\}_{i=1}^n$, the goal is then to recover the signal $a$. 
 It has attracted a lot of attention especially in X-ray crystallography (\cite{klibanov1995phase}). In the recent years, \cite{candes2013phaselift} and \cite{candes2014solving} proposed computationally efficient approaches for solving this problem by convex optimization. The basic idea is to transform a measurement vector $x_i$ into a measurement matrix via ${X}_i=x_ix_i^{\top}\in\RR^{m\times m}$. Instead of estimating the signal $a$ directly, we can solve for the low-rank matrix ${M}=aa^{\top}$. Consequently, phase retrieval is reformulated as the matrix regression model (\ref{eq:tr_reg}) with rank $r=1$. Although being formulated as  matrix regression, the more important goal of phase retrieval is to recover ${M}$'s column space.  
 
 {\it Blind deconvolution.} Blind deconvolution refers to the problem of recovering two unknown vectors from their circular convolution.  Let $a\in\RR^{m_1}$ and $b\in\RR^{m_2}$ be the two unknown vectors. The measurement can be expressed as a noisy bilinear function of the rank-$1$ matrix ${M}=ab^{\top}$. Given two measurement vectors $s\in\RR^{m_1}, t\in\RR^{m_2}$, the outcome $y=\tr\big({M}^{\top}{X}\big)+\xi$ with the measurement matrix ${X}=st^{\top}$ and the noise $\xi\sim\calN(0,\sigma_\xi^2)$. It is thus translated into the trace regression model (\ref{eq:tr_reg}) and was solved by convex programming in \cite{ahmed2014blind}. Similarly, the ultimate goal of blind deconvolution is to recover ${M}$'s row and column space, rather than the matrix ${M}$ itself. Toward that end, the leading left and right singular vector of the low-rank estimator $\hat{M}$ are usually computed. See \cite{ahmed2014blind} and \cite{ma2017implicit} for more details.  
 \smallskip
 
 In the aforementioned applications of low-rank matrix regression, the underlying ``signal" of interest is the singular subspace of the unknown low-rank matrix. The objective of this article is to propose an approach for constructing the confidence regions of the singular subspace for low-rank matrix regression. Essentially, it is corresponding to the inference of the pure states in quantum state tomography, and the inference of the signal directions in phase retrieval and blind deconvolution. In general,
 the statistical inference for low-rank matrix regression model is subtle. The opening question is to choose the parameters of interest for the investigation. In \cite{carpentier2015uncertainty} and \cite{carpentier2018adaptive}, the authors proposed confidence regions of the matrix ${M}$ with respect to the matrix Frobenius norm. In \cite{cai2016geometric} and \cite{carpentier2015iterative}, the confidence intervals for $M$'s entries are  established. Note that the matrix Frobenius norm is equivalent to the $\ell_2$-norm of the vectorization of a matrix. As a result, the matrix structure can not be directly reflected. Similarly, the individual matrix entries also do not reflect the matrix geometry. In comparison, the statistical inference of the underlying singular subspace  is more important in the aforementioned applications. 
 
%  we believe that the statistical inference of the underlying singular subspace is more fundamentally important even for the general purposes because knowing the distribution of the empirical singular vectors could assist us in developing methods for the statistical inference of other ${M}$-related parameters such as ${M}$'s entries.  
 
 In this article, we propose a novel approach for the statistical inference of the singular subspace in low-rank matrix regression model. On a high level, the approach consists of two procedures. It begins with a statistically optimal estimator of the underlying low-rank matrix, where, for instance, the nuclear-norm penalized least squares estimator will be implemented. It is followed by a de-biasing treatment which outputs an unbiased estimator of the underlying low-rank matrix. Then, we compute the singular value decomposition and extract the corresponding left and right singular vectors to serve as the final estimator of the singular vectors. The de-biasing procedure is essentially to re-randomize the low-rank estimate from the first step.  
 We characterize the bias of the empirical singular vectors. The bias depends on the sample size, ambient dimension and the inverse of true singular values. All the explicit constant factors are developed. With near-optimal sample size requirement, we prove the asymptotical normality, when the bias is subtracted,  of the joint projection distance between the empirical singular subspace and the true singular subspace. This result is still insufficient for constructing the confidence regions of the singular subspace if the bias is unknown. To ensure sharp estimation of the bias such that its error is dominated by the standard deviation of the joint projection distance, it turns out that we require the sample size to be somewhat larger than the typically optimal conditions. But the sample size requirement is still much smaller than the ambient dimension of the matrix space. Analogous phenomenon also exists in the statistical inference for sparse vector linear regression with unknown design. See \cite{zhang2014confidence} and \cite{cai2017confidence} for more details. Based on the normal approximation of a novel data-dependent statistics, we construct the confidence region of $M$'s singular subspace  which achieves any pre-determined confidence level asymptotically. Numerical experiments show that the proposed method works extremely well. 

The statistical inference of low-dimensional structure in (ambient) high-dimensional space has attracted tremendous amount of attention in the recent years, especially for the statistical inference of sparse vector linear regression. 
Statistically efficient procedures have been developed in the recent decades to handle the challenges posed by 
the high dimensionality.  For instance, it includes the $\ell_1$-penalization for sparse linear regression (see \cite{tibshirani1996regression}, \cite{yuan2006model}, \cite{zou2005regularization}, \cite{candes2008enhancing} and references therein) and the matrix nuclear-norm penalization for low-rank matrix regression (see \cite{candes2010matrix}, \cite{negahban2011estimation},\cite{rohde2011estimation}, \cite{koltchinskii2011nuclear} and references therein). Under certain regularity conditions, those methods are guaranteed to be statistically efficient. It means that the minimax optimal rates of the estimation error, usually relevant to the degrees of freedom, are attainable.  However, the statistical inference for the aforementioned high dimensional problems is difficult.
  In several recent papers \cite{javanmard2014confidence}, \cite{zhang2014confidence}, \cite{lockhart2014significance}, \cite{cai2017confidence}, a post-processing approach was proposed which de-biases the $\ell_1$-penalized least squares estimator. It was shown that the statistical inference usually requires stronger conditions for the sparse vector linear regression. In recent years, the statistical inference of the singular subspaces and eigen subspaces is popular in the statistics community. 
A lot of efforts have been put into studying the asymptotic property of the {\it principle component analysis} (PCA). For instance,
 the normal approximation of the eigenvectors of the sample covariance matrix has been studied in \cite{koltchinskii2016asymptotics}, \cite{koltchinskii2015normal} and \cite{koltchinskii2017new}. In both \cite{koltchinskii2015normal} and \cite{koltchinskii2017new}, the data splitting trick is applied for estimating the bias of the empirical eigenvectors, which is critical since the bias of the empirical eigenvectors usually significantly dominates its standard deviation. 
In addition, a Bayesian approach for constructing  the confidence regions of the principle components is studied in \cite{silin2017bayesian}. A more sophisticated bias reduction framework by iterative bootstrap for the inference of PCA is proposed in \cite{koltchinskii2017asymptotically}.

 %\newpage
 %which achieves any pre-determined confidence level $\alpha\in(0,1)$ asymptotically.
%The motivation for quantifying the uncertainty of singular vectors instead of other parameters is to take advantage of matrix structure. Since matrix regression is equivalent to vector regression when both the unknown matrix and measurement matrices are vectorized, it is crucial to dig out the fundamental difference between matrix regression and vector regression: the ultimate goal of matrix regression is to recover the column and row spaces, in addition to the signal strength (singular values).  
The rest of the paper is organized as follows. In Section~\ref{sec:problem+main_results}, we explain important notations and introduce the basic assumptions. An overview of our main results is also provided in Section~\ref{sec:problem+main_results}. The two-step procedure for estimating the singular subspaces is given in Section~\ref{sec:estimator}. We present the theoretical performance of the proposed method in Section~\ref{sec:theory} where we provide the normal approximation of the joint projection distance between the empirical singular subspace and the true singular subspace. In Section~\ref{sec:CI}, we define the data-dependent confidence region which attains the pre-determined confidence level asymptotically. Numerical simulations are displayed in Section~\ref{sec:simulation}. In Section~\ref{sec:discussion}, we discuss about the dealing with unknown ranks and the double-sample-splitting trick which can avoid the loss of efficiency due to the data splitting. The proofs are postponed to Section~\ref{sec:proofs} and Section~\ref{sec:proof_additional}.

\section{Overview of main results}\label{sec:problem+main_results}
\subsection{Notations}\label{notationsec}
%We use boldfaced upper-case letters to denote matrices, and use the same letter in normal font with indices to denote its entries. 
For a matrix ${A}\in\RR^{m_1\times m_2}$, we denote by $\|{A}\|_{\rm F}$ its Frobenius norm and $\|{A}\|$ its operator norm. The nuclear norm of ${A}$ is denoted by $\|{A}\|_{\star}$, i.e., the sum of its singular values. 
Let ${\rm vec}({A})\in\RR^{m_1m_2}$ denote its vectorized version. Similarly, we denote by $\calM(\cdot)$ the inverse of ${\rm vec}(\cdot)$ such that $\calM\big({\rm vec}({A})\big)={A}$. Given ${B}\in\RR^{m_1\times m_2}$, we denote $\langle{A},{B} \rangle=\tr({A}^{\top}{B})$. We use $c_1,c_2,C_1,C_2,\cdots$ to represent absolute constants which might vary lines from lines during the proof and the statement of theorems. For two sequences of random variables $\{a_n\}_n, \{b_n\}_n$ which are positive {\it almost surely}, we write $a_n=O_P(b_n)$ to represent that there exists an absolute constant $C_1>0$ such that $\lim_{n\to\infty}\PP\big(a_n/b_n\geq C_1\big)=0$. We denote by $\OO^{m\times r}$ the set of $m\times r$ matrices whose columns are orthonormal. We write $\bar{m}=\max\{m_1,m_2\}$. 

\subsection{Nuclear-norm penalized low-rank estimation}
%Let ${M}\in\RR^{m_1\times m_2}$ be an unknown low-rank matrix with $r=\rank({M})\ll \min(m_1,m_2)$
%whose singular value decomposition is written as ${M}={U}{\Lambda} {V}^\top$ with ${U}\in\RR^{m_1\times r}, {V}\in\mathbb{R}^{m_2\times r}$ being ${M}$'s left and right singular vectors. The diagonal matrix ${\Lambda}={\rm diag}(\lambda_1,\ldots,\lambda_r)$ where $\lambda_1\geq\cdots\geq \lambda_r>0$ are its corresponding singular values. The goal of matrix regression is to recover ${M}$ from a set of measurements and noisy outcomes, which was intensively investigated in the last decade especially when ${M}$ has low-rank. See, e.g., \cite{candes2011tight}, \cite{koltchinskii2011nuclear}, \cite{koltchinskii2015optimal}, \cite{negahban2011estimation}, \cite{rohde2011estimation}, \cite{gross2011recovering} and references therein. In particular, this paper is focused on the trace regression model
%which is characterized by a random pair $({X},y)$ such that
%\begin{gather}\label{eq:tr_reg}
%y=\tr({M}^{\top}{X})+\xi
%\end{gather}
%where the noise $\xi$ is independent of ${X}$ and $\xi\sim\calN(0,\sigma_\xi^2)$. 
Given i.i.d. copies $\{({X}_i,y_i)\}_{i=1}^{2n}$ satisfying (\ref{eq:tr_reg}), 
it was shown in \cite{negahban2011estimation} and \cite{candes2011tight} that the matrix nuclear-norm penalized least squares estimator eq. (\ref{eq:nuclear_pen_est}), denoted as $\hat{M}^{\rm nuc}$, achieves the statistically optimal convergence rate:
\begin{equation}\label{eq:nuclear}
\|\hat{M}^{\rm nuc}-{M}\|_{\rm F}^2=O_P\bigg(\frac{\sigma_\xi^2r(m_1+m_2)}{n}\bigg)
\end{equation}
if the so-called {\it restricted isometry property} (RIP) or {\it restricted strong convexity} (RSC) hold.  Specifically, if ${X}$ has a sub-gaussian distribution, it was shown in \cite{candes2011tight} and \cite{negahban2011estimation} that the RIP and RSC hold with probability at least $1-c_1e^{-c_2\bar{m}}$ as long as $n\geq C_1r(m_1+m_2)\log\bar{m}$ where $\bar{m}=\max\{m_1,m_2\}$ and $c_1,c_2,C_1$ are absolute constants.

%Clearly, the rate (\ref{eq:nuclear}) is statistically optimal, where, however, leaves a fundamental unresolved question.
%As explained in Section~\ref{sec:intro}, the fundamental difference between matrix regression and vector regression is that, in matrix regression, the ultimate objective is to recover the unknown column space and row space of ${M}$, in addition to ${\Lambda}$, which is not directly reflected in
%the matrix Frobenius norm (\ref{eq:nuclear}) since it neglects the matrix structure. 
Our goal in this article is to estimate the singular subspace of ${M}$, i.e., the column space of ${U}$ and ${V}$, and conduct the statistical inference. An immediate approach is to take the singular vectors of $\hat {M}^{\rm nuc}$.
If we apply the famous Wedin's $\sin\Theta$ theorem \cite{wedin1972perturbation} or Davis-Kahan theorem \cite{davis1970rotation}, by (\ref{eq:nuclear}), we can get a naive bound 
\begin{align}\label{eq:U_nuclear}
{\rm dist}^2\big[(\hat{U}^{\rm nuc}, \hat{V}^{\rm nuc}), ({U},{V})\big]=O_P\bigg(\frac{\sigma_\xi^2}{\lambda_r^2}\cdot\frac{r(m_1+m_2)}{n}\bigg)
\end{align}
where $\hat{U}^{\rm nuc}$ and $\hat{V}^{\rm nuc}$ are $\hat{M}^{\rm nuc}$'s top-$r$ left and right singular vectors and
\begin{align}
{\rm dist}^2\big[(\hat{U}^{\rm nuc}, \hat{V}^{\rm nuc}),& ({U},{V})\big]\label{eq:def_dist}\\
=&\|\hat{U}^{\rm nuc}(\hat{U}^{\rm nuc})^{\top}-{U}{U}^{\top}\|_{\rm F}^2+\|\hat{V}^{\rm nuc}(\hat{V}^{\rm nuc})^{\top}-{V}{V}^{\top}\|_{\rm F}^2.\nonumber
\end{align}
The loss ${\rm dist}^2\big[(\hat{U}^{\rm nuc}, \hat{V}^{\rm nuc}), ({U},{V})\big]$ is usually called the (squared) joint projection distance between the empirical singular subspace and the true singular subspace. 
The naive bound (\ref{eq:U_nuclear}) is sub-optimal especially when $\lambda_1\geq \cdots\geq \lambda_{r-1}\gg \lambda_r$ in which case the inhomogeneity of the singular values is not reflected in (\ref{eq:U_nuclear}). Moreover, the bound (\ref{eq:U_nuclear}) is insufficient for constructing the confidence region of the singular vectors ${U}$ and ${V}$. We note that, from eq. (\ref{eq:U_nuclear}), $\hat{U}^{\rm nuc}$ and $\hat{V}^{\rm nuc}$ are nontrivial if $n\gg \frac{\sigma_\xi^2}{\lambda_r^2}\cdot r(m_1+m_2)$. In light of the standard sample size requirement $n\geq C_1r(m_1+m_2)$ for estimating ${M}$, it is therefore convenient for us to focus on the scenario that $\frac{\sigma_\xi}{\lambda_r}=O(1)$ for simplicity. Otherwise, we shall adjust the baseline of sample size requirement accordingly which involves $\beta:=\frac{\sigma_\xi}{\lambda_r}$. 

\subsection{Estimating the singular subspaces}
To construct the confidence region for ${U}$ and ${V}$, we propose a two-step procedure for estimating the column spaces of ${U}$ and ${V}$. In particular,
we focus on the standard Gaussian design where ${X}$ has i.i.d. standard Gaussian entries, i.e., $X_{ij}\stackrel{{\rm i.i.d.}}{\sim}\calN(0,1)$ for all $(i,j)\in[m_1]\times [m_2]$ where $[m]:=\{1,\ldots,m\}$.
The two-step procedure starts with a nuclear-norm penalized estimator $\hat{M}^{\rm nuc}$
and it is followed by a de-biasing step which produces a new estimator $\hat{M}$. Even though $\hat{M}$ loses the low-rank property, it is an unbiased estimator of ${M}$.
 Then, we compute $\hat{U}$ and $\hat{V}$ from the top-$r$ left and right singular vectors of $\hat{M}$ which serve as the final estimators of ${U}$ and ${V}$. %Denote by 
%$$
%\calP_{{U}{V}}=\left(
%\begin{array}{cc}
%{U}{U}^{\top}&{0}\\
%{0}&{V}{V}^{\top}
%\end{array}
%\right)
%\quad {\rm and}\quad 
%\calP_{\hat{U}\hat{V}}=\left(
%\begin{array}{cc}
%\hat{U}\hat{U}^{\top}&{0}\\
%{0}&\hat{V}\hat{V}^{\top}
%\end{array}
%\right).
%$$
%Therefore, 
%$$
%\|\calP_{{U}{V}}-\calP_{\hat{U}\hat{V}}\|_{\rm F}^2=\|\hat{U}\hat{U}^{\top}-{U}{U}^{\top}\|_{\rm F}^2+\|\hat{V}\hat{V}^{\top}-{V}{V}^{\top}\|_{\rm F}^2.
%$$
The joint space spanned by the columns of $\hat U$ and $\hat V$ is called the empirical singular subspace. The de-biasing procedure is critical for characterizing the distribution of the empirical singular subspace. Note that the initial estimate $\hat M^{\rm nuc}$ is already close to the true low-rank matrix in the Euclidean norm. However, it is exceedingly difficult to characterize the distribution of $\hat M^{\rm nuc}$'s singular vectors since it comes from an output of a convex programming. Therefore, the de-biasing step is essentially to re-randomize $\hat M^{\rm nuc}$ with pre-determined distributions. The benefit of re-randomization is that we are able to characterize the distribution of the empirical singular subspaces without affecting the convergence rates. 

Since the empirical singular vectors $\hat U$ and $\hat V$ are determined up to the multiplication of an orthonormal matrix, we focus on the empirical spectral projectors $\hat U\hat U^{\top}$ and $\hat V\hat V^{\top}$. We study the loss function ${\rm dist}^2[(\hat U,\hat V), (U,V)]$ which is as defined in (\ref{eq:def_dist}) with
$$
{\rm dist}^2[(\hat U,\hat V), (U,V)]=\|\hat U\hat U^{\top}-UU^{\top}\|_{\rm F}^2+ \|\hat V\hat V^{\top}-VV^{\top}\|_{\rm F}^2.
$$ 
We will characterize the expected loss of the empirical singular subspace. 
If $n\geq C_1\big[\beta^2\bar{m}+\bar{m}r\log^2n\big]$ with $\beta=\frac{\sigma_\xi}{\lambda_r}$ for a large enough constant $C_1>0$, then we show that
\begin{align*}
\EE\ {\rm dist}^2\big[(\hat{U},\hat{V}),({U},{V})\big]
=&\sigma_\xi^2\|{\Lambda}^{-1}\|_{\rm F}^2\cdot\frac{2m_{\star}}{n}
+O\bigg((\beta\vee 1)^4\cdot \frac{r^{3/2}\bar{m}^2\log^{1/2}n}{n^2}\bigg)
\end{align*}
where $m_{\star}=m_1+m_2-2r$. Therefore, if $n\gg r^{3/2}\bar{m}\log^{2}\bar{n}$ when $\beta=O(1)$, we simply obtain 
\begin{equation}\label{eq:EEhatU-U}
\EE\ {\rm dist}^2\big[(\hat{U},\hat{V}),({U},{V})\big]=[2+o(1)]\sigma_\xi^2\|\Lambda^{-1}\|_{\rm F}^2\cdot \frac{m_{\star}}{n}.
\end{equation}
Comparing with the naive bound (\ref{eq:U_nuclear}) which is described only by the smallest singular value $\lambda_r$, our bound (\ref{eq:EEhatU-U}) shows that those larger singular values play the same role. 
 Moreover, the bound (\ref{eq:EEhatU-U}) also establishes the exact constant factor. 

We then study the standard deviation of the loss ${\rm dist}^2\big[(\hat{U},\hat{V}),({U},{V})\big]$. Similarly, if $n\gg r\bar{m}\log^2n$ when $\beta=O(1)$,
 we prove that its standard deviation
\begin{align*}
\big|{\rm dist}^2\big[(\hat{U},\hat{V}),({U},&{V})\big]-\EE{\rm dist}^2\big[(\hat{U},\hat{V}),({U},{V})\big] \big|
=O_P\bigg(\sigma_\xi^2\|\Lambda^{-2}\|_{\rm F}\cdot \frac{(r\bar{m}\log n)^{1/2}}{n}\bigg).
\end{align*}
It shows that the standard deviation of ${\rm dist}^2\big[(\hat{U},\hat{V}),({U},{V})\big]$ is only of $O(\bar{m}^{-1/2})$ fraction of its expectation. This typical phenomenon of the empirical singular subspaces is the foremost reason why the statistical inference of the singular subspaces is difficult.  
Indeed, if we investigate the normal approximation of ${\rm dist}^2\big[(\hat{U},\hat{V}),({U},{V})\big]$, we can show that 
\begin{align}\label{eq:supx_hatP-P_normal}
\sup_{x\in\RR}\bigg|\PP\bigg\{&\frac{{\rm dist}^2\big[(\hat{U},\hat{V}),({U},{V})\big]-\EE{\rm dist}^2\big[(\hat{U},\hat{V}),({U},{V})\big]}{\sqrt{8}\sigma_{\xi}^2\|{\Lambda}^{-2}\|_{\rm F}\cdot \frac{m_{\star}^{1/2}}{n}}\leq x\bigg\}-\Phi(x) \bigg|\longrightarrow 0
\end{align}
as long as  $\bar{m},n\to\infty$ and $\frac{r^2\bar{m}\log^3n}{n}\to 0$ when $\beta=O(1)$. Here, $\Phi(x)$ represents the {\it cumulative distribution function} of the standard normal distribution. The sample size is optimal up to the rank $r$ and the logarithmic factor. 
By eq. (\ref{eq:supx_hatP-P_normal}), it suffices to estimate the expected loss and prove the normal approximation of ${\rm dist}^2[(\hat U,\hat V),(U,V)]$ with explicit centering and normalization terms. However, since the centering term is much larger (with a factor of $\bar{m}^{1/2}$) than the nomalization term, a small error (such as the constant factor in eq. (\ref{eq:EEhatU-U}) ) in estimating $\EE{\rm dist}^2\big[(\hat{U},\hat{V}),({U},{V})\big]$ is likely to ruin the overall asymptotical distribution. 

In view of eq. (\ref{eq:EEhatU-U}), we replace $\EE{\rm dist}^2\big[(\hat{U},\hat{V}),({U},{V})\big]$ with $\sigma_\xi^2\|{\Lambda}^{-1}\|_{\rm F}^2\cdot \frac{2m_{\star}}{n}$ in (\ref{eq:supx_hatP-P_normal}). Surprisingly, we will prove that
\begin{align}\label{eq:supx_hatP-sigma_normal}
\sup_{x\in\RR}\bigg|\PP\bigg\{&\frac{{\rm dist}^2\big[(\hat{U},\hat{V}),({U},{V})\big]-\sigma_\xi^2\|{\Lambda}^{-1}\|_{\rm F}^2\cdot \frac{2m_{\star}}{n}}{\sqrt{8}\sigma_{\xi}^2\|{\Lambda}^{-2}\|_{\rm F}\cdot \frac{m_{\star}^{1/2}}{n}}\leq x\bigg\}-\Phi(x) \bigg|\longrightarrow 0
\end{align}
as long as  $\bar{m},n\to\infty$ and $\frac{r^{3/2}\bar{m}^{3/2}\log n}{n}\to 0$ when $\beta=O(1)$. Therefore, we obtain the normal approximation of the loss with explicit centering and normalization terms. 
Of course, the sample size requirement for (\ref{eq:supx_hatP-sigma_normal}) turns out to be stronger than that for (\ref{eq:supx_hatP-P_normal}). In Section~\ref{sec:theory}, we will explain that the sample size requirement in the formulation (\ref{eq:supx_hatP-sigma_normal}) indeed needs to grow as $O(\bar{m}^{3/2})$. Put it differently, this sample size requirement is generally not improvable unless sharper estimates are derived for the expected loss. Finally, we propose data-dependent estimates for $\sigma_\xi^2$, $\|\Lambda^{-1}\|_{\rm F}^2$ and $\|\Lambda^{-2}\|_{\rm F}$ and will prove the normal approximation of ${\rm dist}^2[(\hat U,\hat V),(U,V)]$ with data-dependent centering and normalization terms. 

\section{Methodology: de-biasing and  re-randomization}\label{sec:estimator}
Our method consists of two steps, each of which is implemented on one independent data sample. 
Suppose that i.i.d. copies $\{({X}_i,y_i)\}_{i=1}^{2n}$ satisfying (\ref{eq:tr_reg}) are available where the underlying matrix ${M}={U}{\Lambda}{V}^{\top}$ is unknown and has rank $r=\rank({M})\ll \min(m_1,m_2)$, our goal is to design estimators of ${U}$ and ${V}$. We will split the data into two independent samples: $\{(X_i, y_i)\}_{i=1}^{n}$ and $\{(X_i,y_i)\}_{i=n+1}^{2n}$. The first data sample is used to obtain an efficient low-rank estimate of $M$ and the second data sample is used for the bias correction. In Section~\ref{sec:discussion}, we discuss a simple double-sample-splitting trick which could avoid the loss of efficiency due to the sample splitting. 

Our first step is implemented by the famous nuclear-norm penalized least squares estimator (\cite{candes2011tight} and \cite{negahban2011estimation}). Given the first data sample $\{(X_i,y_i)\}_{i=1}^n$, the estimator is defined as 
\begin{equation}\label{eq:nuclear_pen_est}
\hat{M}^{\rm nuc}:=\arg\min_{{A}\in\RR^{m_1\times m_2}} \frac{1}{n}\sum_{i=1}^n\big(y_i-\tr({A}^{\top}{X}_i)\big)^2+\lambda \|{A}\|_{\star}
\end{equation}
where $\|\cdot\|_{\star}$ denotes the matrix nuclear norm.  The estimator relies on solving the convex program where the nuclear norm penalization promotes low-rank solutions. It was shown in \cite{candes2011tight} and \cite{negahban2011estimation} that if $n\gg \bar{m}r$ with $\bar m=\max(m_1,m_2)$ and $\lambda=C_1\sigma_\xi\sqrt{\frac{\bar{m}}{n}}$ for some absolute constant $C_1>0$, then $\hat M^{\rm nuc}$
achieves the minimax optimal convergence rate in the matrix Frobenius norm (see eq. (\ref{eq:nuclear})). We note that the initial low-rank estimator is unnecessary to be always fixed to $\hat{M}^{\rm nuc}$. Actually, in the first step, any estimator $\hat{M}^{\rm init}$, such as the projection estimator \cite{xia2016estimation}, \cite{klopp2011rank} and the matrix Dantzig estimator \cite{xia2017estimation}, \cite{candes2011tight}, which achieves the statistically optimal convergence rate are all qualified to replace $\hat{M}^{\rm nuc}$. 

Although the estimator $\hat M^{\rm nuc}$ is statistically optimal, it is usually biased. Our second step is to implement the bias correction for $\hat M^{\rm nuc}$. 
Indeed, by utilizing the second data sample $\{({X}_i,y_i)\}_{i=n+1}^{2n}$, we compute a de-biased version of $\hat M^{\rm nuc}$ as
\begin{equation}\label{eq:debias_M}
\hat{M}=\hat{M}^{\rm nuc}+\frac{1}{n}\sum_{i=n+1}^{2n}\big(y_i-\tr({X}_i^{\top}\hat{M}^{\rm nuc})\big){X}_i
\end{equation}
Note that $\hat{M}^{\rm nuc}$ is independent with $\{({X}_i,\xi_i)\}_{i=n+1}^{2n}$, it is straightforward to check that $\EE\hat{M}={M}$ even though $\hat{M}$ has full rank almost surely. The idea of de-biasing was initially proposed for the statistical inference of sparse vector linear regression where the sample splitting  (\ref{eq:tr_reg}) is the simplest approach when the design of ${X}$ is known in advance. See, e.g., \cite{javanmard2014confidence}, \cite{zhang2014confidence}, \cite{lockhart2014significance}, \cite{cai2017confidence} and references therein. 
The de-biasing procedure (\ref{eq:debias_M}) can be viewed as a re-randomization of the initial low-rank estimator $\hat M^{\rm nuc}$. If we denote by $\bDelta={M}-\hat{M}^{\rm nuc}$, we can write 
\begin{equation}\label{eq:hatM=M+Z_1+Z_2}
\hat{M}={M}+\underbrace{\frac{1}{n}\sum_{i={n+1}}^{2n}\xi_i{X}_i}_{{Z}_1}+\underbrace{\Big(\frac{1}{n}\sum_{i={n+1}}^{2n}\tr(\bDelta^{\top}{X}_i){X}_i-\bDelta\Big)}_{{Z}_2}
\end{equation}
where $\bDelta$, $\{\xi_i\}_{i=n+1}^{2n}$ and $\{{X}_i\}_{i=n+1}^{2n}$ are mutually independent but ${Z}_1$ and ${Z}_2$ are dependent. If $\|\Delta\|_{\rm F}=o_P(\sigma_\xi)$ so that $Z_1$ dominates $Z_2$, then we can view $\hat M$ as a random perturbed (with i.i.d. entries) version of $M$. This explicit characterization of $\hat M$ is the reason why we can study the distribution of $\hat M$'s singular subspace.

Finally, we compute the top-$r$ left and right singular vectors of $\hat M$,
denoted by $\hat{U}$ and $\hat{V}$. They are our final estimators of ${U}$ and ${V}$. 

\begin{remark}
The first step in our method is important. Actually, an immediately unbiased estimator of $M$ is $\tilde M:=n^{-1}\sum_{i=1}^n y_i X_i$ which does not rely on any initial estimators. We can write 
$$
\tilde M= M +\underbrace{\frac{1}{n}\sum_{i=1}^n \xi_i X_i}_{\tilde Z_1} +\underbrace{\Big(\frac{1}{n}\sum_{i=1}^n\tr(M^{\top}X_i)X_i-M\Big)}_{\tilde Z_2}
$$
which has an analogous formulation as (\ref{eq:hatM=M+Z_1+Z_2}). However, in this case, the second term $\tilde Z_2$ can dominate $\tilde Z_1$ and the entries of $\tilde Z_2$ are not independent. The distribution of $\tilde M$'s singular subspace is generally more difficult to analyze. Moreover, the spectral norm $\|\tilde Z_2\|$ tends to increase if the signal strength $\lambda_r$ is larger. Therefore, the naive unbiased estimator $\tilde M$ is not a suitable choice. 
\end{remark}

\section{Theory: normal approximation for singular subspaces }\label{sec:theory}
To characterize the empirical singular vectors $\hat{U}$ and $\hat{V}$, we assume that ${X}$ is a standard Gaussian matrix such that its each entry has the standard normal distribution. %The following quantity plays the critical role in the theoretic analysis of performance of $\hat{U}$ and $\hat{V}$. 
%\begin{definition}
%Let $\beta=\frac{\sigma_{\xi}}{\lambda_r}$ denote the signal-to-noise ratio (SNR) of model (\ref{eq:tr_reg}).
%\end{definition}
Even though $\hat{U}$ and $\hat{V}$ are computed from a two-step estimator (\ref{eq:nuclear_pen_est}) and (\ref{eq:debias_M}), it suffices to focus on analyzing the spectral properties of $\hat{M}$. To this end, 
the following proposition is needed which is due to \cite{candes2011tight} and \cite{negahban2011estimation}.
\begin{proposition}(\cite[Corollary~5]{negahban2011estimation} and \cite[Theorem~2.7]{candes2011tight})\label{prop:nuclear_pen}
If $n\geq C_1r\bar{m}$ and $\lambda=C_2\sigma_{\xi}\Big(\frac{\bar{m}}{n}\Big)^{1/2}$ for some universal constants $C_1,C_2>0$, then with probability at least $1-c_1\exp(-c_2\bar{m})$, 
\begin{equation}\label{eq:bDelta-frob}
\|\bDelta\|_{\rm F}^2\leq C_3\sigma_\xi^2\cdot\frac{r(m_1+m_2)}{n}
\end{equation}
for some absolute constants $c_1,c_2,C_3>0$. 
\end{proposition}
We apply the dilation operator to turn asymmetric matrices into symmetric ones. It is a standard technique to treat singular subspaces. 
See \cite{koltchinskii2015perturbation}, \cite{tropp2012user} and \cite{xia2019data} for more details. For any matrix ${A}\in\RR^{m_1\times m_2}$, we define
$$
\mathfrak{D}({A})=
\left(
\begin{array}{cc}
{0}&{A}\\
{A}^{\top}&{0}
\end{array}
\right)
\in \RR^{(m_1+m_2)\times (m_1+m_2)}
$$
which is a symmetric matrix. Then, we write $\hat{N}=\mathfrak{D}(\hat{M})$ and ${N}=\mathfrak{D}({M})$ where
$$
\hat{N}={N}+{E}:={N}+{E}_1+{E}_2
$$
with ${E}_1=\mathfrak{D}({Z}_1)$ and ${E}_2=\mathfrak{D}({Z}_2)$ where $Z_1,Z_2$ are defined in (\ref{eq:hatM=M+Z_1+Z_2}).
\begin{lemma}\label{lem:E_op}
Suppose that $n\geq \log\bar{m}$. 
There exist absolute constants $C_1,C_2>0$ such that
\begin{align*}
\EE\|{E}_1\|\leq C_1\sigma_{\xi}\frac{\bar{m}^{1/2}}{n^{1/2}}\quad {\rm and}\quad 
\EE\|{E}_2\|\leq C_2\|\bDelta\|_{\rm F}\frac{\bar{m}^{1/2}\log^{1/2}\bar{m}}{n^{1/2}}.
\end{align*}
There exist absolute constants $C_3,C_4>0$ such that for all $t\geq 1$, the following bound holds with probability at least $1-3e^{-t}-e^{-n}$,
\begin{align*}
\big|\|{E}_1\|-\EE\|{E}_1\|\big|\leq&C_3\sigma_\xi\cdot\bigg[\frac{t^{1/2}}{n^{1/2}}+\frac{\bar{m}^{1/2}t^{1/2}}{n}\bigg]\\
\big|\|{E}_2\|-\EE\|{E}_2\|\big|\leq& C_4\|\bDelta\|_{\rm F}\cdot\bigg[\frac{t^{1/2}+\log^{1/2}\bar{m}}{n^{1/2}}+\frac{\bar{m}^{1/2}t^{1/2}+t}{n}\bigg].
\end{align*}
\end{lemma}

\subsection{Representation of empirical singular vectors}
We write ${M}={U}{\Lambda}{V}^{\top}$ where ${U}=({u}_1,\ldots,{u}_r)$ and ${V}=({v}_1,\ldots,{v}_r)$ and ${\Lambda}={\rm diag}(\lambda_1,\ldots,\lambda_r)$. It is easy to check that the matrix ${N}$ has $2r$ non-zero eigenvalues which are $\lambda_1\geq\ldots\geq \lambda_r>0\geq \lambda_{-r}\geq\ldots\geq \lambda_{-1}$ where $\lambda_{-k}=-\lambda_{k}$ for $1\leq k\leq r$. The eigenvectors (which might not be unique) corresponding to the eigenvalue $\lambda_k$ and $\lambda_{-k}$ can be written, respectively, as
$$
{\theta}_k=\frac{1}{\sqrt{2}}\left(
\begin{array}{c}
{u}_k^{\top},{v}_k^{\top}
\end{array}
\right)^{\top}\quad{\rm and}\quad
{\theta}_{-k}=\frac{1}{\sqrt{2}}\left(
\begin{array}{c}
{u}_k^{\top},-{v}_k^{\top}
\end{array}
\right)^{\top}.
$$
The spectral projector corresponding to ${N}$ is defined as
$$
\calP_{{U}{V}}=\sum_{1\leq |k|\leq r}{\theta}_k{\theta}_k^{\top}=
\left(
\begin{array}{cc}
{U}{U}^{\top}&{0}\\
{0}&{{V}{V}^{\top}}
\end{array}
\right).
$$
Let $\{\hat{{\theta}}_k\}_k$ and $\{\hat{\theta}_{-k}\}_k$ represent the eigenvectors of $\hat{N}$ corresponding to the $r$ largest and $r$ smallest eigenvalues of $\hat{N}$. Then, we define the empirical spectral projector
$$
\calP_{\hat{U}\hat{V}}=\sum_{1\leq |k|\leq r}\hat{\theta}_k\hat{\theta}_k^{\top}=
\left(
\begin{array}{cc}
\hat{U}\hat{U}^{\top}&{0}\\
{0}&{\hat{V}\hat{V}^{\top}}
\end{array}
\right).
$$
By the definition of ${\rm dist}^2\big[(\hat{U},\hat{V}), ({U},{V})\big]$ in eq. (\ref{eq:def_dist}), we can immediately write
$$
{\rm dist}^2\big[(\hat{U},\hat{V}), ({U},{V})\big]=\|\calP_{{U}{V}}-\calP_{\hat{U}\hat{V}}\|_{\rm F}^2.
$$
We write $\calP_{{U}{V}}^{\perp}$ as the orthogonal projection onto the complement of the image space of $\calP_{{U}{V}}$. More explicitly, we can write
$$
\calP_{{U}{V}}^{\perp}
=\left(
\begin{array}{cc}
{U}_{\perp}{U}_{\perp}^{\top}&0\\
0&{V}_{\perp}{V}_{\perp}^{\top}
\end{array}
\right)
$$
where ${U}_{\perp}$ and ${V}_{\perp}$ are chosen such that $({U},{U}_{\perp})$ and $({V},{V}_{\perp})$ are both orthogonal matrices.
Another important operator is
\begin{gather*}
\calC_{{U}{V}}=\sum_{1\leq |k|\leq r}\frac{1}{\lambda_k}({\theta}_k{\theta}_k^{\top})
=\left(
\begin{array}{cc}
0&{U}{\Lambda}^{-1}{V}^{\top}\\
{V}{\Lambda}^{-1}{U}^{\top}&0
\end{array}
\right).
\end{gather*}
%For any $\alpha>0$, define the event
%\begin{equation}\label{eq:calE0_alpha}
%\calE_0(\alpha)=\big\{\lambda_r\geq 2(1+\alpha)\|{E}\|\big\}.
%\end{equation}
%In view of Lemma~\ref{lem:E_op}, $\PP(\calE_0(\alpha))\geq 1-e^{-n}-e^{-t}$ for any $t\geq 1$ if 
%$$
%n\geq C_1(1+\alpha)^2\beta^2(\bar{m}+t)+C_2(1+\alpha)^2\frac{\|\bDelta\|_{\rm F}^2}{\lambda_r^2}\big(\bar{m}\log{\bar m}+t\big)
%$$
%where $\beta=\frac{\sigma_{\xi}}{\lambda_r}$  and $C_1,C_2$ are absolute constants.

\begin{lemma}
\label{lem:hatcalP_N}
The following  decomposition of $\calP_{\hat{U}\hat{V}}$ holds
%\begin{align}\label{eq:hatcalP_N-calP_N_op}
%\|\calP_{\hat{U}\hat{V}}-\calP_{{U}{V}}\|\leq \frac{4\|{E}\|}{\lambda_r}
%\end{align}
%and
\begin{eqnarray*} 
\calP_{\hat{U}\hat{V}}-\calP_{{U}{V}}=\calL_{N}({E})+\calS_{N}({E}),
\end{eqnarray*}
where $\calL_{N}({E}):=\calP_{{U}{V}}^{\perp}{E} \calC_{{U}{V}}+\calC_{{U}{V}} {E} \calP_{{U}{V}}^{\perp}$ and 
\begin{eqnarray*}
\|\calL_{N}({E})\|\leq \frac{2\|{E}\|}{\lambda_r}\quad{\rm and}\quad\|\calS_{N}({E})\|\leq 80\cdot\Big(\frac{\|{E}\|}{\lambda_r}\Big)^2.
\end{eqnarray*}
\end{lemma}
\begin{remark}
The representation formula of the joint spectral projectors $\hat U\hat U^{\top}$ and $\hat V\hat V^{\top}$ in Lemma~\ref{lem:hatcalP_N} is interesting because there will be no eigen-gap requirements on the distinct singular values $\lambda_1,\cdots, \lambda_r$ in the subsequent results in the next sections. If we directly apply the existing methods and results in the literature (\cite{koltchinskii2015perturbation}, \cite{koltchinskii2015normal}), then we require that the eigen-gaps: $\min_{1\leq i\leq r-1}|\lambda_{i}-\lambda_{i+1}|\gg \|E\|$ and $\lambda_r\gg \|E\|$ which is unnecessary in Lemma~\ref{lem:hatcalP_N}.  In other words, Lemma~\ref{lem:hatcalP_N} allows the singular values to have multiplicity larger than $1$.
\end{remark}

\subsection{Normal approximation of projection distance}
In this section, we will prove the normal approximation of the loss ${\rm dist}^2[(\hat U,\hat V),(U,V)]$. By Lemma~\ref{lem:hatcalP_N}, we immediately have 
\begin{align*}
{\rm dist}^2[(\hat U,\hat V),(U,V)]=&\|\calP_{\hat U\hat V}-\calP_{UV}\|_{\rm F}^2\\
=&\|\calL_N(E)\|_{\rm F}^2+\|\calS_N(E)\|_{\rm F}^2+2\big<\calL_N(E),\calS_N(E)\big>.
\end{align*}
We begin with the linear term $\calL_{N}({E})=\calP_{UV}^{\perp}E\calC_{UV}+\calC_{UV}E\calP_{UV}^{\perp}$. In particular, the variance of ${\rm dist}^2[(\hat U,\hat V),(U,V)]$ can be  characterized by the variance of $\|\calL_N(E)\|_{\rm F}^2$. 
\begin{theorem}\label{thm:con_LNE}
Denote by $\beta=\frac{\sigma_\xi}{\lambda_r}$.
Suppose that $n\geq C\big(\beta^2\bar{m}+r\bar{m}\log^2n\big)$ and $n\leq C^{-1}e^{\bar{m}}$ for a large enough constant $C>0$. Then,
there exist absolute constants $c_1,c_2,C_5,C_6>0$ such that with probability at least $1-\frac{2n+5}{n^2}-2e^{-n}-c_1ne^{-c_2\bar{m}}$, 
\begin{align*}
\big|\|\calL_{N}({E})\|_{\rm F}^2&-\EE\|\calL_{N}({E})\|_{\rm F}^2\big|
\leq C_5\sigma_\xi^2\|{\Lambda}^{-2}\|_{\rm F}\cdot\frac{\bar{m}^{1/2}\log^{1/2}n}{n}
+C_6\sigma_\xi^2\|{\Lambda}^{-1}\|_{\rm F}^2\cdot\frac{r\bar{m}\log n}{n^{3/2}}
\end{align*}
 and
\begin{align*}
\EE\|\calL_{N}({E})\|_{\rm F}^2=&\sigma_\xi^2\|{\Lambda}^{-1}\|_{\rm F}^2\cdot \frac{2m_{\star}}{n}
+O\Big(\sigma_\xi^2\|{\Lambda}^{-1}\|_{\rm F}^2\cdot \frac{2r\bar{m}^2}{n^2}\Big).
\end{align*}
\end{theorem}

In Theorem~\ref{thm:hatPN-PN-frob-con} and Theorem~\ref{thm:hatP-normal_approx}, we prove the concentration of the loss $\|\calP_{\hat U\hat V}-\calP_{UV}\|_{\rm F}^2$ and its related normal approximation. Eq. (\ref{eq:E_hatP-P}) implies that the dominating term of the expected loss $\EE\|\calP_{\hat U\hat V}-\calP_{UV}\|_{\rm F}^2$ is determined by $\EE\|\calL_N(E)\|_{\rm F}^2$ (see Theorem~\ref{thm:con_LNE}). Similarly, the dominating term in the variance of $\|\calP_{\hat U\hat V}-\calP_{UV}\|_{\rm F}^2$ is also determined by the variance of $\|\calL_N(E)\|_{\rm F}^2$. 

\begin{theorem}\label{thm:hatPN-PN-frob-con}
Denote by $\beta=\frac{\sigma_\xi}{\lambda_r}$. Suppose that $n\geq C_6\big(\beta^2\bar{m}+r\bar{m}\log^2n\big)$ and $n\leq C_6^{-1}e^{\bar{m}}$ for some large enough absolute constant $C_6>0$. Then, there exist absolute constants $c_1,c_2,C_7,C_8>0$ such that with probability at least $1-\frac{2n+9}{n^2}-3e^{-n}-c_1ne^{-c_2\bar{m}}$,
\begin{align*}
\big|\|\calP_{\hat{U}\hat{V}}-&\calP_{{U}{V}}\|_{\rm F}^2-\EE\|\calP_{\hat{U}\hat{V}}-\calP_{{U}{V}}\|_{\rm F}^2\big|\\
\leq&  C_7\sigma_\xi^2\|{\Lambda}^{-2}\|_{\rm F}\cdot\frac{\bar{m}^{1/2}\log^{1/2}n}{n}+C_8\Big[\Big(\frac{\sigma_\xi}{\lambda_r}\Big)^3+\sigma_\xi^2\|\Lambda^{-1}\|_{\rm F}^2\Big]\cdot\frac{r\bar{m}\log^{1/2}n}{n^{3/2}}.
\end{align*}
%If $\frac{r^2\log n+\log^3n}{n}\to 0$, then,
%$$
%\|\calP_{\hat{U}\hat{V}}-\calP_{{U}{V}}\|_{\rm F}^2\to \EE\|\calP_{\hat{U}\hat{V}}-\calP_{{U}{V}}\|_{\rm F}^2,\quad \textrm{ in probability}.
%$$
and
\begin{align}
\EE\|\calP_{\hat{U}\hat{V}}-&\calP_{{U}{V}}\|_{\rm F}^2=\sigma_{\xi}^2\|{\Lambda}^{-1}\|_{\rm F}^2\cdot\frac{2m_{\star}}{n}+O\Big((\beta\vee 1)^4\cdot \frac{r^{3/2}\bar{m}^2\log^{1/2}n}{n^2}\Big).\label{eq:E_hatP-P}
\end{align}
\end{theorem}
\begin{remark}
The most important conclusion in Theorem~\ref{thm:hatPN-PN-frob-con}  is that the second order term in $\EE\|\calP_{\hat U\hat V}-\calP_{UV}\|_{\rm F}^2$ is of the order $\frac{\bar{m}^2}{n^2}$ rather than the order $\frac{\bar{m}^{3/2}}{n^{3/2}}$ (if by a naive analysis). This improvement comes from the second order analysis on the perturbation formula of the empirical spectral projectors. Basically, we treat $\EE\langle\calL_N(E), \calS_N(E) \rangle$ more sophisticatedly and will prove that the term involving $\frac{\bar{m}^{3/2}}{n^{3/2}}$ vanishes. 
\end{remark}

\begin{theorem}\label{thm:hatP-normal_approx}
Suppose the conditions in Theorem~\ref{thm:hatPN-PN-frob-con} hold and $n\geq C_1r^2\bar{m}$  for a large enough absolute constant $C_1>0$. Let $\Phi(\cdot)$ denote the {\it cumulative distribution function} of the standard normal distribution. Then,
\begin{align*}
\sup_{x\in\RR}\bigg|\PP\Big\{&\frac{\|\calP_{\hat{U}\hat{V}}-\calP_{{U}{V}}\|_{\rm F}^2-\EE\|\calP_{\hat{U}\hat{V}}-\calP_{{U}{V}}\|_{\rm F}^2}{\sqrt{8}\sigma_{\xi}^2\|{\Lambda}^{-2}\|_{\rm F}\cdot \frac{m_{\star}^{1/2}}{n}}\leq x\Big\}-\Phi(x) \bigg|\\
\leq&C_7(\beta\vee 1)\frac{r\bar{m}^{1/2}\log^{3/2}n}{n^{1/2}}+c_1ne^{-c_2m_{\star}}+\frac{2n+7}{n^2}+ \frac{C_8}{\bar{m}^{1/2}}.
\end{align*}
for absolute constants $c_1,c_2,C_7,C_8>0$.
\end{theorem}
\begin{remark}\label{rmk:conf_reg}
Theorem~\ref{thm:hatP-normal_approx} implies that if $\bar{m},n\to\infty$ and $\frac{r^2\bar{m}\log^3 n}{n}\to 0$ when $\beta=O(1)$, then 
$$
\frac{\|\calP_{\hat{U}\hat{V}}-\calP_{{U}{V}}\|_{\rm F}^2-\EE\|\calP_{\hat{U}\hat{V}}-\calP_{{U}{V}}\|_{\rm F}^2}{\sqrt{8}\sigma_{\xi}^2\|{\Lambda}^{-2}\|_{\rm F}\cdot \frac{m_{\star}^{1/2}}{n}}\stackrel{{\rm d}}{\longrightarrow}\calN(0,1).
$$
The sample size requirement $n\gg r^2\bar{m}\log^3n$ is optimal up to the rank $r$ and the logarithmic factor. It also implies that the ``ideal" $100(1-\alpha)\%$  confidence region of $(\hat U,\hat V)$ is 
\begin{align*}
\calC_{\alpha}:=\Big\{(X,Y&): X\in\OO^{m_1\times r}, Y\in\OO^{m_2\times r},\\
& \big|{\rm dist}^2[(X,Y),(\hat U,\hat V)]-\EE\|\calP_{\hat U\hat V}-\calP_{UV}\|_{\rm F}^2\big|\leq \sqrt{8}z_{\alpha/2}\sigma_\xi^2\|\Lambda^{-2}\|_{\rm F}\cdot \frac{m_{\star}^{1/2}}{n}\Big\}
\end{align*}
where $z_{\alpha}=\Phi^{-1}(1-\alpha)$. 
It is ``ideal" because the centering term $\EE\|\calP_{\hat U\hat V}-\calP_{UV}\|_{\rm F}^2$ is not completely determined yet. Clearly by the differential property on Grassmannians (see \cite{xia2019tensor} and \cite{edelman1998geometry}), the diameter of the ``ideal" confidence region $\calC_{\alpha}$ in the (squared) projection distance has the same order as $\EE\|\calP_{\hat U\hat V}-\calP_{UV}\|_{\rm F}^2$. By eq. (\ref{eq:E_hatP-P}), we can conclude that the diameter of the ``ideal" confidence region has the order $\sigma_\xi^2\|\Lambda^{-1}\|_{\rm F}^2\cdot \frac{m_{\star}}{n}$ as long as $n\gg r^2\bar{m}\log^3n$.  Because the stochastic deviation $\sigma_\xi^2\|\Lambda^{-2}\|_{\rm F}\cdot m_{\star}^{1/2}/n$ is much smaller than the bias $\EE\|\calP_{\hat U\hat V}-\calP_{UV}\|_{\rm F}^2$,  the ``ideal" confidence region also implies that the minimax optimal diameter of the confidence regions for $(U,V)$ has the order $\sigma_\xi^2\|\Lambda^{-1}\|_{\rm F}^2\cdot \frac{m_{\star}}{n}$. 
\end{remark}

By Remark~\ref{rmk:conf_reg}, the confidence region of the true singular subspace can be constructed if we can completely determine the expected loss $\EE\|\calP_{\hat U\hat V}-\calP_{UV}\|_{\rm F}^2$. 
Now, we replace $\EE\|\calP_{\hat{U}\hat{V}}-\calP_{{U}{V}}\|_{\rm F}^2$ with its first order approximation $\sigma_\xi^2\|{\Lambda}^{-1}\|_{\rm F}^2\cdot \frac{2m_{\star}}{n}$ from Theorem~\ref{thm:hatPN-PN-frob-con} and obtain the following normal approximation of ${\rm dist}^2[(\hat U,\hat V),(U,V)]$ with the explicit centering and normalization terms. By using only the first order approximation of $\EE\|\calP_{\hat{U}\hat{V}}-\calP_{{U}{V}}\|_{\rm F}^2$, we need a larger sample size requirement for the asymptotical normality (compared with Theorem~\ref{thm:hatP-normal_approx}). 
\begin{corollary}\label{cor:hatP-normal_approx}
Suppose the conditions in Theorem~\ref{thm:hatPN-PN-frob-con} hold and $n\geq C_1r^2\bar{m}$ for a large enough absolute constant $C_1>0$. Let $\Phi(\cdot)$ denote the {\it cumulative distribution function} of the standard normal distribution. Then, 
\begin{align*}
\sup_{x\in\RR}\bigg|\PP\Big\{&\frac{\|\calP_{\hat{U}\hat{V}}-\calP_{{U}{V}}\|_{\rm F}^2-\sigma_\xi^2\|{\Lambda}^{-1}\|_{\rm F}^2\cdot \frac{2m_{\star}}{n}}{\sqrt{8}\sigma_{\xi}^2\|{\Lambda}^{-2}\|_{\rm F}\cdot \frac{m_{\star}^{1/2}}{n}}\leq x\Big\}-\Phi(x) \bigg|\\
\leq&C_7(\beta\vee 1)\frac{r\bar{m}^{1/2}\log^{3/2}n}{n^{1/2}}+C_8(\beta\vee 1)^2\frac{r^{3/2}\bar{m}^{3/2}\log^{1/2}n}{n}+c_1ne^{-c_2m_{\star}}+\frac{3n+6}{n^2}+ \frac{C_{9}}{\bar{m}^{1/2}}.
\end{align*}
for absolute constants $c_1,c_2,C_7,C_8,C_9>0$.
\end{corollary}
\begin{remark}
Corollary~\ref{cor:hatP-normal_approx} implies that if $\bar{m},n\to\infty$ and $\frac{r^{3/2}\bar{m}^{3/2}\log^{1/2} n}{n}\to 0$ when $\beta=O(1)$, then 
$$
\frac{\|\calP_{\hat{U}\hat{V}}-\calP_{{U}{V}}\|_{\rm F}^2-\sigma_\xi^2\|{\Lambda}^{-1}\|_{\rm F}^2\cdot \frac{2m_{\star}}{n}}{\sqrt{8}\sigma_{\xi}^2\|{\Lambda}^{-2}\|_{\rm F}\cdot \frac{m_{\star}^{1/2}}{n}}\stackrel{{\rm d}}{\longrightarrow}\calN(0,1).
$$
We note that the sample size requirement $n\gg \bar{m}^{3/2}$ is optimal for the above normal approximation. The reason is that the approximation error by Lemma~\ref{lemma:PN-E2-CN-frob}, conditioned on $\Delta$, is 
$$
\Big|\EE\|\calP_{\hat U\hat V}-\calP_{UV}\|_{\rm F}^2- 2\sigma_\xi^2\|\Lambda^{-1}\|_{\rm F}^2\cdot m_{\star}/n\Big|\geq c_1\frac{m_{\star}}{n}\cdot \|\Lambda^{-1}\|_{\rm F}^2\|\Delta\|_{\rm F}^2
$$
for some absolute constant $c_1>0$. By the minimax optimal lower bounds of low-rank matrix regression (\cite{candes2011tight} and \cite{koltchinskii2011nuclear}), $\|\Delta\|_{\rm F}^2$ is lower bounded by the rate $\sigma_\xi^2rm_{\star}/n$ with probability at least $c_2$ for some constant $c_2>0$. Together with Theorem~\ref{thm:hatP-normal_approx},  it is easy to check that the asymptotical normality in Corollary~\ref{cor:hatP-normal_approx} holds only when $n\gg \bar{m}^{3/2}$. 
\end{remark}

\begin{remark}
Let's compare with the PCA results in \cite{koltchinskii2017new} where the limiting distribution is a Cauchy distribution. In \cite{koltchinskii2017new}, a data-dependent estimator of the expected loss is designed whose error follows a Gaussian distribution with the standard deviation comparable with the normalization term (that is $\frac{\bar{m}^{1/2}}{n}$ in our problem).  As a result, they end up with a Cauchy distribution. However, our estimation error of $\EE\|\calP_{\hat U\hat V}-\calP_{UV}\|_{\rm F}^2$ by Theorem~\ref{thm:hatPN-PN-frob-con} is of the order $\frac{\bar{m}^2}{n^2}$. Therefore, if $n\gg \bar{m}^{3/2}$, the limiting distribution we get is a Gaussian distribution. 
\end{remark}

\section{Data-dependent confidence regions of singular subspaces}\label{sec:CI}
In this section, we apply the limiting distributions established in Theorem~\ref{thm:hatP-normal_approx} to construct the confidence regions of ${U}$ and ${V}$. We assume that the true rank $r$ is known. In Section~\ref{sec:discussion}, we will discuss about a simple method to estimate the true rank if $r$ is not given in advance. 
In view of Theorem~\ref{thm:hatP-normal_approx}, it suffices to estimate $\sigma_{\xi}^2$, $\|{\Lambda}^{-1}\|_{\rm F}^2$ and $\|{\Lambda}^{-2}\|_{\rm F}$. Recall the definition of $\hat{M}^{\rm nuc}$, we estimate the noise variance by $\hat{M}^{\rm nuc}$'s goodness of fitting data $\{({X}_i, y_i)\}_{i=n+1}^{2n}$. More exactly, we define
\begin{equation}\label{eq:hat_sigma}
\hat\sigma_\xi^2:= \frac{1}{n}\sum_{i=n+1}^{2n}\big(y_i-\tr({X}_i^{\top}\hat{M}^{\rm nuc})\big)^2.
\end{equation}
Recall that the singular values of $\hat{M}$ are denoted by $\hat\lambda_k$. 
To this end,  we define
\begin{equation}\label{eq:bias_Bn}
\hat{B}_n:=\sum_{k=1}^r \tilde\lambda_k^{-2}
\end{equation}
, where $\tilde\lambda_k^2:=\hat\lambda_k^2-\frac{2m_{\star}}{n}\cdot \hat\sigma_\xi^2$. The shrinkage estimators $\{\tilde\lambda_k\}_{k\geq 1}^r$ are inspired by random matrix theory (\cite{ding2017high}). 
 Similarly, we define the estimator of $\|{\Lambda}^{-2}\|_{\rm F}^2$ as
\begin{equation}\label{eq:std_Vn}
\hat{V}_n=\sum_{k=1}^r \tilde\lambda_k^{-4}.
\end{equation}
Lemma~\ref{lem:V_n-B_n} provides the accuracy of $\hat B_n$ and $\hat V_n$. We note that sharper characterization of $\hat B_n$ and $\hat V_n$ might be possible, but the bounds in Lemma~\ref{lem:V_n-B_n} are sufficient for the objectives of this article. 
\begin{lemma}\label{lem:V_n-B_n}
Denote by $\beta=\frac{\sigma_\xi}{\lambda_r}$. 
Suppose that $n\geq Cr\bar{m}$ for a large enough constant $C>0$. Then, with probability at least $1-\frac{1}{\bar{m}^2}-c_1e^{-c_2\bar{m}}$,
\begin{align*}
\big|\hat B_n-\|\Lambda^{-1}\|_{\rm F}^2 \big|\leq C_2(\beta\vee 1)^2\|\Lambda^{-1}\|_{\rm F}^2\cdot \frac{r^{1/2}\bar{m}\log^{1/2}\bar{m}}{n}
\end{align*}
and
\begin{align*}
\big|\hat V_n-\|\Lambda^{-2}\|_{\rm F}^2 \big|\leq C_2(\beta\vee 1)^2\|\Lambda^{-2}\|_{\rm F}^2\cdot \frac{ r^{1/2}\bar{m}\log^{1/2}\bar{m}}{n}
\end{align*}
for some absolute constant $c_1,c_2,C_2>0$. In addition, with probability at least $1-\frac{2}{n^2}-c_1e^{-c_2\bar{m}}$, 
$$
\big|\hat\sigma_\xi^2-(\sigma_\xi^2+\|\bDelta\|_{\rm F}^2) \big|\leq C_7\sigma_\xi^2\cdot \frac{\log n}{n^{1/2}}
$$
for some absolute constant $c_1,c_2,C_7>0$. 
\end{lemma}
We define a new statistics:
\begin{equation*}
\hat T_{{U}{V}}:=\frac{\|\calP_{\hat{U}\hat{V}}-\calP_{{U}{V}}\|_{\rm F}^2-(\hat B_n\hat \sigma_\xi^2)\cdot \frac{2m_{\star}}{n}}{\sqrt{8}\hat V_n^{1/2}\hat \sigma_\xi^2\cdot \frac{m_{\star}^{1/2}}{n}}
\end{equation*}
 and prove the normal approximation of  $\hat T_{{U}{V}}$ in Theorem~\ref{thm:hatT_bUbV-normal-approx}. 
\begin{theorem}\label{thm:hatT_bUbV-normal-approx}
Denote by $\beta=\frac{\sigma_\xi}{\lambda_r}$. Suppose  that $n\geq C\big(\beta^2\bar{m}+r\bar{m}\log^2n\big)$ and $n\leq C^{-1}e^{\bar{m}}$ for a large enough constant $C>0$. Then, 
\begin{align*}
\sup_x\Big|\PP\big\{&\hat T_{{U}{V}}\leq x\big\}-\Phi(x)\Big|\\
\leq& C_7(\beta\vee 1)^4\cdot \bigg(\frac{r\bar{m}^{1/2}\log^{3/2}n}{n^{1/2}}+\frac{r^{3/2}\bar{m}^{3/2}\log n}{n}\bigg)+c_1ne^{-m_{\star}}+\frac{C_8}{\bar{m}^{1/2}}
\end{align*}
for absolute constants $c_1,c_2,C_7, C_8>0$.
\end{theorem}
\begin{remark}
In Theorem~\ref{thm:hatT_bUbV-normal-approx}, if $\beta=O(1)$ and
$$
\frac{r^2\bar{m}\log^3 n+r^{3/2}\bar{m}^{3/2}\log n}{n}\stackrel{\bar{m},n\to\infty}{\longrightarrow}0
$$
when $\beta=O(1)$, then $\hat T_{{U}{V}}\stackrel{{\rm d}}{\longrightarrow}\calN(0,1)$ as $\bar{m},n\to\infty$. In the case $r\ll \bar{m}$, it suffices to require the sample size $n\gg r^{3/2}\bar{m}^{3/2}\log n$.
\end{remark}
We apply the normal approximation of $\hat T_{{U}{V}}$ to construct confidence regions of ${U}$ and ${V}$. The following corollary is an immediate result from Theorem~\ref{thm:hatT_bUbV-normal-approx}.
\begin{corollary}\label{cor:CI_alpha}
Suppose the conditions of Theorem~\ref{thm:hatT_bUbV-normal-approx} hold and suppose that
$$
\lim_{\bar{m},n\to\infty}\frac{r^2\bar{m}\log^3n+r^{3/2}\bar{m}^{3/2}\log n}{n}=0.
$$
For any $\alpha\in(0,1)$, denote by $z_{\alpha}=\Phi^{-1}(1-\alpha)$. Define the confidence region
\begin{align*}
{\rm CR}_{\alpha}:=\bigg\{({X},{Y}): X\in&\OO^{m_1\times r}, Y\in\OO^{m_2\times r} \textrm{ such that }\\
,&\Big|{\rm dist}^2[(X,Y),(\hat U,\hat V)]- \hat B_n\hat \sigma_\xi^2\cdot \frac{2m_{\star}}{n}\Big|\leq \sqrt{8}z_{\alpha/2}\hat V_n^{1/2}\hat\sigma_\xi^2\cdot \frac{m_{\star}^{1/2}}{n}\bigg\}
\end{align*}
where $\hat B_n$ and $\hat V_n$ are defined as (\ref{eq:bias_Bn}) and (\ref{eq:std_Vn}). 
If $\beta=O(1)$, then, 
$$
\lim_{\bar{m},n\to\infty}\PP\Big(({U},{V})\in{\rm CR}_{\alpha}\Big)=\alpha.
$$
\end{corollary}

\section{Numerical experiments}\label{sec:simulation}
In this section, we present some numerical results. In these simulations, the underlying low-rank matrix ${M}\in\RR^{m\times m}$ has $\rank({M})=r$ and the thin singular value decomposition ${M}={U}{\Lambda}{V}^{\top}$ where $\lambda_k=2^{r-k+1}$ for $1\leq k\leq r$. The condition number of ${M}$ is $2^{r-1}$ growing fast with respect to $r$. The singular vectors ${U}$ and ${V}$ are generated from the singular subspace of Gaussian random matrices. 
The initial estimator $\hat{M}^{\rm nuc}$ is solved by the famous {\it alternating direction method of multipliers} (ADMM) algorithm. See \cite{boyd2011distributed} for more details.

\begin{figure}
\centering
\begin{subfigure}{0.49\textwidth}
\includegraphics[scale=0.45]{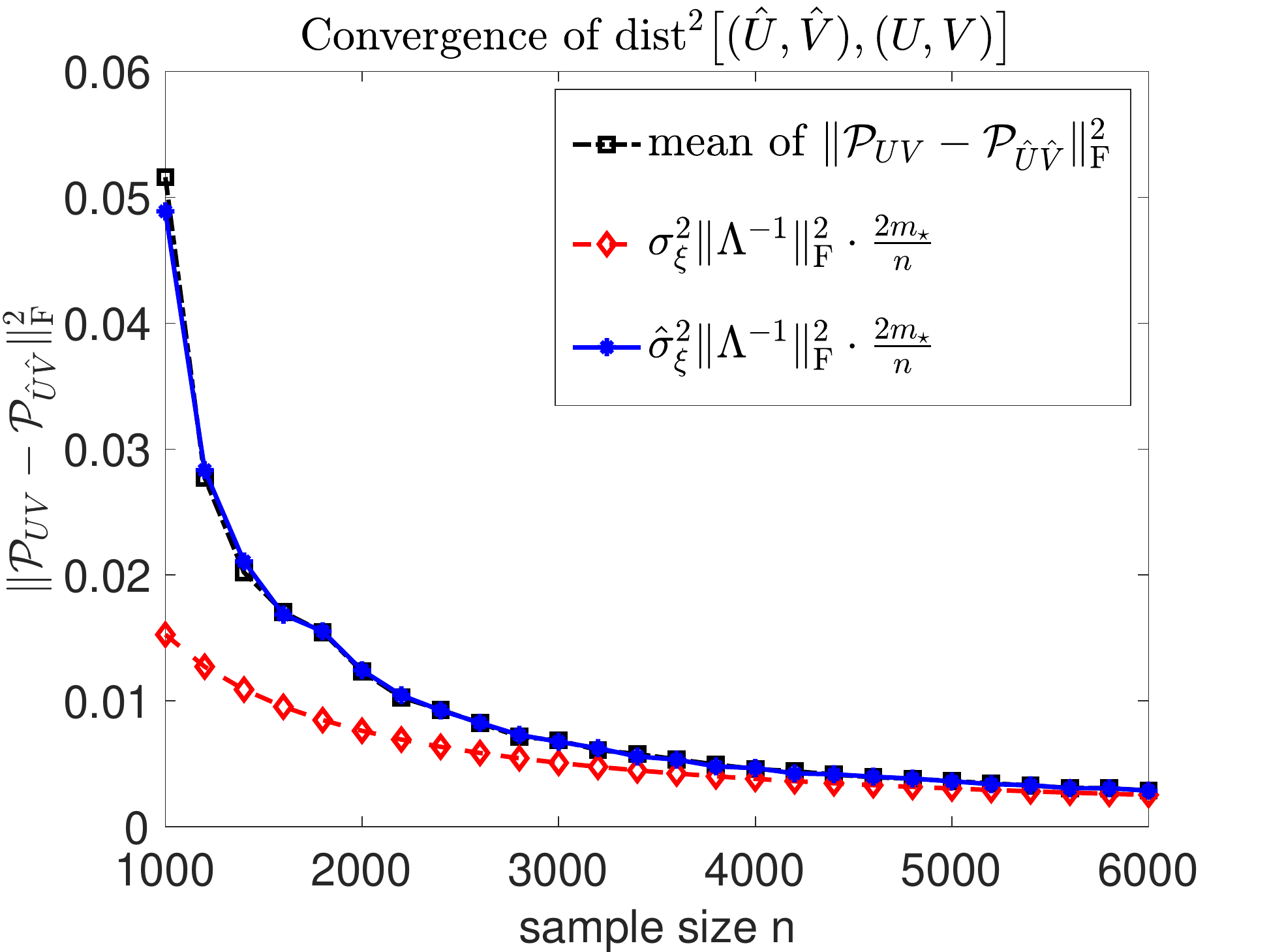}
\caption{$m=50,r=4$ and $\sigma=0.5$}
\label{fig:conv_1}
\end{subfigure}
\begin{subfigure}{0.49\textwidth}
\includegraphics[scale=0.45]{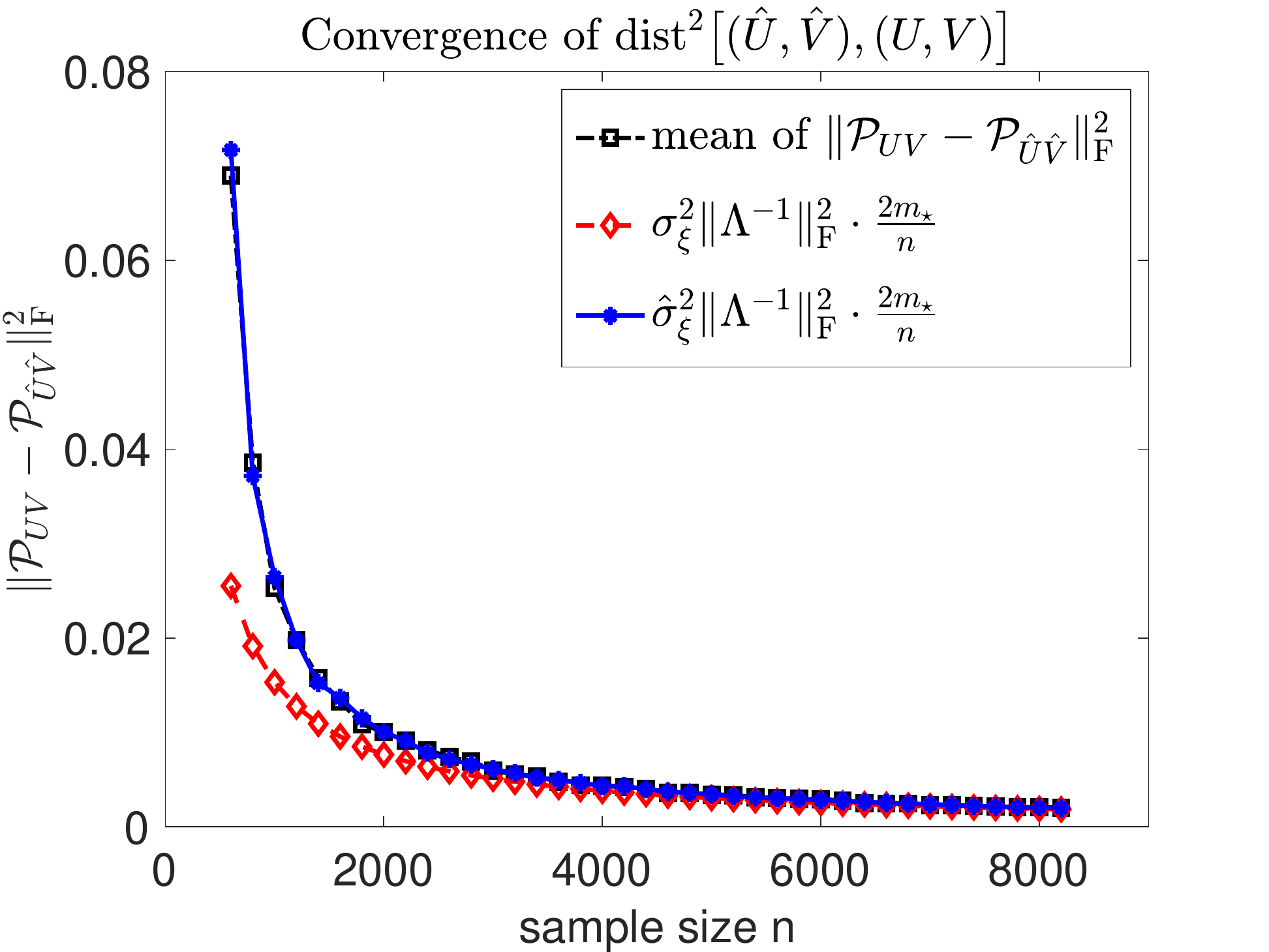}
\caption{$m=100, r=4$ and $\sigma=0.5$}
\label{fig:conv_2}
\end{subfigure}
\caption{Comparison of $\EE\|\calP_{{U}{V}}-\calP_{\hat{U}\hat{V}}\|_{\rm F}^2$ and $\sigma_\xi^2\|{\Lambda}^{-1}\|_{\rm F}^2\cdot \frac{2m_{\star}}{n}$ and $(\EE\hat\sigma_\xi^2)\|{\Lambda}^{-1}\|_{\rm F}^2\cdot \frac{2m_{\star}}{n}$ with respect to the sample size. For each $n$, the mean of $\|\calP_{{U}{V}}-\calP_{\hat{U}\hat{V}}\|_{\rm F}^2$ and $\hat\sigma_\xi^2$ are obtained by the average of $50$ independent simulations. In theory, we have $\EE\hat\sigma_\xi^2=[1+o(1)]\sigma_\xi^2$ when $n\gg r\bar{m}$. However, when the sample size $n$ is moderately large, we observe that $\EE \hat\sigma_\xi^2$ is more accurate for characterizing $\EE\|\calP_{{U}{V}}-\calP_{\hat{U}\hat{V}}\|_{\rm F}^2$.}
\label{fig:conv}
\end{figure}

First, we compare $\EE\|\calP_{{U}{V}}-\calP_{\hat{U}\hat{V}}\|_{\rm F}^2$ with $\sigma_\xi^2\|{\Lambda}^{-1}\|_{\rm F}^2\cdot \frac{2m_{\star}}{n}$ as $n$ grows. In addition, we also compare $\EE\|\calP_{{U}{V}}-\calP_{\hat{U}\hat{V}}\|_{\rm F}^2$ with $(\EE \hat\sigma_\xi^2)\|{\Lambda}^{-1}\|_{\rm F}^2\cdot \frac{2m_{\star}}{n}$ with $\hat\sigma_\xi^2$ being defined as in (\ref{eq:hat_sigma}). In theory, we have $\EE\hat\sigma_\xi^2=[1+o(1)]\cdot \sigma_\xi^2$ as long as $n\gg r\bar{m}$.  However, in simulations when $n$ is only moderately large, we observe that $(\EE \hat\sigma_\xi^2)\|{\Lambda}^{-1}\|_{\rm F}^2\cdot \frac{2m_{\star}}{n}$ is more accurate for estimating $\EE\|\calP_{{U}{V}}-\calP_{\hat{U}\hat{V}}\|_{\rm F}^2$.
Two scenarios are implemented with $m=50, r=4, \sigma_\xi=0.5$ and $m=100, r=4, \sigma_\xi=0.5$, respectively. For each $n$, the algorithm is repeated for $50$ times on independently sampled data and the average of $\|\calP_{{U}{V}}-\calP_{\hat{U}\hat{V}}\|_{\rm F}^2$ is recorded. The empirical mean of $\|\calP_{{U}{V}}-\calP_{\hat{U}\hat{V}}\|_{\rm F}^2$, the theoretical bound $\sigma_\xi^2\|{\Lambda}^{-1}\|_{\rm F}^2\cdot \frac{2m_{\star}}{n}$ and the empirical bound $(\EE \hat\sigma_\xi^2)\|{\Lambda}^{-1}\|_{\rm F}^2\cdot \frac{2m_{\star}}{n}$ are displayed in Figure~\ref{fig:conv}. 

\begin{figure}
\centering
\begin{subfigure}{0.49\textwidth}
\includegraphics[scale=0.45]{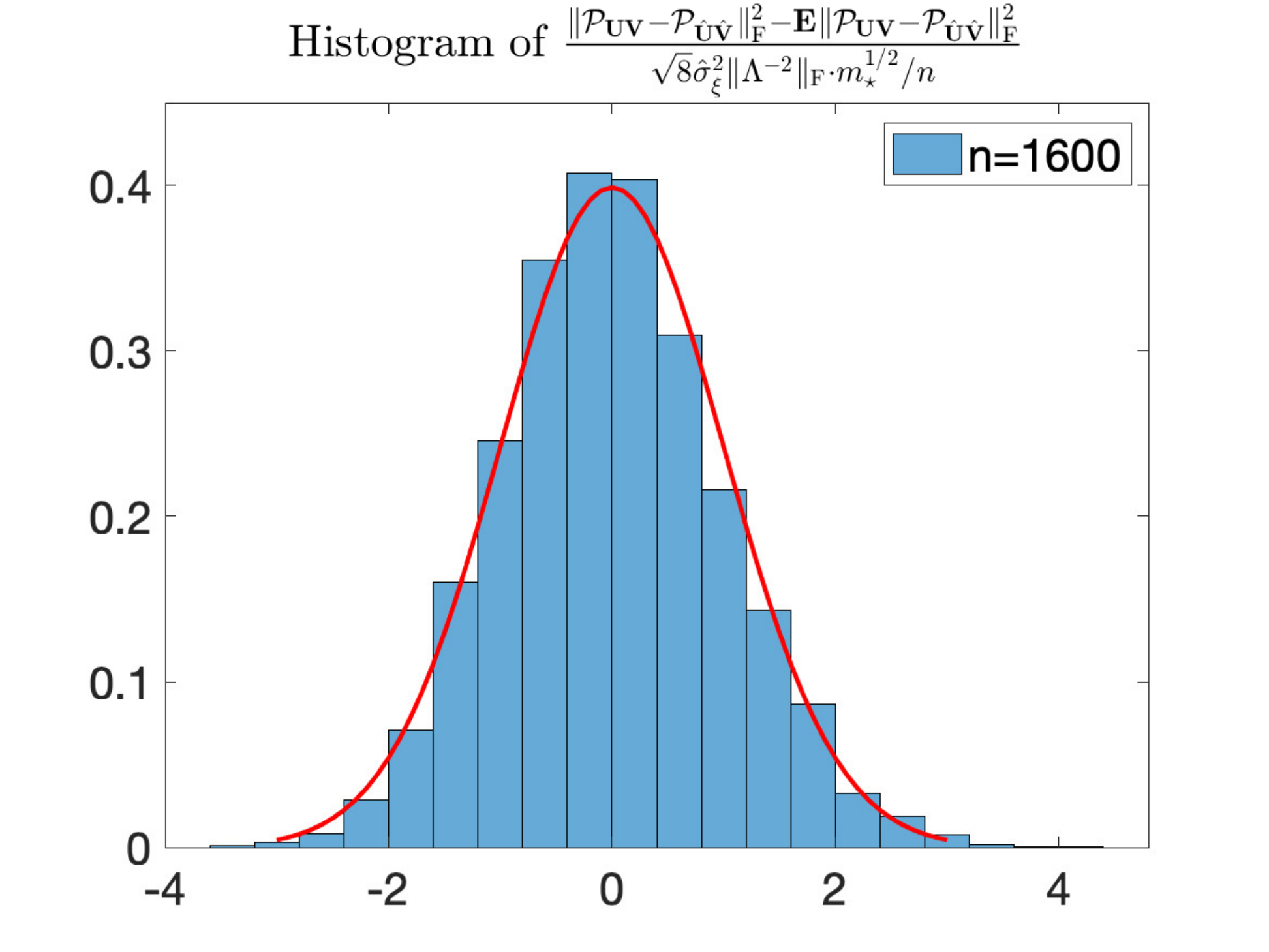}
\end{subfigure}
\begin{subfigure}{0.49\textwidth}
\includegraphics[scale=0.45]{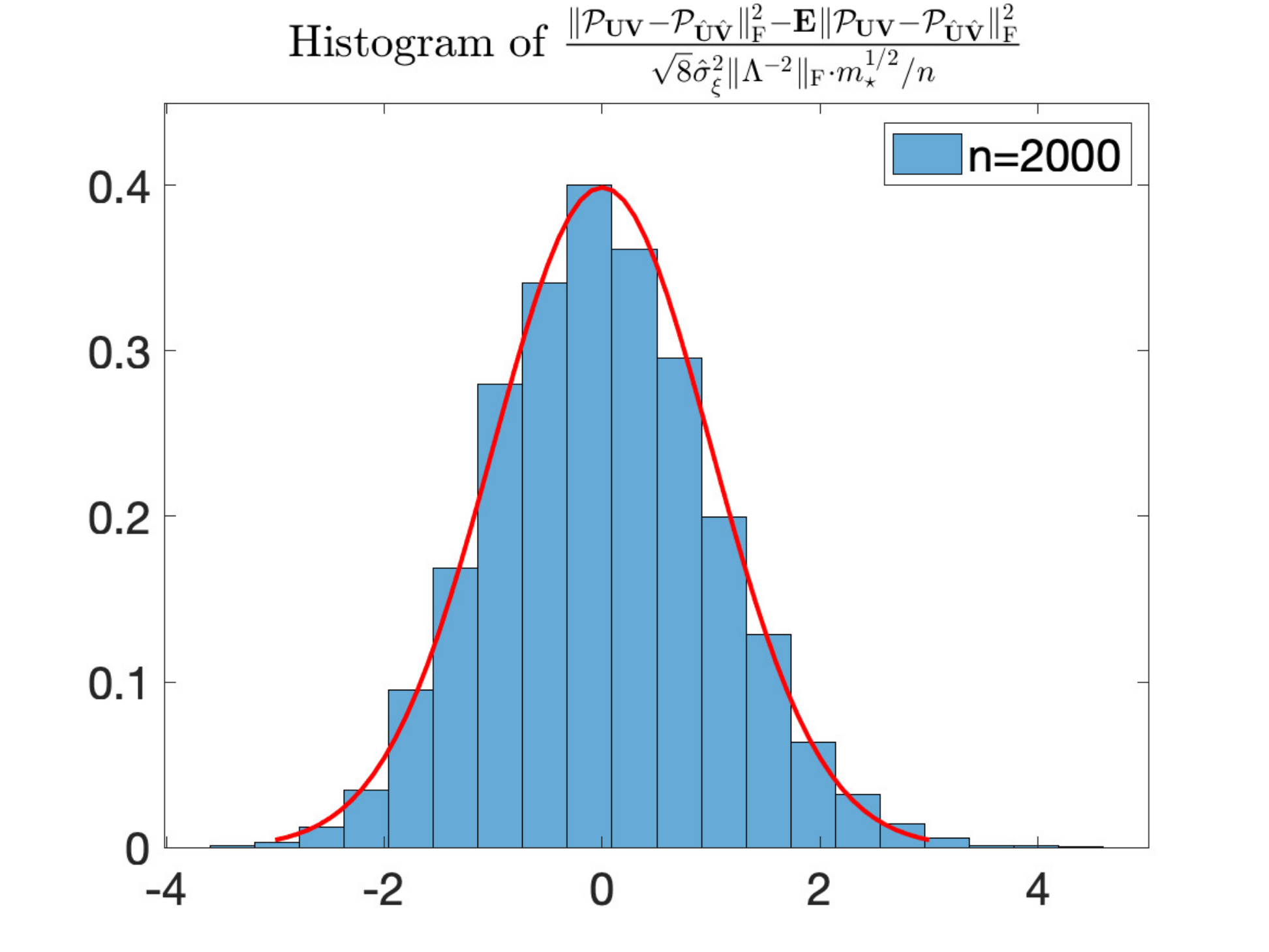}
\end{subfigure}

\begin{subfigure}{0.49\textwidth}
\includegraphics[scale=0.45]{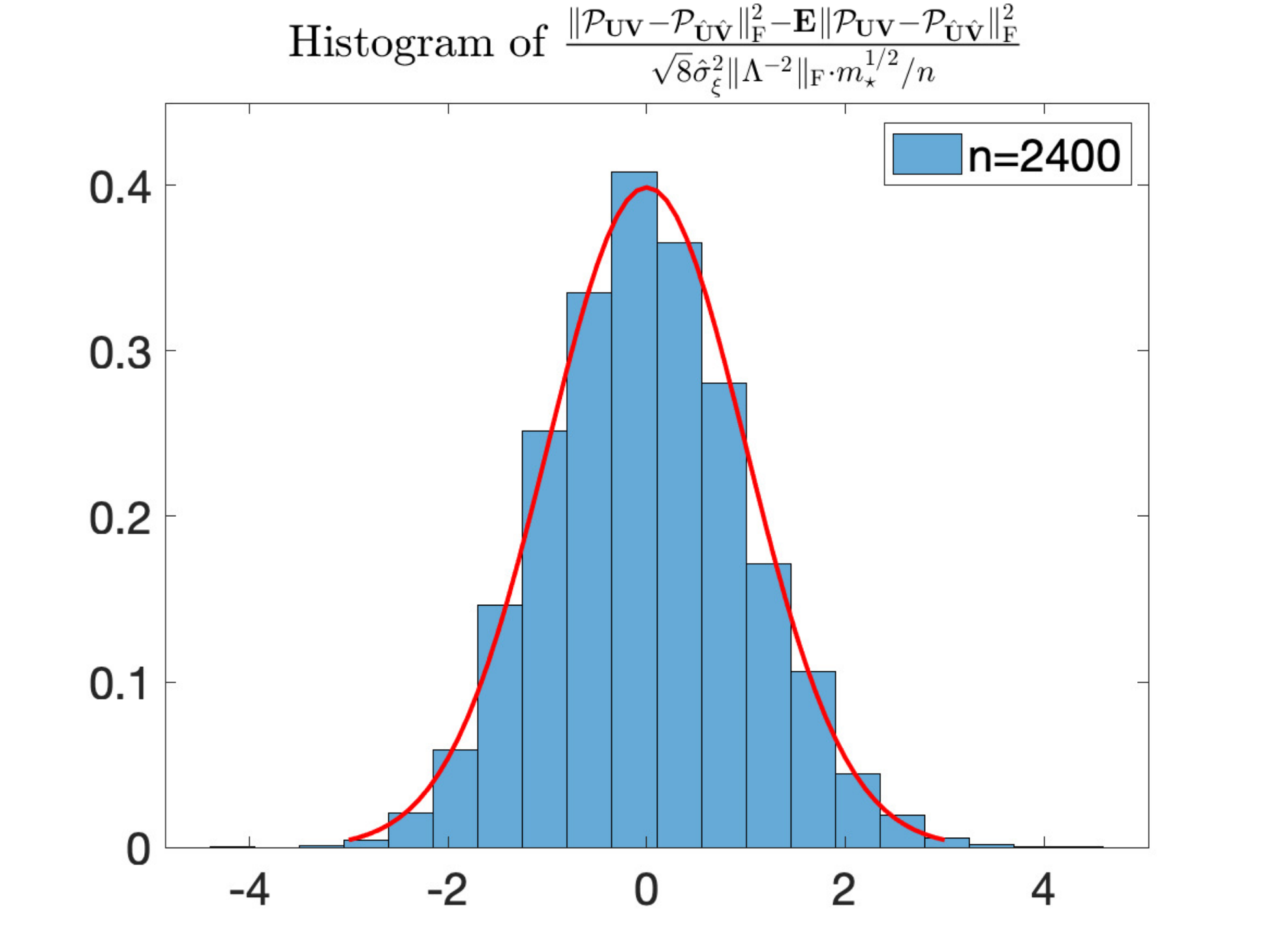}
\end{subfigure}
\begin{subfigure}{0.49\textwidth}
\includegraphics[scale=0.45]{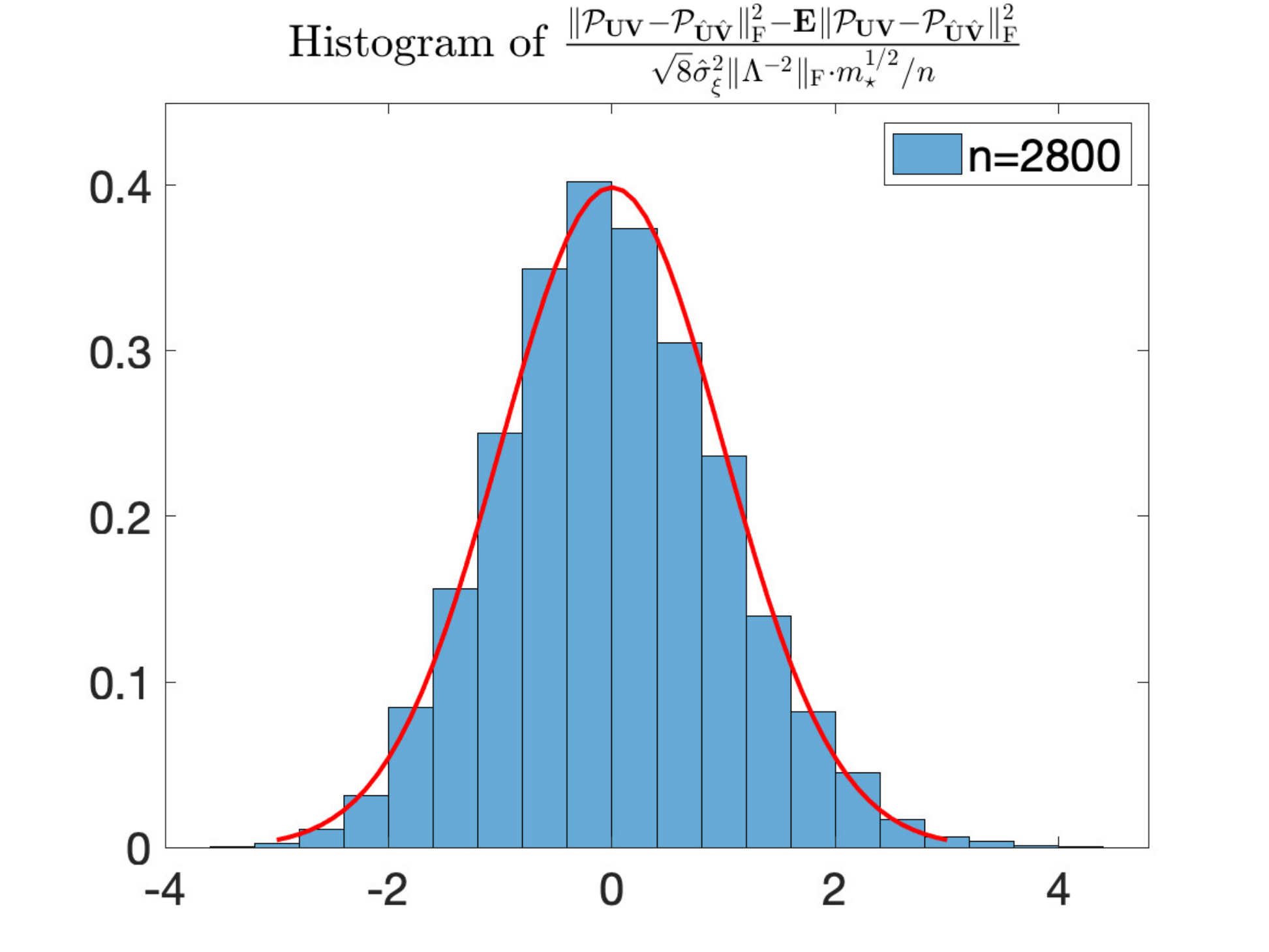}
\end{subfigure}
\caption{Normal approximation of $\frac{\|\calP_{{U}{V}}-\calP_{\hat{U}\hat{V}}\|_{\rm F}^2-\EE\|\calP_{{U}{V}}-\calP_{\hat{U}\hat{V}}\|_{\rm F}^2}{\sqrt{8}\hat\sigma_\xi^2\|{\Lambda}^{-2}\|_{\rm F}\cdot m_{\star}^{1/2}/n}$ with $m_1=m_2=100, r=4$ and $\sigma_\xi=0.1$. For each $n$, the density histogram is based on $10000$ repetitions whose average is used to estimate $\EE\|\calP_{{U}{V}}-\calP_{\hat{U}\hat{V}}\|_{\rm F}^2$. The empirical noise variance $\hat\sigma_\xi^2$ is calculated as in (\ref{eq:hat_sigma}). 
The red curve represents the {\it probability density function} of the standard normal distribution. %In the top two figures, the true $\sigma_\xi$ is used; in the bottom two figures, the empirical $\hat\sigma_\xi$ is used. Even though they require the same sample size for the asymptotical normality, 
}
\label{fig:P-EP}
\end{figure}

Second, we fix $m=100, r=4, \sigma_\xi=0.1$ and show the normal approximation of 
$$
\frac{\|\calP_{{U}{V}}-\calP_{\hat{U}\hat{V}}\|_{\rm F}^2-\EE\|\calP_{{U}{V}}-\calP_{\hat{U}\hat{V}}\|_{\rm F}^2}{\sqrt{8}\hat\sigma_\xi^2\|{\Lambda}^{-2}\|_{\rm F}\cdot m_{\star}^{1/2}/n}.
$$
We record $\EE\|\calP_{\hat U\hat V}-\calP_{UV}\|_{\rm F}^2$ by the average of 10000 simulations. The empirical noise variance $\hat\sigma_\xi^2$ is calculated from eq. (\ref{eq:hat_sigma}). 
For each $n=1600,2000,2400,2800$,  we record the statistics from $10000$ independent simulations and draw the density histogram. 
 The density histogram and the {\it probability density function} of the standard normal distribution are displayed in Figure~\ref{fig:P-EP}.  It shows that the normal approximation is actually very good. 

\begin{figure}
\centering

\begin{subfigure}{0.49\linewidth}
\includegraphics[scale=0.45]{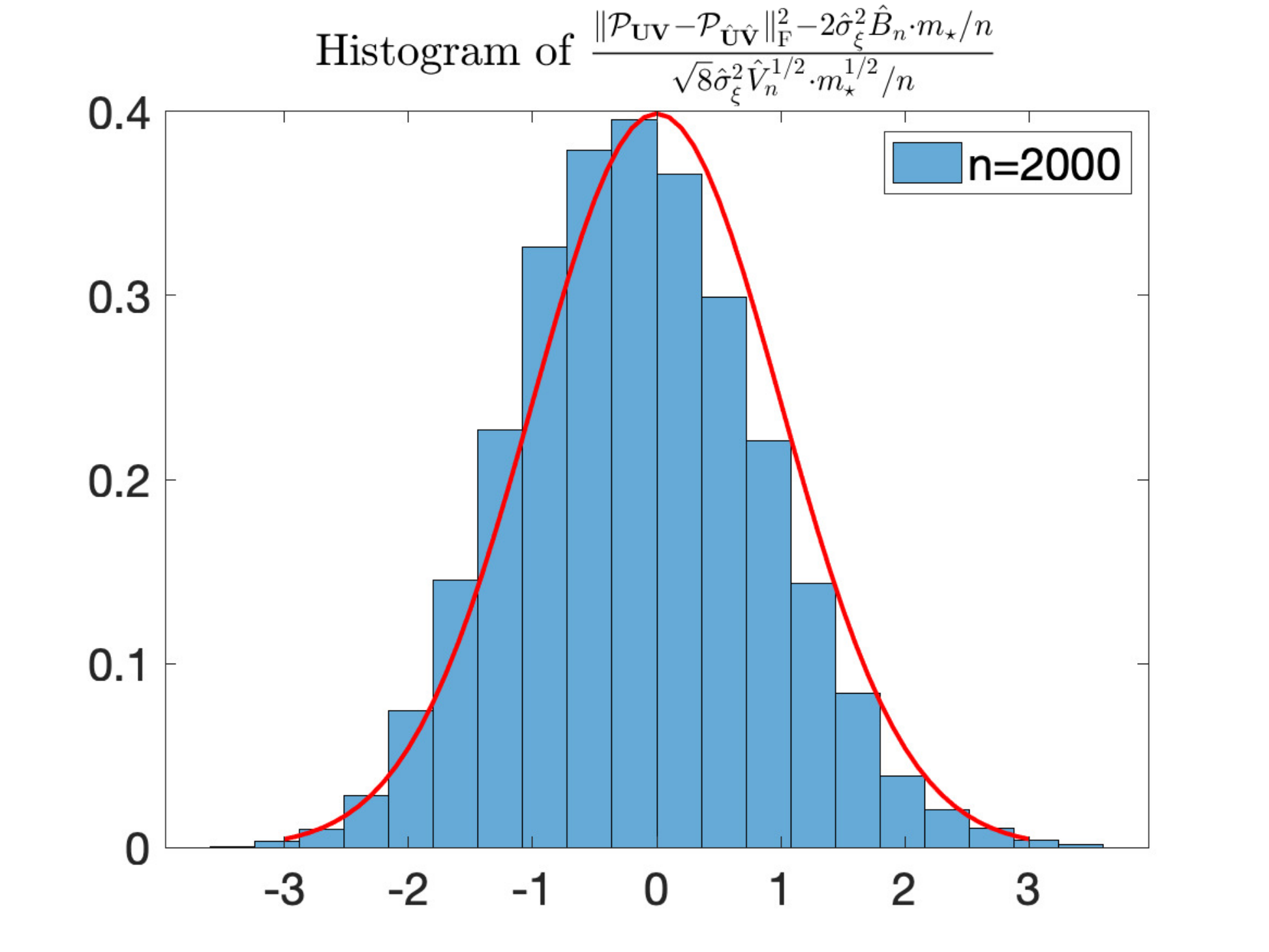}
\end{subfigure}
\begin{subfigure}{0.49\linewidth}
\includegraphics[scale=0.45]{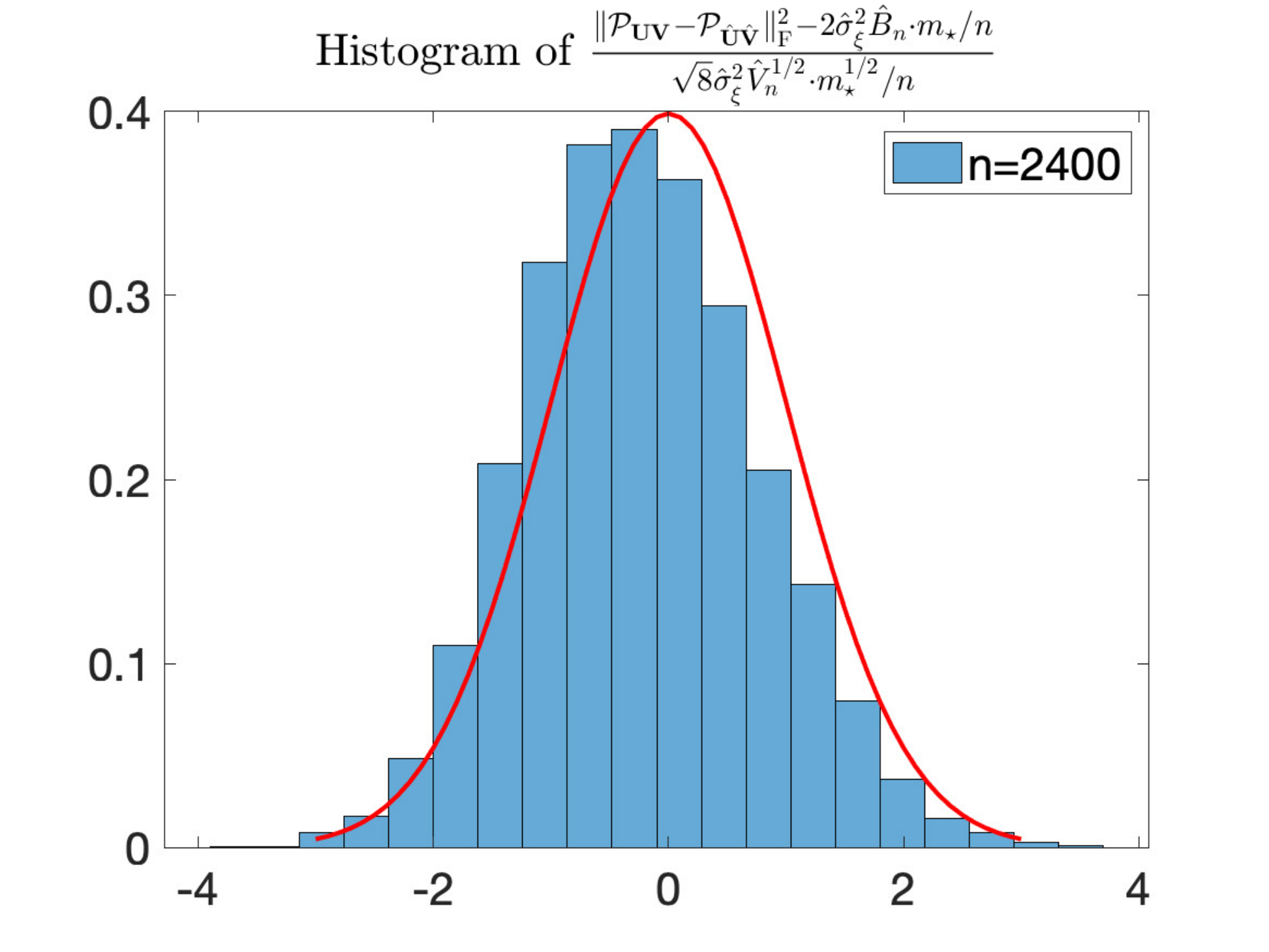}
\end{subfigure}

\begin{subfigure}{0.49\linewidth}
\includegraphics[scale=0.45]{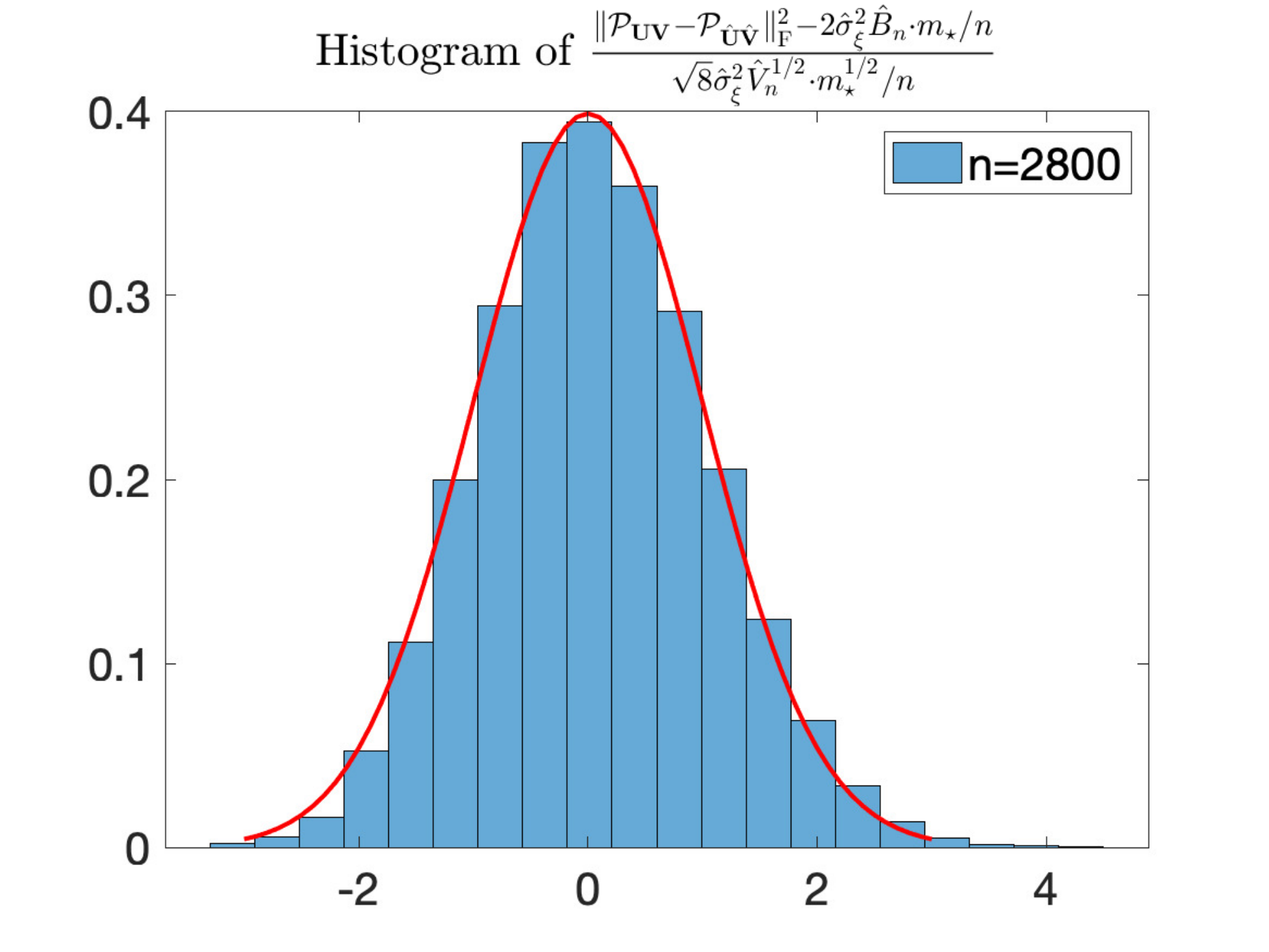}
\end{subfigure}
\begin{subfigure}{0.49\linewidth}
\includegraphics[scale=0.45]{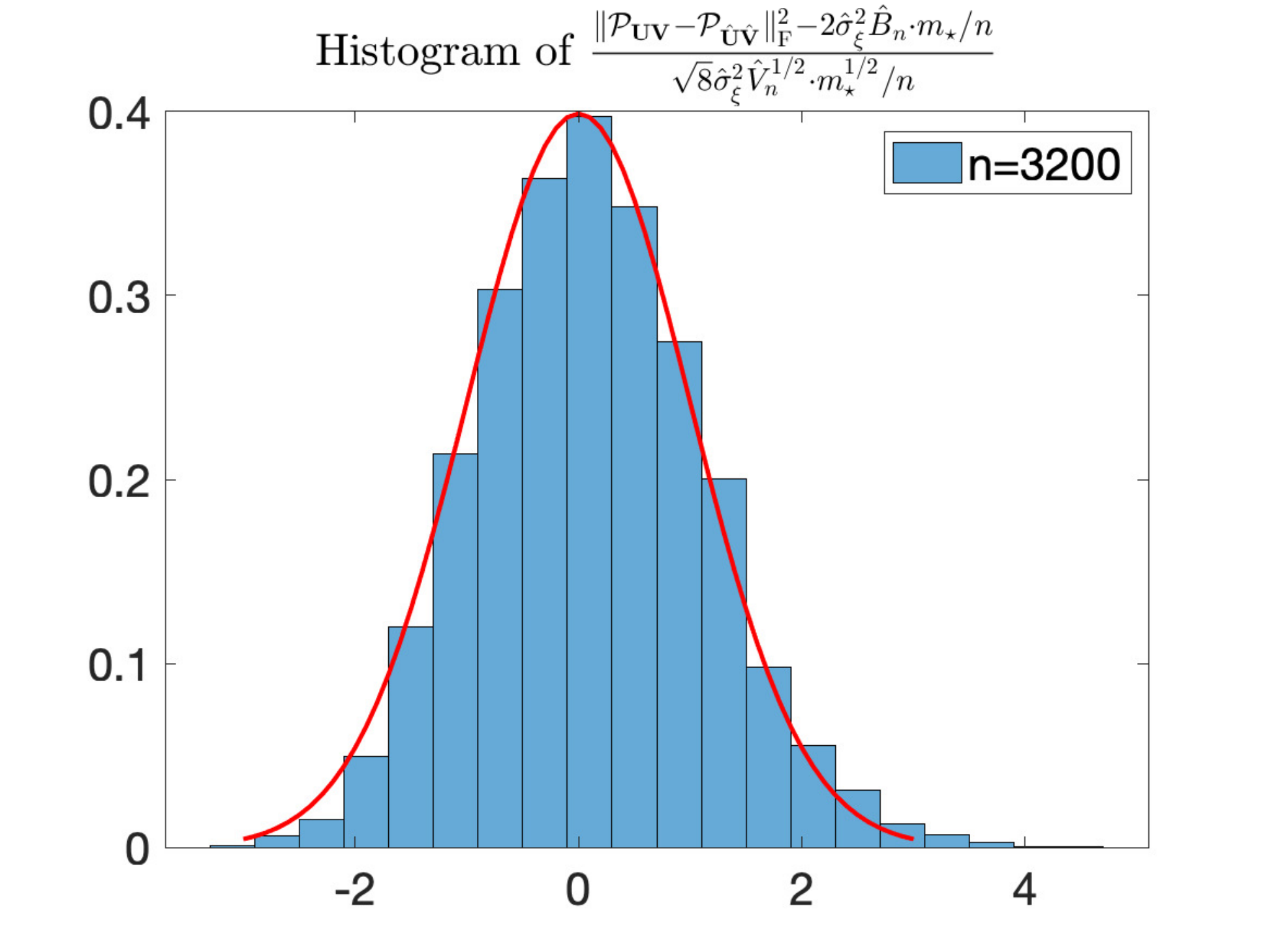}
\end{subfigure}
\caption{Normal approximation of $\frac{\|\calP_{{U}{V}}-\calP_{\hat{U}\hat{V}}\|_{\rm F}^2-\hat\sigma_\xi^2\hat B_n\cdot 2m_{\star}/n}{\sqrt{8}\hat\sigma_\xi^2\hat V_n^{1/2}\cdot m_{\star}^{1/2}/n}$ with $m_1=m_2=100, r=4$ and $\sigma_\xi=0.1$. For each $n$, the density histogram is based on $10000$ repetitions. The  estimators $\hat\sigma_\xi^2$, $\hat B_n$ and $\hat V_n$ are defined as (\ref{eq:hat_sigma}), eq.  (\ref{eq:bias_Bn}) and eq. (\ref{eq:std_Vn}).  The red curve represents the {\it probability density function} of the standard normal distribution. }
\label{fig:P-B_n}
\end{figure}

Third, we fix $m=100, r=4, \sigma_\xi=0.1$ and show the normal approximation of 
$$
\frac{\|\calP_{{U}{V}}-\calP_{\hat{U}\hat{V}}\|_{\rm F}^2-\hat\sigma_\xi^2\hat B_n\cdot 2m_{\star}/n}{\sqrt{8}\hat \sigma_\xi^2\hat V_n^{1/2}\cdot m_{\star}^{1/2}/n}
$$
where $\hat\sigma_\xi^2$, $\hat B_n$ and $\hat V_n$ are estimators as in (\ref{eq:hat_sigma}), (\ref{eq:bias_Bn}) and (\ref{eq:std_Vn}). The simulation is repeated for $10000$ times and the statistics are recorded. 
The density histogram and the {\it probability density function} of the standard normal distribution are displayed in Figure~\ref{fig:P-B_n}. The normal approximation looks very good even when $n=2800$. 

\section{Discussion}\label{sec:discussion}
%For constructing the confidence region in Section~\ref{sec:CI}, we assumed that the noise variance $\sigma_\xi^2$ is known in advance, which can be easily violated in practice.  However, $\sigma_\xi^2$ can be sharply estimated under the settings of Section~\ref{sec:CI}. One plug-in approach is 
%\begin{align*}
%\hat\sigma_\xi^2:=&\frac{1}{n}\sum_{i=n+1}^{2n}\big(y_i-\tr({X}_i^{\top}\hat{M}^{\rm nuc})\big)^2\\
%=&\frac{1}{n}\sum_{i=n+1}^{2n}\xi_i^2+\frac{1}{n}\sum_{i=n+1}^{2n}\big<\bDelta,{X}_i\big>^2+\frac{2}{n}\sum_{i=n+1}^{2n}\xi_i\cdot\big<\bDelta,{X}_i\big>.
%\end{align*}
%By Proposition~\ref{prop:nuclear_pen}, it is straightforward to show
%\begin{align*}
%\bigg|\frac{\hat\sigma_\xi^2-\sigma_\xi^2}{\sigma_\xi^2}\bigg|=O_p\bigg(\frac{r^{1/2}\bar{m}^{1/2}}{n^{1/2}}\bigg)=O_p\bigg(\frac{1}{r^{3/4}\bar{m}^{1/4}}\bigg)
%\end{align*}
%where the last equality is due to the requirement $n\gg r^{5/2}\bar{m}^{3/2}$ in Corollary~\ref{cor:CI_alpha}. Therefore, if $\bar{m}\to\infty$, we get $\Big|\frac{\hat\sigma_\xi^2}{\sigma_\xi^2}-1\Big|\to 0$ in probability.
In this paper, we assume that the rank $r$ is known. 
Otherwise, the rank $r$ can be exactly estimated from data under the similar settings. Indeed, by Lemma~\ref{lem:E_op}, we get that, if $n\gg r\bar{m}\log \bar{m}$, then with probability at least $1-c_1e^{-c_2\bar{m}}$,
$$
\sup_{1\leq k\leq \min(m_1,m_2)} \big|\hat\lambda_k-\lambda_k \big|\leq C_1\sigma_\xi\frac{\bar{m}^{1/2}}{n^{1/2}}
$$
for absolute constants $c_1,c_2,C_1>0$ where $\lambda_k=0$ for $k>r$. By choosing
$$
\hat r={\rm Card}\Big(\Big\{\hat\lambda_k: \hat\lambda_k\geq 2C_1\hat\sigma_\xi\cdot\frac{\bar{m}^{1/2}}{n^{1/2}}, 1\leq k\leq \min(m_1,m_2)\Big\}\Big)
$$
, then we have $\PP\big(\hat r=r\big)\geq 1-c_1e^{-c_2\bar{m}}$ as long as $n\gg \beta^2\bar{m}+r\bar{m}\log \bar{m}$. 

To construct the unbiased estimator $\hat{M}$ as in eq. (\ref{eq:debias_M}), our procedure splits the data $\{({X}_i,y_i)\}_{i=1}^{2n}$ into two independent samples which might be inefficient when $n$ has a moderate size. This loss of efficiency can be overcame by applying the ``double-sample-splitting" trick introduced in \cite{chernozhukov2018double}. The core idea is to  flip the role of the main and auxiliary samples to obtain a second version of the estimator. To be more specific, we compute $\hat M_1^{\rm nuc}$ from the first data sample and de-bias it using the second data sample which produces $\hat M^{(1)}$. Then, we repeat the process but using the second data sample for computing $\hat M_2^{\rm nuc}$ and the first data sample for de-biasing which produces $\hat M^{(2)}$. Finally, we calculate the average $\hat M=(\hat M^{(1)}+\hat M^{(2)})/2$. Then, we have 
\begin{align*}
\hat M=M+\frac{1}{2n}\sum_{i=1}^{2n}\xi_i X_i +&\frac{1}{2}\Big(\frac{1}{n}\sum_{i=n+1}^{2n}\tr(X_i^{\top}\Delta^{(1)})X_i-\Delta^{(1)}\Big)\\
+&\frac{1}{2}\Big(\frac{1}{n}\sum_{i=1}^{n}\tr(X_i^{\top}\Delta^{(2)})X_i-\Delta^{(2)}\Big)
\end{align*}
where $\Delta^{(1)}=M-\hat M_1^{\rm nuc}$ and $\Delta^{(2)}=M-\hat M_2^{\rm nuc}$. As a result, we can regain the full efficiency. The normal approximation of $\hat M$'s singular subspace can be proved in a similar fashion and will not be pursued in this article.

%We are curious about the performance if we reuse the same data to debias $\hat{M}^{\rm nuc}$, i.e.,
%\begin{align*}
%\doublehat{{M}}:=&\hat{M}^{\rm nuc}+\frac{1}{n}\sum_{i=1}^{n}\big(y_i-\tr({X}_i^{\top}\hat{M}^{\rm nuc})\big){X}_i\\
%=&{M}+\underbrace{\frac{1}{n}\sum_{i=1}^{n}\xi_i{X}_i}_{{Z}_1}+\underbrace{\Big[\frac{1}{n}\sum_{i=1}^n\tr(\bDelta^{\top}{X}_i){X}_i-\bDelta\Big]}_{{Z}_2}
%\end{align*}
%where $\{\xi_i\}_i$ are independent with $\{{X}_i\}_i$. However, $\bDelta$ depends on $\{\xi_i\}_i$ and $\{{X}_i\}_i$. By a covering argument, it is easy to show
%$$
%\|{Z}_2\|=O_P\bigg(\|\bDelta\|_{\rm F}\cdot\frac{\bar{m}^{1/2}\rank^{1/2}(\bDelta)}{n^{1/2}}\bigg)
%$$
%while 
%$$
%\|{Z}_1\|=O_P\bigg(\sigma_\xi\cdot \frac{\bar{m}^{1/2}}{n^{1/2}}\bigg).
%$$
%If $n\geq r^{5/2}\bar{m}^{3/2}$ as required in Corollary~\ref{cor:CI_alpha}, then $\|{Z}_2\|_{\rm F}=O_p\Big(\frac{\sigma_\xi}{r^{3/4}\bar{m}^{1/4}}\cdot \frac{\bar{m}^{1/2}\rank^{1/2}(\bDelta)}{n^{1/2}}\Big)$ implying that $\|{Z}_2\|\ll \|{Z}_1\|$ with high probability if $\rank(\bDelta)=O_P(r)$. Therefore, the dominating term of $\doublehat{M}-{M}$ is ${Z}_1$ where the technique tools in this paper is still applicable. However, it requires nontrivial effort to investigate the concentration of quantities involving ${Z}_2$. 
\section{Proofs}\label{sec:proofs}

%\begin{assump}\label{assump:cond}
%The condition number of $M$, $\kappa(M)=\frac{\sigma_1(M)}{\sigma_r(M)}\leq \kappa_0$.
%\end{assump}

%\begin{assump}\label{assump:SNR}
%For some $\gamma\in(0,1)$,
%$$
%\sigma_r(M)\geq 2(1+\gamma)\sigma_{\xi}\sqrt{\frac{2m\log m}{n}}
%$$
%\end{assump}

\subsection{Proof of Lemma~\ref{lem:hatcalP_N}, Theorem~\ref{thm:con_LNE}, Theorem~\ref{thm:hatPN-PN-frob-con}}

\subsubsection{Proof of Lemma~\ref{lem:hatcalP_N}}
First, we focus on the event $\calE_0:=\{\lambda_r\geq 5\|{E}\|\}$. We apply the representation formula of empirical spectral projectors developed in \cite{xia2019data}. For notational simplicity, we write
$$
\calP_{{U}{V}}^{-p}=\calC_{{U}{V}}^{-p}:=\sum_{1\leq |k|\leq r}\frac{1}{\lambda_k^p}({\theta}_k{\theta}_k^{\top})=
\begin{cases}
\left(\begin{array}{cc}
{0}&{U}{\Lambda}^{-p}{V}^{\top}\\
{V}{\Lambda}^{-p}{U}^{\top}&{0}
\end{array}
\right),& \textrm{ if } p \textrm{ is odd}\\
&\\
\left(\begin{array}{cc}
{U}{\Lambda}^{-p}{U}^{\top}& \bf{0}\\
{0}&{V}{\Lambda}^{-p}{V}^{\top}
\end{array}
\right),& \textrm{ if } p \textrm{ is even}. 
\end{cases}
$$
Therefore, by \cite[Theorem~1]{xia2019data}, we get
\begin{align*}
\calP_{\hat{U}\hat{V}}-\calP_{{U}{V}}=\sum_{k\geq 1}\calS_{{N},k}({E})
\end{align*}
where $\calS_{{N},k}({E})$ is defined as
\begin{align}\label{eq:SN_k}
\calS_{{N},k}({E})=\sum_{\bs: s_1+\cdots+s_{k+1}=k} (-1)^{1+\tau(\bs)}\cdot \calP_{{U}{V}}^{-s_1}{E}\calP_{{U}{V}}^{-s_2}{E}\cdots{E}\calP_{{U}{V}}^{-s_{k+1}}
\end{align}
where $s_1,\cdots,s_{k+1}\geq 0$ and we denote $\calP_{{U}{V}}^{0}=\calP_{{U}{V}}^{\perp}$ and 
$
\tau(\bs):=\sum_{j\geq 1}^{k+1}{\bf 1}(s_{j}>0).
$
 Now, we define
$$
\calL_{N}({E})=\calS_{N,1}({E})\quad {\rm and}\quad \calS_{{N}}({E})=\sum_{k\geq 2}\calS_{N,k}({E})
$$
so that $\calP_{\hat U\hat V}-\calP_{UV}=\calL_N(E)+\calS_N(E)$. Clearly, we have
\begin{align*}
\big\|\calS_{{N}}({E}) \big\|\leq \sum_{k\geq 2} {2k\choose k}\cdot \frac{\|{E}\|^k}{\lambda_r^k}\leq \sum_{k\geq 2}\frac{(4\|{E}\|)^k}{\lambda_r^k}\leq \frac{80\|{E}\|^2}{\lambda_r^2}
\end{align*}
where the last inequality is due to the fact $\lambda_r\geq 5\|{E}\|$ on event $\calE_0$. On the other hand, on event $\calE_0^{\rm c}$, we have $\lambda_r<5\|{E}\|$. Then, we have
\begin{align*}
\big\|\calP_{\hat{U}\hat{V}}-\calP_{{U}{V}}-\calL_{{N}}({E}) \big\|\leq& 2+\frac{2\|{E}\|}{\lambda_r}\leq 2\cdot \frac{25\|{E}\|^2}{\lambda_r^2}+\frac{2\|{E}\|}{\lambda_r}\cdot \frac{5\|{E}\|}{\lambda_r}\\
\leq&\frac{60\|{E}\|^2}{\lambda_r^2}
\end{align*}
where the claimed bound holds immediately. 

\subsubsection{Supporting lemmas}
The proof of Theorem~\ref{thm:con_LNE} and Theorem~\ref{thm:hatPN-PN-frob-con} involves several lemmas.
Observe that $\|\calP_{\hat{U}\hat{V}}-\calP_{{U}{V}}\|_{\rm F}^2=\|\calP_{\hat{U}\hat{V}}\|_{\rm F}^2+\|\calP_{{U}{V}}\|_{\rm F}^2-2\big<\calP_{\hat{U}\hat{V}},\calP_{{U}{V}}\big>$. By the definitions of $\calP_{\hat{U}\hat{V}}$ and $\calP_{{U}{V}}$, we have $\|\calP_{\hat{U}\hat{V}}\|_{\rm F}^2\equiv \|\calP_{{U}{V}}\|_{\rm F}^2=2r$.
Therefore, we get
\begin{align*}
\|\calP_{\hat{U}\hat{V}}-\calP_{{U}{V}}\|_{\rm F}^2-\mathbb{E}\|\calP_{\hat{U}\hat{V}}-\calP_{{U}{V}}\|_{\rm F}^2
&=-2\big<\calP_{\hat{U}\hat{V}},\calP_{{U}{V}}\big>
+2\mathbb{E}\big<\calP_{\hat{U}\hat{V}},\calP_{{U}{V}}\big>\\
&=2\big<\EE\calP_{\hat{U}\hat{V}}-\calP_{{U}{V}}, \calP_{{U}{V}}\big>
\end{align*}
Recall the formula from Lemma~\ref{lem:hatcalP_N} that
\begin{equation*}
 \calP_{\hat{U}\hat{V}}=\calP_{{U}{V}}+\calL_{N}({E})+\calS_{N}({E})
\end{equation*}
where $\calL_{N}({E})=\calP_{{U}{V}}^{\perp}{E}\calC_{{U}{V}}+\calC_{{U}{V}}{E}\calP_{{U}{V}}^{\perp}$. Then, we write
\begin{equation*}
 \big<\mathbb{E}\calP_{\hat{U}\hat{V}}-\calP_{{U}{V}},\calP_{{U}{V}}\big>=\big<\mathbb{E}\calS_{N}({E})-\calS_{N}({E})-\calL_{N}({E}),\calP_{{U}{V}}\big>
\end{equation*}
where we used the fact $\EE\calL_{N}({E})={0}$. By the definition of $\calC_{UV}$, we have
\begin{gather*}
 \langle \calL_{N}({E}),\calP_{{U}{V}}\rangle=\langle \calP_{{U}{V}}^{\perp}{E}\calC_{{U}{V}}+\calC_{{U}{V}}{E}\calP_{{U}{V}}^{\perp},\calP_{{U}{V}}\rangle=0.
\end{gather*}
Finally, we conclude that
\begin{gather}\label{eq:hatP-P_ES-S}
 \|\calP_{\hat{U}\hat{V}}-\calP_{{U}{V}}\|_{\rm F}^2-\mathbb{E}\|\calP_{\hat{U}\hat{V}}-\calP_{{U}{V}}\|_{\rm F}^2=2\big<\mathbb{E}\calS_{N}({E})-\calS_{N}({E}),\calP_{{U}{V}}\big>.
\end{gather}
%For any $t\geq 1$, and the constants $C_1,C_2,C_3,C_4$ in Lemma~\ref{lem:E_op}, define
%\begin{align}
%\delta_n(t):=&C_1\|\bDelta\|_{\rm F}\frac{(\bar{m}\log\bar{m})^{1/2}}{n^{1/2}}+C_2\sigma_{\xi}\frac{\bar{m}^{1/2}}{n^{1/2}}\nonumber\\
%&+C_3\|\bDelta\|_{\rm F}\bigg(\frac{t^{1/2}}{n^{1/2}}+\frac{t}{n}\bigg)+C_4\sigma_{\xi}\frac{t^{1/2}}{n^{1/2}}.\label{eq:delta_n_t}
%\end{align}
\begin{lemma}\label{HSlemma}
Suppose that $\lambda_r\geq 10\EE\|{E}\|$ and $n\geq C_1r\bar{m}$ for some large enough absolute constant $C_1>0$, the following bound holds  with probability at least at least $1-c_1e^{-c_2\bar{m}}-e^{-n}$,
\begin{align*}
\big|\|\calP_{\hat{U}\hat{V}}-&\calP_{{U}{V}}\|_{\rm F}^2-\|\calL_{N}({E})\|_{\rm F}^2\big|\leq 20r\cdot\Big(\frac{8\EE\|{E}\|}{\lambda_r}\Big)^3
\end{align*}
where $c_1,c_2>0$ are absolute constants.
Moreover, with probability at least $1-4e^{-t}-2e^{-n}-c_1e^{-c_2\bar{m}}$ for all $t\in[1,n]$,
\begin{align*}
 \Big|\Big(\|\calP_{\hat{U}\hat{V}}-\calP_{{U}{V}}\|_{\rm F}^2&-\mathbb{E}\|\calP_{\hat{U}\hat{V}}-\calP_{{U}{V}}\|_{\rm F}^2\Big)-\Big(\|\calL_{N}({E})\|_{\rm F}^2-\EE\|\calL_{N}({E})\|_{\rm F}^2\Big)\Big|\\
 \leq&C_7r\frac{\sigma_\xi}{\lambda_r}\Big(\frac{9\EE\|{E}\|}{\lambda_r}\Big)^2\cdot\frac{t^{1/2}+\log^{1/2}{\bar m}}{n^{1/2}}
\end{align*}
for an absolute constant $c_1,c_2, C_7>0$.
\end{lemma}
It is thus sufficient to investigate the normal approximation of $\|\calL_{N}({E})\|_{\rm F}^2$. By the definition of $\calL_{N}({E})$, we get
\begin{equation*}
\begin{split}
\|\calL_{N}({E})\|_{\rm F}^2=&\|\calP_{{U}{V}}^{\perp}{E}\calC_{{U}{V}}+\calC_{{U}{V}}{E}\calP_{{U}{V}}^{\perp}\|_{\rm F}^2\\
=&\|\calP_{{U}{V}}^{\perp}{E}\calC_{{U}{V}}\|_{\rm F}^2+2\big<\calP_{{U}{V}}^{\perp}{E}\calC_{{U}{V}},\calC_{{U}{V}}{E}\calP_{{U}{V}}^{\perp}\big>+\|\calC_{{U}{V}}{E}\calP_{{U}{V}}^{\perp}\|_{\rm F}^2
=2\|\calP_{{U}{V}}^{\perp}{E}\calC_{{U}{V}}\|_{\rm F}^2\\
=&2\|\calP_{{U}{V}}^{\perp}{E}_1\calC_{{U}{V}}\|_{\rm F}^2+2\|\calP_{{U}{V}}^{\perp}{E}_2\calC_{{U}{V}}\|_{\rm F}^2+4\langle\calP_{{U}{V}}^{\perp}{E}_1\calC_{{U}{V}}, \calP_{{U}{V}}^{\perp}{E}_2\calC_{{U}{V}}\rangle,
\end{split}
\end{equation*}
where the third equality is due to the fact that $\calP_{{U}{V}}^{\perp}\calC_{{U}{V}}={0}$. Recall the definitions of $\calP_{{U}{V}}^{\perp}$ and $\calC_{{U}{V}}$, we can write $\calP_{{U}{V}}^{\perp}{E}_1\calC_{{U}{V}}$ explicitly as
\begin{gather*}
\calP_{{U}{V}}^{\perp}{E}_1\calC_{{U}{V}}=\left(\begin{array}{cc}{U}_{\perp}{U}_{\perp}^{\top}&{0}\\
{0}&{V}_{\perp}{V}_{\perp}^{\top}\end{array}\right)
\left(\begin{array}{cc}{0}& {Z}_1\\
{Z}_1^{\top}&{0}\end{array}\right)
\left(\begin{array}{cc}{0}&{U}{\Lambda}^{-1}{V}^{\top}\\
{V}{\Lambda}^{-1}{U}^{\top}&{0}\end{array}\right)\\
=\left(
\begin{array}{cc}
{U}_{\perp}{U}_{\perp}^{\top}{Z}_1{V}{\Lambda}^{-1}{U}^{\top}&{0}\\
{0}&{V}_{\perp}{V}_{\perp}^{\top}{Z}_1^{\top}{U}{\Lambda}^{-1}{V}^{\top}
\end{array}
\right)
\end{gather*}
implying that
\begin{gather*}
\|\calP_{{U}{V}}^{\perp}{E}_1\calC_{{U}{V}}\|_{\rm F}^2=\big\|{U}_{\perp}{U}_{\perp}^{\top}{Z}_1{V}{\Lambda}^{-1}{U}^{\top}\big\|_{\rm F}^2+\big\|{V}_{\perp}{V}_{\perp}^{\top}{Z}_1^{\top}{U}{\Lambda}^{-1}{V}^{\top}\big\|_{\rm F}^2.
\end{gather*}
\begin{lemma}\label{lem:PN-E-CN}
Let $\{z_k^2\}_{k=1}^{r}$ be i.i.d. Chi-squared random variables with degrees of freedom $m_{\star}$ where $m_{\star}=m_1+m_2-2r$. Then,
\begin{gather*}
\|\calP_{{U}{V}}^{\perp}{E}_1\calC_{{U}{V}}\|_{\rm F}^2\stackrel{{\rm d}}{=}\Big(\frac{1}{n^2}\sum_{i=n+1}^{2n}\xi_i^2\Big)\sum_{k=1}^{r}\frac{z_k^2}{\lambda_k^2}.
\end{gather*}
Therefore, $\EE \|\calP_{{U}{V}}^{\perp}{E}_1\calC_{{U}{V}}\|_{\rm F}^2=\sigma_{\xi}^2\Big(\frac{m_{\star}}{n}\Big)\|\Lambda^{-1}\|_{\rm F}^2$. Meanwhile, for any $t\geq 1$, we get that with probability at least $1-e^{-t}-e^{-n}$,
\begin{align*}
\big|\|\calP_{{U}{V}}^{\perp}{E}_1\calC_{{U}{V}}\|_{\rm F}^2-&\EE\|\calP_{{U}{V}}^{\perp}{E}_1\calC_{{U}{V}}\|_{\rm F}^2 \big|
\leq C_1\frac{\sigma_\xi^2}{n}\cdot\max\bigg\{\|{\Lambda}^{-2}\|_{\rm F}m_{\star}^{1/2}t^{1/2}, \frac{t}{\lambda_r^2}\bigg\}
\end{align*}
for some absolute constant $C_1>0$. 
\end{lemma}

\begin{lemma}\label{lemma:PN-E2-CN-frob}
Under the conditions of Lemma~\ref{HSlemma}, 
the following bounds hold with probability at least $1-3e^{-t}-c_1e^{-c_2\bar{m}}-2e^{-n}$, 
\begin{align*}
\Big|\|\calP_{{U}{V}}^{\perp}{E}_2\calC_{{U}{V}}\|_{\rm F}^2&-\EE\|\calP_{{U}{V}}^{\perp}{E}_2\calC_{{U}{V}}\|_{\rm F}^2\Big|
\leq C_7\sigma_\xi^2\|{\Lambda}^{-1}\|_{\rm F}^2\cdot \frac{t^{1/2}+\log^{1/2}{\bar m}}{n^{1/2}}\cdot \frac{r\bar{m}}{n}%+C_2\|{\Lambda}^{-1}\|_{\rm F}^2\|\bDelta\|_{\rm F}^2\frac{\bar{m}t^{1/2}\log^{1/2}{\bar m}}{n^{3/2}}
\end{align*}
for absolute constants $c_1,c_2,C_7>0$.
Meanwhile, the following bound holds
$$
\EE\|\calP_{{U}{V}}^{\perp}{E}_2\calC_{{U}{V}}\|_{\rm F}^2\leq C_1\frac{r\bar{m}^2}{n^2}\sigma_\xi^2\cdot \|{\Lambda}^{-1}\|_{\rm F}^2+C_1\frac{r\bar{m}}{n^2}\cdot \frac{\sigma_\xi^2}{\lambda_r^2}.
$$
Meanwhile, conditioned on $\{(X_i, y_i)\}_{i=1}^n$, we have $\EE\|\calP_{{U}{V}}^{\perp}{E}_2\calC_{{U}{V}}\|_{\rm F}^2\geq c_1m_{\star}\|\Delta\|_{\rm F}^2\|\Lambda^{-1}\|_{\rm F}^2/n$ for some absolute constant $c_1>0$. 
\end{lemma}
By the independence between $\{\xi_i\}_{i}$ and $\{X_i\}_{i}$, we have $\EE\langle\calP_{{U}{V}}^{\perp}{E}_1\calC_{{U}{V}}, \calP_{{U}{V}}^{\perp}{E}_2\calC_{{U}{V}}\rangle=0$. 

\begin{lemma}\label{lem:PN-E1-CN-PN-E2-CN}
The following bound holds with probability at least $1-(2n+1)e^{-t}-ne^{-m_{\star}}-c_1e^{-c_2\bar{m}}$ for all $t\geq 2\log n$,
\begin{align*}
\big|\langle\calP_{{U}{V}}^{\perp}&{E}_1\calC_{{U}{V}}, \calP_{{U}{V}}^{\perp}{E}_2\calC_{{U}{V}}\rangle\big|
\leq C_6r^{1/2}\sigma_\xi^2\|{\Lambda}^{-2}\|_{\rm F}\cdot \frac{t^{3/2}\bar{m}}{n^{3/2}}+C_7r^{1/2}\sigma_\xi^2\|{\Lambda}^{-1}\|_{\rm F}^2\cdot \frac{t\bar{m}^{3/2}}{n^2}
\end{align*}
for absolute constants $c_1,c_2,C_6,C_7>0$.
\end{lemma}

By combining Lemma~\ref{lem:PN-E-CN}, Lemma~\ref{lemma:PN-E2-CN-frob} and Lemma~\ref{lem:PN-E1-CN-PN-E2-CN} with $t=2\log n$, we will prove the concentration of $\|\calL_{N}({E})\|_{\rm F}^2$.

\subsubsection{Proof of Theorem~\ref{thm:con_LNE}}
By putting together the bounds in Lemma~\ref{lem:PN-E-CN}, Lemma~\ref{lemma:PN-E2-CN-frob} and Lemma~\ref{lem:PN-E1-CN-PN-E2-CN} with $t=2\log n$, we immediate obtain, with probability at least $1-\frac{2n+5}{n^2}-2e^{-n}-n e^{-m_{\star}}-c_1e^{-c_2\bar{m}}$ that
\begin{align*}
\big|\|\calL_{N}({E})\|_{\rm F}^2&-\EE\|\calL_{N}({E})\|_{\rm F}^2\big|\\
\leq& C_1\sigma_\xi^2\|{\Lambda}^{-2}\|_{\rm F}\cdot\bigg[\frac{\bar{m}^{1/2}\log^{1/2}n}{n}+\frac{r^{1/2}\bar{m}\log^{3/2}n}{n^{3/2}}\bigg]
+C_2\sigma_\xi^2\|{\Lambda}^{-1}\|_{\rm F}^2\cdot\frac{r\bar{m}\log n}{n^{3/2}}\\
\leq&C_1\sigma_\xi^2\|{\Lambda}^{-2}\|_{\rm F}\cdot\frac{\bar{m}^{1/2}\log^{1/2}n}{n}+C_2\sigma_\xi^2\|{\Lambda}^{-1}\|_{\rm F}^2\cdot\frac{r\bar{m}^{3/2}\log n}{n^2}
\end{align*}
where the last inequality is due to $n\geq Cr\bar{m}\log^2n$.

Since $\EE\|\calL_{N}({E})\|_{\rm F}^2=2\EE\|\calP_{{U}{V}}^{\perp}{E}_1\calC_{{U}{V}}\|_{\rm F}^2+2\EE\|\calP_{{U}{V}}^{\perp}{E}_2\calC_{{U}{V}}\|_{\rm F}^2$, we immediately obtain the second claim from Lemma~\ref{lem:PN-E-CN} and Lemma~\ref{lemma:PN-E2-CN-frob}. 

\subsubsection{Proof of Theorem~\ref{thm:hatPN-PN-frob-con}}
By Lemma~\ref{HSlemma} and setting $t=2\log n$, with probability at least $1-\frac{4}{n^2}-e^{-n}-c_1e^{-c_2\bar{m}}$,
\begin{align*}
 \Big|\Big(\|\calP_{\hat{U}\hat{V}}-\calP_{{U}{V}}\|_{\rm F}^2&-\mathbb{E}\|\calP_{\hat{U}\hat{V}}-\calP_{{U}{V}}\|_{\rm F}^2\Big)-\Big(\|\calL_{N}({E})\|_{\rm F}^2-\EE\|\calL_{N}({E})\|_{\rm F}^2\Big)\Big|\\
 \leq&C_7r\frac{\sigma_\xi}{\lambda_r}\Big(\frac{9\EE\|{E}\|}{2\lambda_r}\Big)^2\cdot\frac{\log^{1/2}n}{n^{1/2}}
 \leq C_7\Big(\frac{\sigma_\xi}{\lambda_r}\Big)^3\cdot \frac{r\bar{m}\log^{1/2}n}{n^{3/2}}.
\end{align*}
Together with the concentration of $\|\calL_{N}({E})\|_{\rm F}^2-\EE\|\calL_{N}({E})\|_{\rm F}^2$ in Theorem~\ref{thm:con_LNE}, we obtain the first claimed bound. 

We now prove the second bound. Recall that $\calP_{\hat{U}\hat{V}}-\calP_{{U}{V}}=\calL_{N}({E})+\calS_{{N}}({E})$. Therefore,
\begin{align*}
\EE \|\calP_{\hat{U}\hat{V}}-\calP_{{U}{V}}\|_{\rm F}^2=\EE\|\calL_{N}({E})\|_{\rm F}^2+\EE\|\calS_{N}({E})\|_{\rm F}^2+2\cdot\EE\big<\calL_{N}({E}),\calS_{N}({E})\big>.
\end{align*}
By Theorem~\ref{thm:con_LNE}, we have
$$
\EE\|\calL_{N}({E})\|_{\rm F}^2=\sigma_{\xi}^2\|{\Lambda}^{-1}\|_{\rm F}^2\cdot \frac{2m_{\star}}{n}+O\Big(\sigma_{\xi}^2\|{\Lambda}^{-1}\|_{\rm F}^2\cdot \frac{2r\bar{m}^2}{n^2}\Big).
$$
Since $\rank\big(\calS_{N}({E})\big)\leq \rank\big(\calP_{\hat{U}\hat{V}}-\calP_{{U}{V}}-\calL_{N}({E})\big)\leq 8r$, by Lemma~\ref{lem:hatcalP_N},
\begin{align*}
\EE\|\calS_{N}({E})\|_{\rm F}^2\leq& 8r\cdot\EE\|\calS_{N}({E})\|^2\leq C\cdot r\frac{\|{E}\|^4}{\lambda_r^4}
\leq C_1\Big(\frac{\sigma_\xi}{\lambda_r}\Big)^4\cdot \frac{r \bar{m}^2}{n^2}.
\end{align*}
The upper bound of $\EE\big<\calL_{N}({E}),\calS_{N}({E})\big>$ requires more delicate treatments. To this end, denote the event $\calE_1:=\{\|{E}\|\leq 2\EE\|{E}\|\}$ with $\PP\big(\calE_1\big)\geq 1-e^{-c_1\bar{m}}$ by the proof of Lemma~\ref{HSlemma}.

 Recall the definition of $\calS_{{N},k}({E})$ in eq. (\ref{eq:SN_k}). Then, we get
 \begin{align*}
 \EE\big<\calL_{{N}}({E}), \calS_{{N}}({E})\big>{\bf 1}_{\calE_1}=\EE\big<\calL_{{N}}({E}),\calS_{{N},2}({E})\big>{\bf 1}_{\calE_1}+\sum_{k\geq 3}\EE\big<\calL_{{N}}({E}),\calS_{{N},k}({E})\big>{\bf 1}_{\calE_1}.
 \end{align*} 
As in the proof of Lemma~\ref{lem:hatcalP_N}, we get 
\begin{align*}
\Big|\sum_{k\geq 3}\EE\big<\calL_{{N}}({E}),\calS_{{N},k}({E})\big>{\bf 1}_{\calE_1}\Big|\leq& 2r\sum_{k\geq 3}\EE \|\calL_{{N}}({E})\| \|\calS_{{N},k}({E})\|{\bf 1}_{\calE_1}\\
\leq&2r\sum_{k\geq 3}4^k\EE \frac{\|{E}\|^{k+1}}{\lambda_r^{k+1}}\leq C_1\beta^4\cdot \frac{r\bar{m}^2}{n^2}
\end{align*}
where the last inequality holds as long as $n\geq C_2\beta^2\bar{m}$ for some large enough constant $C_2>0$. 
 %
%$$
%\calS_{N}({E}){\bf 1}_{\calE_1}=-\frac{1}{2\pi i}\Big(\oint_{\gamma_{N}^+}+\oint_{\gamma_{{N}^-}}\Big)\sum_{j\geq 2}(-1)^j\big[\calR_{N}(\eta){E}\big]^j\calR_{{N}}(\eta){\bf 1}_{\calE_1}d\eta.
%$$
Now, we bound $\big| \EE\big<\calL_{{N}}({E}),\calS_{{N},2}({E})\big>{\bf 1}_{\calE_1}\big|$. In view of eq. (\ref{eq:SN_k}), W.L.O.G., 
we bound 
$$
\big| \EE\big<\calL_{{N}}({E}),\calP_{{U}{V}}^{\perp}{E}\calC_{{U}{V}}^{-2}{E}\calP_{{U}{V}}^{\perp}\big>{\bf 1}_{\calE_1}\big|.
$$
Recall that ${E}={E}_1+{E}_2$ where ${E}_1=n^{-1}\sum_i\xi_i\mathfrak{D}({X}_i)$ and ${E}_2=n^{-1}\sum_i \langle\bDelta,{X}_i \rangle\mathfrak{D}({X}_i)-\mathfrak{D}(\bDelta)$. We can write
%\begin{align*}
%\EE \Big<\big[\calR_{N}(\eta)&{E}\big]^2\calR_{{N}}(\eta), \calL_{N}({E})\Big>=\EE\tr\Big(\big[\calR_{N}(\eta){E}\big]^2\calR_{{N}}(\eta)\calL_{N}({E})\Big)\\
%=&\EE\tr\Big(\calR_{N}({E}_2)\calR_{N}(\eta){E}\calR_{N}(\eta)\calL_{N}({E})\Big)\\
%&+\EE\tr\Big(\calR_{N}({E}_1)\calR_{N}(\eta){E}_2\calR_{N}(\eta)\calL_{N}({E})\Big)\\
%&+\tr\Big(\calR_{N}({E}_1)\calR_{N}(\eta){E}_1\calR_{N}(\eta)\calL_{N}({E}_2)\Big)\\
%&+\tr\Big(\calR_{N}({E}_1)\calR_{N}(\eta){E}_1\calR_{N}(\eta)\calL_{N}({E}_1)\Big).
%\end{align*}
\begin{align*}
\EE\big<\calL_{{N}}({E}),\calP_{{U}{V}}^{\perp}{E}&\calC_{{U}{V}}^{-2}{E}\calP_{{U}{V}}^{\perp}\big>=\EE\big<\calL_{{N}}({E}_1),\calP_{{U}{V}}^{\perp}{E}_1\calC_{{U}{V}}^{-2}{E}_1\calP_{{U}{V}}^{\perp}\big>\\
+&\EE\big<\calL_{{N}}({E}),\calP_{{U}{V}}^{\perp}{E}\calC_{{U}{V}}^{-2}{E}\calP_{{U}{V}}^{\perp}\big>-\EE\big<\calL_{{N}}({E}_1),\calP_{{U}{V}}^{\perp}{E}_1\calC_{{U}{V}}^{-2}{E}_1\calP_{{U}{V}}^{\perp}\big>.
\end{align*}
Since $\xi_i\sim\calN(0,\sigma_\xi^2)$ are i.i.d., we have $\EE\big<\calL_{{N}}({E}_1),\calP_{{U}{V}}^{\perp}{E}_1\calC_{{U}{V}}^{-2}{E}_1\calP_{{U}{V}}^{\perp}\big>=0$. Together with Lemma~\ref{lem:E_op}, we get
\begin{align*}
\Big|\EE\big<\calL_{{N}}&({E}),\calP_{{U}{V}}^{\perp}{E}\calC_{{U}{V}}^{-2}{E}\calP_{{U}{V}}^{\perp}\big>\Big|\\
\leq& 4r\cdot \frac{\EE\|{E}\|^2\|{E}_2\|+\EE\|{E}_1\|\|{E}_2\|\|{E}\|+\EE\|{E}_1\|^2\|{E}_2\|}{\lambda_r^3}\leq C\frac{\sigma_\xi^3}{\lambda_r^3}\cdot \frac{r^{3/2}\bar{m}^2\log^{1/2}\bar{m}}{n^2}
\end{align*}
where we used the fact $\EE\|\bDelta\|_{\rm F}\leq C_1\sigma_\xi\frac{(r\bar{m})^{1/2}}{n^{1/2}}$. Therefore, we get
\begin{align*}
\Big|\EE \big<\calL_{{N}}&({E}),\calP_{{U}{V}}^{\perp}{E}\calC_{{U}{V}}^{-2}{E}\calP_{{U}{V}}^{\perp}\big>{\bf 1}_{\calE_1}\Big|\\
\leq& \Big|\EE  \big<\calL_{{N}}({E}),\calP_{{U}{V}}^{\perp}{E}\calC_{{U}{V}}^{-2}{E}\calP_{{U}{V}}^{\perp}\big>\Big|+\Big|\EE \big<\calL_{{N}}({E}),\calP_{{U}{V}}^{\perp}{E}\calC_{{U}{V}}^{-2}{E}\calP_{{U}{V}}^{\perp}\big>{\bf 1}_{\calE_1^{\rm c}}\Big|\\
\leq&C\frac{\sigma_\xi^3}{\lambda_r^3}\cdot \frac{r^{3/2}\bar{m}^2\log^{1/2}\bar{m}}{n^2}+Ce^{-c_1\bar{m}}\frac{\sigma_\xi^3}{\lambda_r^3}\cdot \frac{r\bar{m}^{3/2}}{n^{3/2}}
\leq C_2\frac{\sigma_\xi^3}{\lambda_r^3}\cdot \frac{r^{3/2}\bar{m}^2\log^{1/2}\bar{m}}{n^2}.
\end{align*}
To this end, we conclude that
\begin{align*}
\big|\EE\big<\calL_{N}({E}),&\calS_{N}({E})\big>{\bf 1}_{\calE_1}\big|\leq  C_1\beta^3\cdot \frac{r^{3/2}\bar{m}^2\log^{1/2}\bar{m}}{n^2}+C_3\beta^4\cdot \frac{r\bar{m}^2}{n^2}
\end{align*}
Similarly, we can show that 
\begin{align*}
\big|\EE\big<\calL_{N}({E}),&\calS_{N}({E})\big>{\bf 1}_{\calE_1^{\rm c}}\big|\leq  C_1e^{-c_1\bar{m}}\Big(\frac{\sigma_\xi}{\lambda_r}\Big)^3\cdot\frac{r\bar{m}^{3/2}}{n^{3/2}}. 
\end{align*}
As long as $e^{-c_1\bar{m}}\leq \sqrt{\frac{\bar{m}}{n}}$, we get 
\begin{align*}
\big|\EE\big<\calL_{N}({E}),\calS_{N}({E})\big>\big|\leq& C_2\beta^3\cdot \frac{r^{3/2}\bar{m}^2\log^{1/2}\bar{m}}{n^2}+C_3\beta^4\cdot \frac{r\bar{m}^2}{n^2}\\
&\leq C_3(\beta\vee 1)^4.\frac{r^{3/2}\bar{m}^2\log^{1/2}\bar{m}}{n^2}.
\end{align*}

\subsection{Proof of supporting lemmas}

\subsubsection{Proof of Lemma~\ref{HSlemma}}
Since $n\geq \bar{m}$ and let $t\leq \bar{m}$ in Lemma~\ref{lem:E_op}, we obtain that
\begin{align*}
\PP\bigg(\|{E}\|-\EE\|{E}\|\geq C_3\sigma_\xi\cdot \frac{t^{1/2}}{n^{1/2}}+C_4\|\bDelta\|_{\rm F}\cdot \frac{t^{1/2}+\log^{1/2}\bar{m}}{n^{1/2}}\bigg)\leq 1-3e^{-t}-e^{-n}.
\end{align*}
By setting $t=c_1\bar{m}$ with small enough absolute constant $c_1>0$, we conclude that with probability at least $1-3e^{-c_1\bar{m}}-e^{-n}$, $\|{E}\|\leq \frac{9}{4}\EE\|{E}\|$. 

Denote $\bar{\delta}=2\EE\|{E}\|$
%By definition of $\delta_n(t)$ in (\ref{eq:delta_n_t}) and Lemma~\ref{lem:E_op},
%$$
%\mathbb{P}\big(\|{E}\|\leq \delta_n(t)\big)\geq 1-e^{-n}-e^{-t}.
%$$
and the event $\calE_1:=\{\|{E}\|\leq \frac{9}{4}\EE\|{E}\|\}$
on which $\frac{4\|{E}\|}{\lambda_r}\leq \frac{9\EE\|{E}\|}{\lambda_r}\leq \frac{9}{10}$ and $\PP\big(\calE_1\big)\geq 1-e^{-n}-3e^{-c_1\bar{m}}$.
We define a Lipschitz function $\phi(\cdot)$ on $\mathbb{R}_+$ such that
$$
\phi(s)=
\begin{cases}
1,&\textrm{ if } 0\leq s\leq 1.\\
1-8(s-1),&\textrm{ if } 1\leq s\leq \frac{9}{8}.\\
0,&\textrm{ if } s\geq \frac{9}{8}.
\end{cases}
$$
Clearly, $\phi(s)$ is Lipschitz with constant $8$.  By the definition of $\phi(\cdot)$, we have $\phi\big(\|{E}\|/\bar\delta\big)=1$ if $\|{E}\|\leq 2\EE\|{E}\|$. 
Observe that on event $\calE_2:=\big\{\|{E}\|\leq 2\EE\|{E}\|\big\}$ with $\PP(\calE_2)\geq 1-e^{-n}-c_1e^{-c_2\bar{m}}$ for some $c_1,c_2>0$, we have
\begin{align}
\label{eq:proofHSlemmaineq1}
  \big|&\big(\|\calP_{\hat U\hat V}-\calP_{UV}\|_{\rm F}^2-\EE\|\calP_{\hat U\hat V}-\calP_{UV}\|_{\rm F}^2\big)-\big(\|\calL_N(E)\|_{\rm F}^2-\EE\|\calL_N(E)\|_{\rm F}^2\big)\big|\nonumber\\
=&  \big|\big(\|\calP_{\hat U\hat V}-\calP_{UV}\|_{\rm F}^2-\EE\|\calP_{\hat U\hat V}-\calP_{UV}\|_{\rm F}^2\big)-\big(\|\calL_N(E)\|_{\rm F}^2-\EE\|\calL_N(E)\|_{\rm F}^2\big)\big|\phi\Big(\frac{\|{E}\|}{\bar\delta}\Big).
\end{align}
It suffices to focus on the concentration of the right hand side of (\ref{eq:proofHSlemmaineq1}).
%We note that $\phi\big(\frac{\|{E}\|}{\delta_n(t)}\big)$ is nonzero if and only if
%$$
%\|{E}\|\leq (1+\frac{\gamma}{2})\delta_n(t)\leq (1+\frac{\gamma}{2})(1-\frac{\gamma}{2})\frac{\sigma_r(M)}{2}=\big(1-\frac{\gamma^2}{4}\big)\frac{\sigma_r(M)}{2}.
%$$
Eq. (\ref{eq:proofHSlemmaineq1}) is equivalent to the 
concentration of $\big(\|\calP_{\hat U\hat V}-\calP_{UV}\|_{\rm F}^2-\|\calL_N(E)\|_{\rm F}^2\big)\phi\big(\frac{\|{E}\|}{\bar\delta}\big)$ around its expectation. 

 Since everything is trivial on the event $\calE_1^{\rm c}$,  the following analysis shall be focused on event $\calE_1$. 
By eq. (\ref{eq:SN_k}), we can write
\begin{align*}
\big(\|\calP_{\hat U\hat V}-\calP_{UV}\|_{\rm F}^2-\|\calL_N(E)\|_{\rm F}^2\big)\phi\Big(\frac{\|{E}\|}{\bar\delta}\Big)=-2\sum_{k\geq 3}\big<\calS_{{N},k}({E}), \calP_{{U}{V}}\big>\phi\Big(\frac{\|{E}\|}{\bar{\delta}}\Big)
\end{align*}
where we used the fact $-2\big<\calS_{{N},2}({E}),\calP_{{U}{V}}\big>=\|\calL_{{N}}({E})\|_{\rm F}^2$. 

\paragraph*{Proof of first claim} By the above representation, on the event $\calE_2$, we have
\begin{align*}
\big|\|\calP_{\hat{U}\hat{V}}-&\calP_{{U}{V}}\|_{\rm F}^2-\|\calL_{N}({E})\|_{\rm F}^2\big|
\leq4r\sum_{k\geq 3}\big\|\calS_{{N},k}({E})\big\|\leq 4r\sum_{k\geq 3}{2k \choose k}\cdot \frac{\|{E}\|^k}{\lambda_r^{k}}\\
\leq&4r\Big(\frac{4\bar\delta}{\lambda_r}\Big)^3\sum_{k\geq 0}\Big(\frac{8\EE\|E\|}{\lambda_r}\Big)^k\leq 4r\Big(\frac{4\bar\delta}{\lambda_r}\Big)^3\sum_{k\geq 0}\Big(\frac{8}{10}\Big)^k\leq 20r\Big(\frac{4\bar\delta}{\lambda_r}\Big)^3.
\end{align*}

\paragraph*{Proof of second claim} It suffices to prove the concentration inequality for the following functions,
\begin{equation*}
 \varphi_{\bar{\delta}}({E}):=\sum_{k\geq 3} \varphi_{k,\bar{\delta}}({E})
\end{equation*}
\begin{equation*}
 \varphi_{k,\bar{\delta}}({E}):=-2\big<\calS_{{N},k}({E}),\calP_{{U}{V}}\big>\phi\Big(\frac{\|{E}\|}{\bar{\delta}}\Big).
\end{equation*}
We now view $\varphi_{k,\bar{\delta}}(\cdot)$ as a function on $\RR^{(m_1+m_2)\times (m_1+m_2)}$ and 
abuse the notation here such that ${E}$ is viewed as a point in $\mathbb{R}^{(m_1+ m_2)\times (m_1+ m_2)}$.
\begin{lemma}\label{lemma:varphi-delta_n}
Under the conditions in Lemma~\ref{HSlemma},
for any ${E},{E}'\in\mathbb{R}^{(m_1+m_2)\times (m_1+m_2)}$, the following bounds hold ,
\begin{equation*}
 \big|\varphi_{k,\bar{\delta}}({E})-\varphi_{k,\bar{\delta}}({E}')\big|\leq 8r\cdot\frac{k+9}{\lambda_r}\Big(\frac{9\bar\delta}{2\lambda_r}\Big)^{k-1}\|{E}-{E}'\|\quad \forall k\geq 3
\end{equation*}
and
\begin{equation*}
 \big|\varphi_{\bar{\delta}}({E})-\varphi_{\bar{\delta}}({E}')\big|\leq \frac{C_5r}{\lambda_r}\Big(\frac{9\bar\delta}{2\lambda_r}\Big)^2\|{E}-{E}'\|
\end{equation*}
for an absolute constant $C_5>0$.
In other words,  $\varphi_{\bar{\delta}}(\cdot)$ and $\varphi_{k,\bar{\delta}}(\cdot)$ are both Lipschitz functions.
\end{lemma}

According to Lemma~\ref{lemma:varphi-delta_n}, we can write
\begin{gather*}
 \Big| \big(\|\calP_{\hat U\hat V}-\calP_{UV}\|_{\rm F}^2-\|\calL_N(E)\|_{\rm F}^2\big)\phi\Big(\frac{\|{E}\|}{\bar{\delta}}\Big)-\mathbb{E}  \big(\|\calP_{\hat U\hat V}-\calP_{UV}\|_{\rm F}^2-\|\calL_N(E)\|_{\rm F}^2\big)\phi\Big(\frac{\|{E}\|}{\bar{\delta}}\Big)\Big|\\
=\Big|\varphi_{\bar{\delta}}({E})-\mathbb{E}\varphi_{\bar{\delta}}({E})\Big|
\end{gather*}
where the function $\varphi_{\bar{\delta}}(\cdot)$ is a Lipschitz function with respect ${E}$ with constant $C_{6}\frac{r}{\lambda_r}\big(\frac{9\bar{\delta}}{2\lambda_r}\big)^2$. Since ${E}$ is a function of $\bDelta,\{\xi_i\}_{i=n+1}^{2n}, \{{X}_i\}_{i=n+1}^{2n}$, we will apply the Gaussian concentration inequality (Lemma~\ref{lem:gaussian_con}).  To this end,
recall  that ${E}={E}_1+{E}_2$ with
$$
{E}_1=\mathfrak{D}({Z}_1)=\frac{1}{n}\sum_{i=n+1}^{2n}\xi_i\mathfrak{D}({X}_i)
$$
and
$$
{E}_2=\calT({Z}_2)=\frac{1}{n}\sum_{i=n+1}^{2n}\langle\bDelta, {X}_i\rangle\mathfrak{D}({X}_i)-\mathfrak{D}(\bDelta).
$$
We denote by ${\rm vec}(\bDelta)$ the vectorization of $\bDelta$ and $\calM(\cdot)$ the matricization of vectors such that $\calM({\rm vec}(\bDelta))=\bDelta$. We also denote
$\calP_{{\rm vec}(\bDelta)}$ the orthogonal projection onto ${\rm vec}(\bDelta)$. More exactly, we can write
$$
\calP_{{\rm vec}(\bDelta)}{\rm vec}({X}_i)={\rm vec}(\bDelta)\cdot \frac{\langle\bDelta,{X}_i \rangle}{\|\bDelta\|_{\rm F}^2}.
$$
For each $i$, we write 
$$
{\rm vec}({X}_i)=\calP_{{\rm vec}(\bDelta)}{\rm vec}({X}_i)+\calP_{{\rm vec}(\bDelta)}^{\perp}{\rm vec}({X}_i).
$$
Clearly, $\calP_{{\rm vec}(\bDelta)}{\rm vec}({X}_i)$ is independent with $\calP_{{\rm vec}(\bDelta)}^{\perp}{\rm vec}({X}_i)$. We view $\varphi_{\bar{\delta}}({E})$ as a function of 
$$
\{\xi_i\}_{i=n+1}^{2n},\quad \{\calP_{{\rm vec}(\bDelta)}{\rm vec}({X}_i)\}_{i=n+1}^{2n}\quad {\rm and}\quad \{\calP_{{\rm vec}(\bDelta)}^{\perp}{\rm vec}({X}_i)\}_{i=n+1}^{2n}
$$
which are mutually independent. 
Let $\big\{\calP_{{\rm vec}(\Delta)}^{\perp}{\rm vec}(X_i')\big\}_{i=n+1}^{2n}$ be independent copies of \newline
$\big\{\calP_{{\rm vec}(\Delta)}^{\perp}{\rm vec}(X_i)\big\}_{i=n+1}^{2n}$.  
Conditioned on $\{\xi_i\}_{i=n+1}^{2n}$, $\{\calP_{{\rm vec}(\bDelta)}{\rm vec}({X}_i)\}_{i=n+1}^{2n}$, we have
\begin{align*}
\|{E}-{E}'\|_{\rm F}\leq& \Big\|\frac{1}{n}\sum_{i=n+1}^{2n} \xi_i \mathfrak{D}\circ\calM\big(\calP^{\perp}_{{\rm vec}(\bDelta)}{\rm vec}({X}_i)-\calP^{\perp}_{{\rm vec}(\bDelta)}{\rm vec}({X}_i')\big)\Big\|_{\rm F}\\
&+\Big\|\frac{1}{n} \sum_{i=n+1}^{2n}\langle\bDelta, {X}_i\rangle\mathfrak{D}\circ\calM\big(\calP_{{\rm vec}(\bDelta)}^{\perp}{\rm vec}({X}_i)-\calP_{{\rm vec}(\bDelta)}^{\perp}{\rm vec}({X}_i')\big)\Big\|_{\rm F}\\
\leq& 2*n^{-1}\Big(\sum_{i=n+1}^{2n} \xi_i^2\Big)^{1/2}\Big(\sum_{i=n+1}^{2n}\Big\|\calP_{{\rm vec}(\bDelta)}^{\perp}{\rm vec}({X}_i)-\calP_{{\rm vec}(\bDelta)}^{\perp}{\rm vec}({X}_i')\Big\|_{\rm F}^2\Big)^{1/2}\\
&+2*n^{-1}\Big(\sum_{i=n+1}^{2n} \langle\bDelta,{X}_i \rangle^2\Big)^{1/2}\Big(\sum_{i=n+1}^{2n}\Big\|\calP_{{\rm vec}(\bDelta)}^{\perp}{\rm vec}({X}_i)-\calP_{{\rm vec}(\bDelta)}^{\perp}{\rm vec}({X}_i')\Big\|_{\rm F}^2\Big)^{1/2}
\end{align*}
implying that $\|{E}\|_{\rm F}$ is a Lipschitz function with respect to $\big\{\calP_{{\rm vec}(\bDelta)}^{\perp}{\rm vec}({X}_i)\big\}_{i=n+1}^{2n}$ with constant $\frac{2}{n}\Big[\big(\sum_i\xi_i^2\big)^{1/2}+\big(\sum_i\langle \bDelta,{X}_i\rangle^2\big)^{1/2}\Big]$. % Note that conditional on $\{\calP_{{\rm vec}(\Delta)}X_i\}_i$, $\{\langle\Delta,X_i \rangle\}_i$ are also fixed. 
Therefore, by Lemma~\ref{lemma:varphi-delta_n}, we get
\begin{align*}
\big|\varphi_{\bar{\delta}}({E})-\varphi_{\bar{\delta}}({E}') \big|\leq& \frac{C_6r}{n\lambda_r}\Big(\frac{9\bar{\delta}}{2\lambda_r}\Big)^2\Big[\big(\sum_i\xi_i^2\big)^{1/2}+\big(\sum_i\langle \bDelta,{X}_i\rangle^2\big)^{1/2}\Big]\\
&\times\Big(\sum_{i=n+1}^{2n}\Big\|\calP_{{\rm vec}(\bDelta)}^{\perp}{\rm vec}({X}_i)-\calP_{{\rm vec}(\bDelta)}^{\perp}{\rm vec}({X}_i')\Big\|_{\rm F}^2\Big)^{1/2}.
\end{align*}
By the Gaussian isoperimetric inequality (Lemma~\ref{lem:gaussian_con}), conditioned on $\{\xi_i\}_i$ and $\{\calP_{{\rm vec}(\bDelta)}{\rm vec}({X}_i)\}_i$, 
with probability at least $1-e^{-t}$ for all $t\geq 1$, we get
\begin{align*}
 \big|\varphi_{\bar{\delta}}({E})-\mathbb{E}_{\{\calP^{\perp}_{{\rm vec}(\bDelta)}{X_i}\}_i}\big[\varphi_{\bar{\delta}}({E})&\big]\big|\leq  \frac{C_7rt^{1/2}}{n\lambda_r}\Big(\frac{9\bar{\delta}}{2\lambda_r}\Big)^2
\cdot \Big[\big(\sum_i\xi_i^2\big)^{1/2}+\big(\sum_i\langle \bDelta,{X}_i\rangle^2\big)^{1/2}\Big].
\end{align*}
Meanwhile, by the concentration inequality of the sum of exponential random variables (\cite{vershynin2010introduction}), with probability at least $1-e^{-n}$,
$$
\big(\sum_i\xi_i^2\big)^{1/2}+\big(\sum_i\langle \bDelta,{X}_i\rangle^2\big)^{1/2}\leq C_1n^{1/2}\big(\sigma_{\xi}+\|\bDelta\|_{\rm F}\big)\leq C_2\sigma_\xi\cdot n^{1/2}
$$
where the last inequality is due to Proposition~\ref{prop:nuclear_pen}. 
Therefore, with probability at least $1-e^{-t}-e^{-n}$,
\begin{align*}
 \Big|\varphi_{\bar{\delta}}({E})-\mathbb{E}_{\{\calP^{\perp}_{{\rm vec}(\bDelta)}{X_i}\}_i}\big[\varphi_{\bar{\delta}}({E})\big]\Big|
 \leq  C_7t^{1/2}\cdot\frac{r\sigma_\xi }{n^{1/2}\lambda_r}\Big(\frac{9\bar{\delta}}{2\lambda_r}\Big)^2.
\end{align*}
Next, we prove the bound of $\big|\mathbb{E}_{\{\calP^{\perp}_{{\rm vec}(\bDelta)}{X}\}}\big[\varphi_{\bar{\delta}}({E})\big]-\EE\big[\varphi_{\bar{\delta}}({E})\big]\big|$. We apply the following lemma whose proof is postponed to the appendix.
\begin{lemma}\label{lem:Eperp-E}
Under the assumptions of Lemma~\ref{lemma:varphi-delta_n} and $n\geq C_1r\bar{m}$ for some large enough absolute constant $C_1>0$, with probability at least $1-2e^{-t}-c_1e^{-c_2\bar{m}}$ for all $t\in[1,n]$, we have
\begin{align*}
 \big|\mathbb{E}_{\{\calP^{\perp}_{{\rm vec}(\bDelta)}{X_i}\}_i}&\big[\varphi_{\bar{\delta}}({E})\big]-\EE\big[\varphi_{\bar{\delta}}({E})\big]\big|
 \leq C_8\frac{r\sigma_\xi}{\lambda_r}\Big(\frac{9\bar{\delta}}{2\lambda_r}\Big)^2\cdot\frac{(t+\log\bar{m})^{1/2}}{n^{1/2}}
\end{align*}
for some absolute constant $c_1,c_2,C_8>0$.
\end{lemma}
We conclude that for all $t\in[1,n]$,  with probability at least $1-3e^{-t}-e^{-n}-c_1e^{-c_2\bar{m}}$,
\begin{align*}
 \big|\varphi_{\bar{\delta}}({E})-\EE&\varphi_{\bar{\delta}}({E})\big|
 \leq C_7\frac{r\sigma_\xi}{\lambda_r}\Big(\frac{9\bar{\delta}}{2\lambda_r}\Big)^2\cdot \frac{(t+\log\bar{m})^{1/2}}{n^{1/2}}.
\end{align*}
Since $\bar\delta=2\EE\|{E}\|$, on event $\calE_2:=\big\{\|{E}\|\leq 2\EE\|{E}\|\big\}$ with $\PP\big(\calE_2\big)\geq 1-e^{-n}-c_1e^{-c_2\bar{m}}$, 
\begin{align*}
\varphi_{\bar{\delta}}({E})=\|\calP_{\hat U\hat V}-\calP_{UV}\|_{\rm F}^2-\|\calL_N(E)\|_{\rm F}^2
\end{align*}
which concludes the proof.

\subsubsection{Proof of Lemma~\ref{lem:PN-E-CN}}

Recall that 
$$
{Z}_1=\frac{1}{n}\sum_{i=n+1}^{2n}\xi_i{X}_i.
$$
Conditional on $\{\xi_i\}_{i=n+1}^{2n}$, ${Z}_1$ has the same distribution as
$$
{Z}_1\stackrel{{\rm d}}{=} \tau\cdot{Z}:=\frac{1}{n}\Big(\sum_{i=n+1}^{2n}\xi_i^2\Big)^{1/2}\cdot {Z}
$$
where ${Z}\in\RR^{m_1\times m_2}$ has i.i.d. standard Gaussian entries. Then,
$$
\|\calP_{{U}{V}}^{\perp}{E}_1\calC_{{U}{V}}\|_{\rm F}^2\stackrel{{\rm d}}{=}\tau^2\cdot\Big(\|{U}_{\perp}{U}_{\perp}^{\top}{Z}{V}{\Lambda}^{-1}{U}^{\top}\|_{\rm F}^2+\|{V}_{\perp}{V}_{\perp}^{\top}{Z}^{\top}{U}{\Lambda}^{-1}{V}^{\top}\|_{\rm F}^2\Big).
$$
Denote by ${z}_1,\ldots,{z}_{m_2}$ the columns of ${Z}$, i.e., ${z}_j\in\calN({0},{I}_{m_1})$ are i.i.d. standard Gaussian vectors. We can write
\begin{gather*}
{U}_{\perp}{U}_{\perp}^{\top}{Z}{V}{\Lambda}^{-1}{U}^{\top}=\sum_{j=1}^{m_2}\big({U}_{\perp}{U}_{\perp}^{\top}{z}_j\big)\otimes\big({U}{\Lambda}^{-1}{V}^{\top}\be_j\big) 
\end{gather*}
where $\{\be_j\}_{j=1}^{m_2}$ denotes the standard basis vectors in $\RR^{m_2}$. In a similar fashion, write
\begin{gather*}
{V}_{\perp}{V}_{\perp}^{\top}{Z}^{\top}{U}{\Lambda}^{-1}{V}^{\top}=\sum_{j=1}^{m_2}\big({V}_{\perp}{V}_{\perp}^{\top}\be_j\big)\otimes\big({V}{\Lambda}^{-1}{U}^{\top}{z}_j\big).
\end{gather*}
We claim that ${U}_{\perp}{U}_{\perp}^{\top}{z}_j$ is independent with ${V}{\Lambda}^{-1}{U}^{\top}{z}_j$. Indeed, their correlation 
$$
\EE\big({U}_{\perp}{U}_{\perp}^{\top}{z}_j\big)\otimes\big({V}{\Lambda}^{-1}{U}^{\top}{z}_j\big)={U}_{\perp}{U}_{\perp}^{\top}{U}{\Lambda}^{-1}{V}^{\top}={0}.
$$
Since both vectors are Gaussian, we conclude that ${U}_{\perp}{U}_{\perp}^{\top}{z}_j$ is independent with ${V}{\Lambda}^{-1}{U}^{\top}{z}_j$ for all $1\leq j\leq m_2$. Therefore,
$$
\big\|{U}_{\perp}{U}_{\perp}^{\top}{Z}{V}\Lambda^{-1}{U}^{\top}\big\|_{\rm F}^2\textrm{ is independent of } \big\|{V}_{\perp}{V}_{\perp}^{\top}{Z}^{\top}{U}{\Lambda}^{-1}{V}^{\top}\big\|_{\rm F}^2.
$$

\paragraph*{Claim 1} Let $\bar {z}_k\in\RR^{m_1-r}$ be i.i.d. standard Gaussian vector independent of ${Z}$ for all $k=1,\ldots,r$. Then, we claim that
$$
\sum_{j=1}^{m_2}\big({U}_{\perp}{U}_{\perp}^{\top}{z}_j\big)\otimes\big({U}{\Lambda}^{-1}{V}^{\top}\be_j\big) \stackrel{{\rm d}}{=}\sum_{k=1}^{r}\big({U}_{\perp}\bar {z}_k\big)\otimes\big(\lambda_k^{-1}{u}_k\big)
$$
where $\{{u}_1,\ldots,{u}_{r}\}$ are the columns of ${U}$. To prove the claim, it suffices to check their covariance. To this end, define the following multilinear mapping:
$$
\calK({u}_1\otimes {u}_2\otimes {u}_3\otimes {u}_4)={u}_1\otimes {u}_3\otimes {u}_2\otimes {u}_4,\quad \forall {u}_1,{u}_2,{u}_3,{u}_4\in\RR^{m},
$$
a technique introduced in \cite{koltchinskii2015normal}. Then, we have
\begin{align*}
{\rm Cov}\big({U}_{\perp}{U}_{\perp}&{Z}{V}{\Lambda}^{-1}{U}^{\top}\big)=\EE \big({U}_{\perp}{U}_{\perp}{Z}{V}{\Lambda}^{-1}{U}^{\top}\big)\otimes \big({U}_{\perp}{U}_{\perp}{Z}{V}{\Lambda}^{-1}{U}^{\top}\big)\\
=&\EE \sum_{j=1}^{m_2}\big({U}_{\perp}{U}_{\perp}^{\top}{z}_j\big)\otimes\big({U}{\Lambda}^{-1}{V}^{\top}\be_j\big)\otimes \big({U}_{\perp}{U}_{\perp}^{\top}{z}_j\big)\otimes\big({U}{\Lambda}^{-1}{V}^{\top}\be_j\big)\\
=&\EE \sum_{j=1}^{m_2}\calK\Big(\big({U}_{\perp}{U}_{\perp}^{\top}{z}_j\big)\otimes \big({U}_{\perp}{U}_{\perp}^{\top}{z}_j\big)\otimes\big({U}{\Lambda}^{-1}{V}^{\top}\be_j\big)\otimes\big({U}{\Lambda}^{-1}{V}^{\top}\be_j\big)\Big)\\
=&\sum_{j=1}^{m_2}\calK\Big({U}_{\perp}{U}_{\perp}\otimes \big({U}{\Lambda}^{-1}{V}^{\top}(\be_j\otimes \be_j){V}{\Lambda}^{-1}{U}^{\top}\big)\Big)
=\calK\Big({U}_{\perp}{U}_{\perp}\otimes\big({U}{\Lambda}^{-2}{U}^{\top}\big)\Big).
\end{align*}
Similarly, we have
\begin{align*}
{\rm Cov}\Big(\sum_{k=1}^{r}\big({U}_{\perp}\bar {z}_k\big)&\otimes\big(\lambda_k^{-1}{u}_k\big)\Big)\\
=&\EE \sum_{k=1}^{r}\Big(\big({U}_{\perp}\bar {z}_k\big)\otimes\big(\lambda_k^{-1}{u}_k\big)\otimes\big({U}_{\perp}\bar {z}_k\big)\otimes\big(\lambda_k^{-1}{u}_k\big)\Big)\\
=&\EE \sum_{k=1}^{r}\calK\Big(\big({U}_{\perp}\bar {z}_k\big)\otimes\big({U}_{\perp}\bar {z}_k\big)\otimes\big(\lambda_k^{-1}{u}_k\big)\otimes\big(\lambda_k^{-1}{u}_k\big)\Big)\\
=&\sum_{k=1}^{r}\calK\Big({U}_{\perp}{U}_{\perp}^{\top}\otimes\big(\lambda_k^{-1}{u}_k\big)\otimes\big(\lambda_k^{-1}{u}_k\big)\Big)=\calK\Big({U}_{\perp}{U}_{\perp}\otimes\big({U}{\Lambda}^{-2}{U}^{\top}\big)\Big).
\end{align*}
which proves the claim. It implies that
\begin{align*}
\big\|{U}_{\perp}{U}_{\perp}{Z}{V}{\Lambda}^{-1}{U}^{\top}\big\|_{\rm F}^2\stackrel{{\rm d}}{=}&\big\|\sum_{k=1}^{r}\big({U}_{\perp}\bar {z}_k\big)\otimes\big(\lambda_k^{-1}{u}_k\big)\big\|_{\rm F}^2
=\sum_{k=1}^{r}\|{U}_{\perp}\bar{{z}}_k\|_{\ell_2}^2\lambda_k^{-2}
\end{align*}
where the last equality is due to the orthogonality of $\{{u}_k\}_{k=1}^{r}$. Clearly, $\|{U}_{\perp}\bar {z}_k\|_{\ell_2}^2$ has a Chi-squared distribution with degrees of freedom $m_1-r$. Therefore, we get that
$$
\big\|{U}_{\perp}{U}_{\perp}{Z}{V}{\Lambda}^{-1}{U}^{\top}\big\|_{\rm F}^2\stackrel{{\rm d}}{=}\sum_{k=1}^{r} \frac{z_{1,k}^2}{\lambda_{k}^2}
$$
where $z_{1,k}^2\sim \calX^2(m_1-r)$ are i.i.d. for $k=1,\ldots,r$.

\paragraph*{Claim 2} Let $\bar {z}_k\in \RR^{m_2-r}$ be i.i.d. standard Gaussian vector independent of ${Z}$ for all $k=1,\ldots,r$. We claim that
$$
{V}_{\perp}{V}_{\perp}^{\top}{Z}^{\top}{U}{\Lambda}^{-1}{V}^{\top}\stackrel{d}{=}\sum_{k=1}^{r}({V}_{\perp}\bar {z}_k)\otimes (\lambda_k^{-1}{v}_k)
$$
where $\{{v}_1,\ldots,{v}_{r}\}$ are the columns of ${V}$. Indeed, if we denote by $\check{{z}}_j^{\top}, 1\leq j\leq m_1$ the rows of ${Z}$. Then, $\check{{z}}_j\sim\calN({0},{I}_{m_2})$ are i.i.d. for all $1\leq j\leq m_1$. We write
\begin{gather*}
{V}_{\perp}{V}_{\perp}^{\top}{Z}^{\top}{U}{\Lambda}^{-1}{V}^{\top}=\sum_{j=1}^{m_1}({V}_{\perp}{V}_{\perp}^{\top}\check{{z}}_j)\otimes ({V}{\Lambda}^{-1}{U}^{\top}\be_j)
\end{gather*}
where $\{\be_1,\ldots,\be_{m_1}\}$ denotes the standard basis vectors in $\RR^{m_1}$. It is straightforward to check that
\begin{gather*}
{\rm Cov}\big({V}_{\perp}{V}_{\perp}^{\top}{Z}^{\top}{U}{\Lambda}^{-1}{V}^{\top}\big)=\calK\Big(({V}_{\perp}{V}_{\perp}^{\top})\otimes {V}{\Lambda}^{-2}{V}^{\top}\Big).
\end{gather*}
Similarly, we obtain
$$
{\rm Cov}\Big(\sum_{k=1}^{r}({V}_{\perp}\bar {z}_k)\otimes (\lambda_k^{-1}{v}_k)\Big)=\calK\Big(({V}_{\perp}{V}_{\perp}^{\top})\otimes {V}{\Lambda}^{-2}{V}^{\top}\Big)
$$
which proves the claim. Thus, we get
\begin{align*}
\big\|{V}_{\perp}{V}_{\perp}^{\top}{Z}^{\top}{U}{\Lambda}^{-1}{V}^{\top}\big\|_{\rm F}^2\stackrel{{\rm d}}{=}&\Big\|\sum_{k=1}^{r}({V}_{\perp}\bar {z}_k)\otimes (\lambda_k^{-1}{v}_k)\Big\|_{\rm F}^2
=\sum_{k=1}^{r}\|{V}_{\perp}\bar{{z}}_k\|_{\ell_2}^2\lambda_k^{-2}\stackrel{d}{=}\sum_{k=1}^{r}\frac{z_{2,k}^2}{\lambda_k^2}
\end{align*}
where $z_{2,k}^2$ are i.i.d. and $z_{2,k}^2\sim\calX^2(m_2-r)$. 

\paragraph*{Finalize the first claim of Lemma~\ref{lem:PN-E-CN}} By Claim~1 and Claim~2, we conclude that 
\begin{gather*}
\|\calP_{{U}{V}}^{\perp}{E}_1\calC_{{U}{V}}\|_{\rm F}^2\stackrel{{\rm d}}{=}\tau^2\cdot \sum_{k=1}^{r}\frac{z_k^2}{\lambda_k^2}
\end{gather*}
where $\{z_k^2\}_{k=1}^{r}$ are i.i.d. Chi-squared random variables with degrees of freedom $m_{\star}=m_1+m_2-2r$.
\paragraph*{Proof of second claim of Lemma~\ref{lem:PN-E-CN}}
 Recall from above that  $\|\calP_{{U}{V}}{E}_1\calC_{{U}{V}}\|_{\rm F}^2$ is 
a sum of sub-exponential random variables. By the standard concentration inequality for the sum of sub-exponential random variables (e.g. \cite[Proposition~5.6]{vershynin2010introduction}), with probability at least $1-e^{-t}$ for all $t\geq \log 2$,
\begin{align*}
\big| \|\calP_{{U}{V}}{E}_1\calC_{{U}{V}}\|_{\rm F}^2-&\EE \|\calP_{{U}{V}}{E}_1\calC_{{U}{V}}\|_{\rm F}^2\big|
\leq C\tau^2\cdot\max\Big\{ \|{\Lambda}^{-2}\|_{\rm F}m_{\star}^{1/2}t^{1/2}, \frac{t}{\lambda_r^2} \Big\}
\end{align*}
which concludes the proof.

\subsubsection{Proof of Lemma~\ref{lemma:PN-E2-CN-frob}}
Recall that
$$
\calP_{{U}{V}}^{\perp}{E}_2\calC_{{U}{V}}=\frac{1}{n}\sum_{i=n+1}^{2n} \langle\bDelta,{X}_i \rangle\calP_{{U}{V}}^{\perp}{X}_i\calC_{{U}{V}}-\calP_{{U}{V}}^{\perp}\bDelta\calC_{{U}{V}}
$$
where $\EE\langle\bDelta,{X}_i \rangle\calP_{{U}{V}}^{\perp}{X}_i\calC_{{U}{V}}=\calP_{{U}{V}}^{\perp}\bDelta\calC_{{U}{V}}$. By the independence between $X_i$ and $X_j$, we have
\begin{align*}
\EE\big\|\calP_{{U}{V}}^{\perp}{E}_2\calC_{{U}{V}}\big\|_{\rm F}^2=&\frac{1}{n^2}\sum_{i=n+1}^{2n}\EE\big\|\langle\bDelta,{X}_i \rangle\calP_{{U}{V}}^{\perp}{X}_i\calC_{{U}{V}}-\calP_{{U}{V}}^{\perp}\bDelta\calC_{{U}{V}}\big\|_{\rm F}^2\\
=&\frac{1}{n^2}\sum_{i=n+1}^{2n}\big(\EE\|\langle\bDelta,{X}_i \rangle\calP_{{U}{V}}^{\perp}{X}_i\calC_{{U}{V}}\|_{\rm F}^2-\|\calP_{{U}{V}}^{\perp}\bDelta\calC_{{U}{V}}\|_{\rm F}^2\big)\\
=& \frac{1}{n}\EE\langle\bDelta,{X} \rangle^2\|\calP_{{U}{V}}^{\perp}{X}\calC_{{U}{V}}\|_{\rm F}^2-\frac{1}{n}\|\calP_{UV}^{\perp}\Delta\calC_{UV}\|_{\rm F}^2.
\end{align*}
We then write
$$
\calP_{{U}{V}}^{\perp}{X}\calC_{{U}{V}}=\calP_{{U}{V}}^{\perp}\calM\big(\calP_{{\rm vec}(\bDelta)}{\rm vec}({X})\big)\calC_{{U}{V}}+\calP_{{U}{V}}^{\perp}\calM\big(\calP_{{\rm vec}(\bDelta)}^{\perp}{\rm vec}({X})\big)\calC_{{U}{V}}.
$$
Since $\calP_{{\rm vec}(\bDelta)}{\rm vec}({X})$ is independent with $\calP_{{\rm vec}(\bDelta)}^{\perp}{\rm vec}({X})$, we obtain
\begin{align*}
\EE\langle\bDelta,{X} \rangle^2\|\calP_{{U}{V}}^{\perp}{X}\calC_{{U}{V}}\|_{\rm F}^2=& \EE\langle\bDelta,{X} \rangle^2\big\|\calP_{{U}{V}}^{\perp}\calM\big(\calP_{{\rm vec}(\bDelta)}{\rm vec}({X})\big)\calC_{{U}{V}}\big\|_{\rm F}^2\\
&\quad+\EE\langle\Delta,X \rangle^2\big\|\calP_{{U}{V}}^{\perp}\calM\big(\calP_{{\rm vec}(\bDelta)}^{\perp}{\rm vec}({X})\big)\calC_{{U}{V}}\big\|_{\rm F}^2\\
=&\frac{\EE\langle\bDelta,{X} \rangle^4}{\|\Delta\|_{\rm F}^2}\cdot \|\calP_{UV}^{\perp}\Delta\calC_{UV}\|_{\rm F}^2+\|\bDelta\|_{\rm F}^2\EE\big\|\calP_{{U}{V}}^{\perp}\calM\big(\calP_{{\rm vec}(\bDelta)}^{\perp}{\rm vec}({X})\big)\calC_{{U}{V}}\big\|_{\rm F}^2\\
=&3\|\Delta\|_{\rm F}^2\cdot  \|\calP_{UV}^{\perp}\Delta\calC_{UV}\|_{\rm F}^2+\|\bDelta\|_{\rm F}^2\EE\big\|\calP_{{U}{V}}^{\perp}\calM\big(\calP_{{\rm vec}(\bDelta)}^{\perp}{\rm vec}({X})\big)\calC_{{U}{V}}\big\|_{\rm F}^2.
\end{align*}
By the proof of Lemma~\ref{lem:PN-E-CN}, we immediately conclude that 
$$
\EE\|\calP_{UV}^{\perp}E_2\calC_{UV}\|_{\rm F}^2\leq  \frac{5}{n}\cdot\frac{\|\bDelta\|_{\rm F}^2}{\lambda_r^2}+\frac{2m_{\star}}{n}\cdot\|{\Lambda}^{-1}\|_{\rm F}^2\|\bDelta\|_{\rm F}^2.
$$
Similarly, it is easy to show that 
$$
\EE\|\calP_{UV}^{\perp}E_2\calC_{UV}\|_{\rm F}^2\geq c_1\frac{m_{\star}}{n}\cdot \|{\Lambda}^{-1}\|_{\rm F}^2\|\bDelta\|_{\rm F}^2
$$
for some absolute constant $c_1>0$. 
It proves the second claim. To prove the first claim, we denote
$$
\bar\delta=2\EE\|{E}_2\|.
$$
Therefore, $\PP\Big(\|{E}_2\|\geq \bar\delta\Big)\leq e^{-n}+c_1e^{-c_2\bar{m}}$ for an absolute constant $c_1, c_2>0$. Let $\phi(\cdot)$ be the Lipschitz function defined in the proof of Lemma~\ref{HSlemma}. We define the function
$$
h_{6}\big(\{{X}_i\}_i\big)=\|\calP_{{U}{V}}^{\perp}{E}_2\calC_{{U}{V}}\|_{\rm F}^2\phi\Big(\frac{\|{E}_2\|}{\bar\delta}\Big).
$$
Since $\big\{\langle\bDelta,{X}_i \rangle\big\}_i$ and $\big\{\calP_{{\rm vec}(\bDelta)}{\rm vec}^{\perp}({X}_i) \big\}_i$ are independent, we view $h_{6}(\{{X}_i\}_i)$ as a function $h_{6}\big(\big\{\langle\bDelta,{X}_i\rangle\big\}_i, \big\{\calP_{{\rm vec}(\bDelta)}^{\perp}{\rm vec}({X}_i)\big\}_i\big)$. Conditional on $\big\{\langle\bDelta,{X}_i\rangle\big\}_i$, similarly as the proof of Lipschitz property in Lemma~\ref{lemma:varphi-delta_n}, we can show that
\begin{align*}
\Big|h_{6}\big(\big\{\langle\bDelta,{X}_i\rangle\big\}_i, &\big\{\calP_{{\rm vec}(\bDelta)}^{\perp}{\rm vec}({X}_i)\big\}_i\big)-h_{6}\big(\big\{\langle\bDelta,{X}_i\rangle\big\}_i, \big\{\calP_{{\rm vec}(\bDelta)}^{\perp}{\rm vec}({X}_i')\big\}_i\big) \Big|\\
\leq& C_3\|\calP_{{U}{V}}^{\perp}({E}_2-{E}_2')\calC_{{U}{V}}\|_{\rm F}\cdot \bar{\delta}\|\calC_{{U}{V}}\|_{\rm F}+C_4\bar{\delta}^2\|\calC_{{U}{V}}\|_{\rm F}^2\cdot \frac{\|{E}_2-{E}_2'\|}{\bar\delta}\\
\leq& C_3\bar{\delta}\|{\Lambda}^{-1}\|_{\rm F}^2\cdot \frac{\big(\sum_i\langle\bDelta,{X}_i \rangle^2\big)^{1/2}}{n}\cdot\Big(\sum_{i=n+1}^{2n}\Big\| \calP_{{\rm vec}(\bDelta)}^{\perp}({\rm vec}({X}_i-{X}_i'))\Big\|_{\rm F}^2\Big)^{1/2}.
\end{align*}
Therefore, conditioned on  $\big\{\langle\bDelta,{X}_i\rangle\big\}_i$, $h_{6}(\cdot)$ is a Lipschitz function. 
%\newpage
%For the first claim, observe that
%\begin{align*}
%\big|\|\calP_{{U}{V}}^{\perp}{E}_2&\calC_{{U}{V}}\|_{\rm F}-\|\calP_{{U}{V}}^{\perp}{E}_2'\calC_{{U}{V}}\|_{\rm F}\big|\leq \|\calP_{{U}{V}}^{\perp}({E}_2-{E}_2')\calC_{{U}{V}}\|_{\rm F}\\
%=&\Big\|\frac{1}{n}\sum_{i=1}^n \langle\bDelta,{X}_i \rangle\calP_{{U}{V}}^{\perp}{X}_i\calC_{{U}{V}}-\frac{1}{n}\sum_{i=1}^n \langle\bDelta,{X}_i' \rangle\calP_{{U}{V}}^{\perp}{X}_i'\calC_{{U}{V}} \Big\|_{\rm F}.
%\end{align*}
%Conditional on $\{\calP_{{\rm vec}(\bDelta)}{\rm vec}({X}_i)\}_i$, we observe
%\begin{align*}
%\big|\|\calP_{{U}{V}}^{\perp}{E}_2\calC_{{U}{V}}\|_{\rm F}&-\|\calP_{{U}{V}}^{\perp}{E}_2'\calC_{{U}{V}}\|_{\rm F}\big|\\
%\leq&\Big\|\frac{1}{n}\sum_{i=1}^n\langle\bDelta,{X}_i \rangle\calP_{{U}{V}}^{\perp}\calM\big(\calP_{{\rm vec}(\bDelta)}^{\perp}{\rm vec}({X}_i)-\calP_{{\rm vec}(\bDelta)}^{\perp}{\rm vec}({X}_i')\big)\calC_{{U}{V}} \Big\|_{\rm F}\\
%\leq& \frac{1}{n\lambda_r}\Big(\sum_{i=1}^n\langle\bDelta,{X}_i \rangle^2\Big)^{1/2}\Big(\sum_{i=1}^n\Big\|\calP_{{\rm vec}(\bDelta)}^{\perp}{\rm vec}({X}_i)-\calP_{{\rm vec}(\bDelta)}^{\perp}{\rm vec}({X}_i')\Big\|_{\rm F}^2\Big)^{1/2}
%\leq& \sqrt{\frac{2}{n}}\frac{\|\bDelta\|_{\rm F}}{\lambda_r}\Big(\sum_{i=1}^n\Big\|\calP_{{\rm vec}(\bDelta)}^{\perp}{\rm vec}({X}_i)-\calP_{{\rm vec}(\bDelta)}^{\perp}{\rm vec}(X_i')\Big\|_{\rm F}^2\Big)^{1/2}
%\end{align*}
 By Gaussian Isoperimetric inequality (Lemma~\ref{lem:gaussian_con}), with probability at least $1-e^{-t}$ for all $t\geq 1$,
\begin{align*}
\Big|\|\calP_{{U}{V}}^{\perp}{E}_2\calC_{{U}{V}}\|_{\rm F}^2\phi\Big(\frac{\|{E}_2\|}{\bar\delta}\Big)-&\EE_{\{\calP_{{\rm vec}(\bDelta)}^{\perp}{\rm vec}({X_i})\}_i}\|\calP_{{U}{V}}^{\perp}{E}_2\calC_{{U}{V}}\|_{\rm F}^2\phi\Big(\frac{\|{E}_2\|}{\bar\delta}\Big)\Big|\\
&\leq  C_3t^{1/2}\bar\delta\|{\Lambda}^{-1}\|_{\rm F}^2\cdot \frac{\big(\sum_i\langle\bDelta,{X}_i \rangle^2\big)^{1/2}}{n}.
\end{align*}
Meanwhile, with probability at least $1-e^{-n}$,
$$
 C_3t^{1/2}\bar\delta\|{\Lambda}^{-1}\|_{\rm F}^2\cdot \frac{\big(\sum_i\langle\bDelta,{X}_i \rangle^2\big)^{1/2}}{n}\leq \frac{C_4t^{1/2}}{n^{1/2}}\cdot \bar\delta \|{\Lambda}^{-1}\|_{\rm F}^2\|\bDelta\|_{\rm F}^2.
$$
Therefore, we get that with probability at least $1-e^{-t}-e^{-n}$ for $t\geq 1$,
\begin{align*}
\Big|\|\calP_{{U}{V}}^{\perp}{E}_2\calC_{{U}{V}}\|_{\rm F}^2\phi\Big(\frac{\|{E}_2\|}{\bar\delta}\Big)-&\EE_{\{\calP_{{\rm vec}(\bDelta)}^{\perp}{\rm vec}({X_i})\}_i}\|\calP_{{U}{V}}^{\perp}{E}_2\calC_{{U}{V}}\|_{\rm F}^2\phi\Big(\frac{\|{E}_2\|}{\bar\delta}\Big)\Big|\\
&\leq C_3t^{1/2}\|\Lambda^{-1}\|_{\rm F}^2\|\Delta\|_{\rm F}^2\cdot \frac{\bar{\delta}}{n^{1/2}}.
\end{align*}
We write ${E}_2={E}_{21}+{E}_{22}$ where
$$
{E}_{21}=\frac{1}{n}\sum_{i=n+1}^{2n}\langle\bDelta,{X}_i \rangle\calM\big(\calP_{{\rm vec}(\bDelta)}{\rm vec}({X}_i)\big)-\bDelta
$$
and
$$
{E}_{22}=\frac{1}{n}\sum_{i=n+1}^{2n}\langle\bDelta,{X}_i \rangle\calM\big(\calP_{{\rm vec}(\bDelta)}^{\perp}{\rm vec}({X}_i)\big).
$$
Then, we get
\begin{align*}
\Big|\EE\|&\calP_{{U}{V}}^{\perp}{E}_2\calC_{{U}{V}}\|_{\rm F}^2\phi\Big(\frac{\|{E}_2\|}{\bar\delta}\Big)-\EE_{\{\calP_{{\rm vec}(\bDelta)}^{\perp}{\rm vec}({X_i})\}_i}\|\calP_{{U}{V}}^{\perp}{E}_2\calC_{{U}{V}}\|_{\rm F}^2\phi\Big(\frac{\|{E}_2\|}{\bar\delta}\Big)\Big|\\
\leq& \Big|\EE\|\calP_{{U}{V}}^{\perp}{E}_2\calC_{{U}{V}}\|_{\rm F}^2\phi\Big(\frac{\|{E}_2\|}{\bar\delta}\Big)-\EE\|\calP_{{U}{V}}^{\perp}{E}_{22}\calC_{{U}{V}}\|_{\rm F}^2\phi\Big(\frac{\|{E}_2\|}{\bar\delta}\Big) \Big|\\
 + &\Big|\EE\|\calP_{{U}{V}}^{\perp}{E}_{22}\calC_{{U}{V}}\|_{\rm F}^2\phi\Big(\frac{\|{E}_2\|}{\bar\delta}\Big)-\EE_{\{\calP_{{\rm vec}(\bDelta)}^{\perp}{\rm vec}({X_i})\}_i}\|\calP_{{U}{V}}^{\perp}{E}_{22}\calC_{{U}{V}}\|_{\rm F}^2\phi\Big(\frac{\|{E}_2\|}{\bar\delta}\Big) \Big|\\
 &+\Big|\EE_{\{\calP_{{\rm vec}(\bDelta)}^{\perp}{\rm vec}({X_i})\}_i}\|\calP_{{U}{V}}^{\perp}{E}_{2}\calC_{{U}{V}}\|_{\rm F}^2\phi\Big(\frac{\|{E}_2\|}{\bar\delta}\Big)-\EE_{\{\calP_{{\rm vec}(\bDelta)}^{\perp}{\rm vec}({X_i})\}_i}\|\calP_{{U}{V}}^{\perp}{E}_{22}\calC_{{U}{V}}\|_{\rm F}^2\phi\Big(\frac{\|{E}_2\|}{\bar\delta}\Big) \Big|.
\end{align*}
Similarly, we can show that the function $\EE_{\{\calP_{{\rm vec}(\bDelta)}^{\perp}{\rm vec}({X_i})\}_i}\|\calP_{{U}{V}}^{\perp}{E}_{22}\calC_{{U}{V}}\|_{\rm F}^2\phi\Big(\frac{\|{E}_2\|}{\delta}\Big)$ is Lipschitz with respect to $\{\langle\bDelta,{X}_i \rangle\}_i$ with constant $C_1n^{-1}\bar{m}^{1/2}\bar\delta\cdot\|{\Lambda}^{-1}\|_{\rm F}^2$. Then,
with probability at least $1-2e^{-t}$,
\begin{align*}
\Big|\EE\|\calP_{{U}{V}}^{\perp}{E}_2&\calC_{{U}{V}}\|_{\rm F}^2\phi\Big(\frac{\|{E}_2\|}{\bar\delta}\Big)-\EE_{\{\calP_{{\rm vec}(\bDelta)}^{\perp}{\rm vec}({X_i})\}_i}\|\calP_{{U}{V}}^{\perp}{E}_2\calC_{{U}{V}}\|_{\rm F}^2\phi\Big(\frac{\|{E}_2\|}{\bar\delta}\Big)\Big|\\
\leq& C_1\bar\delta\|\bDelta\|_{\rm F}\|{\Lambda}^{-1}\|_{\rm F}^2\cdot\frac{t^{1/2}+\log^{1/2}{\bar m}}{n^{1/2}}+C_2\bar\delta\|{\Lambda}^{-1}\|_{\rm F}^2\|\bDelta\|_{\rm F}\frac{\bar{m}^{1/2}t^{1/2}}{n}\\
\leq&C_1\bar\delta\|\bDelta\|_{\rm F}\|{\Lambda}^{-1}\|_{\rm F}^2\cdot\frac{t^{1/2}+\log^{1/2}{\bar m}}{n^{1/2}}
\end{align*}
for some absolute constants $C_1,C_2>0$ where the last inequality is due to $n\geq \bar{m}$. Therefore, we conclude that with probability at least $1-3e^{-t}-e^{-n}$,
\begin{align*}
\Big|\|\calP_{{U}{V}}^{\perp}{E}_2\calC_{{U}{V}}\|_{\rm F}^2&\phi\Big(\frac{\|{E}_2\|}{\bar\delta}\Big)-\EE\|\calP_{{U}{V}}^{\perp}{E}_2\calC_{{U}{V}}\|_{\rm F}^2\phi\Big(\frac{\|{E}_2\|}{\bar\delta}\Big)\Big|\\
\leq& C_1\|\bDelta\|_{\rm F}^2\|{\Lambda}^{-1}\|_{\rm F}^2\frac{(t^{1/2}+\log^{1/2}{\bar m})}{n^{1/2}}%+C_2\|{\Lambda}^{-1}\|_{\rm F}^2\|\bDelta\|_{\rm F}^2\frac{\bar{m}t^{1/2}\log^{1/2}{\bar m}}{n^{3/2}}.
\end{align*}
Since $\PP\Big(\|{E}_2\|\geq \bar\delta\Big)\leq e^{-n}+c_1e^{-c_2\bar{m}}$, we obtain
\begin{align*}
\EE\|\calP_{{U}{V}}^{\perp}{E}_2\calC_{{U}{V}}\|_{\rm F}^2\bigg[1-\phi\Big(\frac{\|{E}_2\|}{\delta}\Big)\bigg]\leq& \EE^{1/2}\|{E}_2\|^4\|{\Lambda}^{-1}\|_{\rm F}^2\cdot \PP^{1/2}\Big(\|{E}_2\|\geq \delta\Big)\\
\leq& C_2\|{\Lambda}^{-1}\|_{\rm F}^2\|\bDelta\|_{\rm F}^2\cdot \frac{\bar{m}\log\bar{m}}{n}\cdot \big(e^{-c_1\bar{m}/2}+e^{-n/2}\big). 
\end{align*}
Together with Proposition~\ref{prop:nuclear_pen}, we get that with probability at least $1-3e^{-t}-c_1e^{-c_2\bar{m}}-2e^{-n}$ for $t\geq 1$, 
\begin{align*}
\Big|\|\calP_{{U}{V}}^{\perp}{E}_2\calC_{{U}{V}}\|_{\rm F}^2&-\EE\|\calP_{{U}{V}}^{\perp}{E}_2\calC_{{U}{V}}\|_{\rm F}^2\Big|
\leq C_1\sigma_\xi^2\|{\Lambda}^{-1}\|_{\rm F}^2\frac{r\bar{m}(t^{1/2}+\log^{1/2}{\bar m})}{n^{3/2}}. %+C_2\|{\Lambda}^{-1}\|_{\rm F}^2\|\bDelta\|_{\rm F}^2\frac{\bar{m}t^{1/2}\log^{1/2}{\bar m}}{n^{3/2}}.
\end{align*}

\subsubsection{Proof of Lemma~\ref{lem:PN-E1-CN-PN-E2-CN}}
We write
\begin{gather*}
\langle\calP_{{U}{V}}^{\perp}{E}_1\calC_{{U}{V}}, \calP_{{U}{V}}^{\perp}{E}_2\calC_{{U}{V}}\rangle=\frac{1}{n}\sum_{i=n+1}^{2n}\xi_i\langle\calP_{{U}{V}}^{\perp}{X}_i\calC_{{U}{V}}, {K}\rangle
\end{gather*}
where 
$$
{K}=\frac{1}{n}\sum_{i=n+1}^{2n}\big(\langle\bDelta, {X}_i\rangle\calP_{{U}{V}}^{\perp}{X}_i\calC_{{U}{V}}-\calP_{{U}{V}}^{\perp}\bDelta\calC_{{U}{V}}\big).
$$
 Conditional on $\{{X}_i\}_{i=n+1}^{2n}$, we get that with probability at least $1-e^{-t}$ for all $t\geq 1$,
\begin{gather*}
\big|\langle\calP_{{U}{V}}^{\perp}{E}_1\calC_{{U}{V}}, \calP_{{U}{V}}^{\perp}{E}_2\calC_{{U}{V}}\rangle\big|\leq C_2\sigma_\xi\frac{ t^{1/2}}{n}\Big(\sum_{i=n+1}^{2n}\langle \calP_{{U}{V}}^{\perp}{X}_i\calC_{{U}{V}}, {K}\rangle^2\Big)^{1/2}.
\end{gather*}
For each $n+1\leq i\leq 2n$, we have
\begin{align*}
\langle\calP_{{U}{V}}^{\perp}&{X}_i\calC_{{U}{V}}, {K} \rangle\\
=&\frac{1}{n}\sum_{j=n+1}^{2n}\Big(\langle \bDelta, {X}_j\rangle\langle\calP_{{U}{V}}^{\perp}{X}_i\calC_{{U}{V}}, \calP_{{U}{V}}^{\perp}{X}_j\calC_{{U}{V}} \rangle-\langle\calP_{{U}{V}}^{\perp}\bDelta\calC_{{U}{V}}, \calP_{{U}{V}}^{\perp}{X}_i\calC_{{U}{V}}\rangle\Big)\\
=&\frac{1}{n}\big(\langle \bDelta,{X}_i\rangle\|\calP_{{U}{V}}^{\perp}{X}_i\calC_{{U}{V}}\|_{\rm F}^2-\langle\calP_{{U}{V}}^{\perp}\bDelta\calC_{{U}{V}}, \calP_{{U}{V}}^{\perp}{X}_i\calC_{{U}{V}}\rangle\big)\\
+&\frac{1}{n}\sum_{j\neq i}\Big(\langle \bDelta, {X}_j\rangle\langle\calP_{{U}{V}}^{\perp}{X}_i\calC_{{U}{V}}, \calP_{{U}{V}}^{\perp}{X}_j\calC_{{U}{V}} \rangle-\langle\calP_{{U}{V}}^{\perp}\bDelta\calC_{{U}{V}}, \calP_{{U}{V}}^{\perp}{X}_i\calC_{{U}{V}}\rangle\Big).
\end{align*}
%By Lemma~\ref{lem:PN-E-CN}, with probability at least $1-e^{-t}$,
%\begin{gather*}
%\big|\frac{1}{n}\big(\langle \Delta,X_i\rangle\|\calP_{{U}{V}}^{\perp}X_i\calC_{{U}{V}}\|_{\rm F}^2-\langle\calP_{{U}{V}}^{\perp}\Delta\calC_{{U}{V}}, \calP_{{U}{V}}^{\perp}X_i\calC_{{U}{V}}\rangle\big)\big|\\
%\leq 
%\frac{C}{n}\|\Delta\|_{\rm F}\big((m_1+m_2-2r_M)\|\Sigma^{-1}\|_{\rm F}^2+\|\Sigma^{-2}\|_{\rm F}t^{1/2}(m_1+m_2-2r_M)^{1/2}\big)\\
%+\frac{C}{n}\|\Delta\|_{\rm F}\sigma_r^{-1}(M)\Big((m_1+m_2-2r_M)^{1/2}\|\Sigma^{-1}\|_{\rm F}+\|\Sigma^{-2}\|_{\rm F}^{1/2}t^{1/4}(m_1+m_2-2r_M)^{1/4}\Big).
%\end{gather*}
Conditioned on ${X}_i$, we apply the concentration inequality of the sum of sub-exponential random variables (\cite{vershynin2010introduction}) and obtain that with probability at least $1-e^{-t}$,
\begin{gather*}
\Big|\frac{1}{n}\sum_{j\neq i}\Big(\langle \bDelta, {X}_j\rangle\langle\calP_{{U}{V}}^{\perp}{X}_i\calC_{{U}{V}}, \calP_{{U}{V}}^{\perp}{X}_j\calC_{{U}{V}} \rangle-\langle\calP_{{U}{V}}^{\perp}\bDelta\calC_{{U}{V}}, \calP_{{U}{V}}^{\perp}{X}_i\calC_{{U}{V}}\rangle \Big)\Big|\\
\leq C_1\|\bDelta\|_{\rm F}\|\calP_{{U}{V}}^{\perp}{X}_i\calC_{{U}{V}}^2\|_{\rm F}\frac{t}{n^{1/2}}.
\end{gather*}
Conditioned on $X_i$, we get that with probability at least $1-e^{-t}$,
\begin{align*}
\big|\langle\calP_{{U}{V}}^{\perp}{X}_i\calC_{{U}{V}}, {K} \rangle \big|&\leq  C_1\|\bDelta\|_{\rm F}\|\calP_{{U}{V}}^{\perp}{X}_i\calC_{{U}{V}}^2\|_{\rm F}\frac{t}{n^{1/2}}+\frac{\|\bDelta\|_{\rm F}\|\calP_{{U}{V}}^{\perp}{X}_i\calC_{{U}{V}}^2\|_{\rm F}}{n}\\
&\quad+\frac{\big|\langle\bDelta,{X}_i \rangle\big|\|\calP_{{U}{V}}^{\perp}{X}_i\calC_{{U}{V}}\|_{\rm F}^2}{n},
\end{align*}
implying that with probability at least $1-ne^{-t}$,
\begin{gather*}
\sum_{i=n+1}^{2n}\big|\langle\calP_{{U}{V}}^{\perp}{X}_i\calC_{{U}{V}}, {K} \rangle \big|^2\leq C_1\|\bDelta\|_{\rm F}^2\frac{t^2}{n}\sum_{i=n+1}^{2n}\|\calP_{{U}{V}}^{\perp}{X}_i\calC_{{U}{V}}^2\|_{\rm F}^2\\
+C_2\frac{\|\bDelta\|_{\rm F}^2}{n^2}\sum_{i=n+1}^{2n}\|\calP_{{U}{V}}^{\perp}{X}_i\calC_{{U}{V}}^2\|_{\rm F}^2+C_2\frac{1}{n^2}\sum_{i=n+1}^{2n}\langle\bDelta,{X}_i \rangle^2\|\calP_{{U}{V}}^{\perp}{X}_i\calC_{{U}{V}}\|_{\rm F}^4.
\end{gather*}
By the proof of Lemma~\ref{lem:PN-E-CN}, with probability at least $1-ne^{-m_{\star}}$ for all $n+1\leq i\leq 2n$, 
$$
\|\calP_{{U}{V}}^{\perp}{X}_i\calC_{{U}{V}}\|_{\rm F}^2\leq C_1m_{\star}\|{\Lambda}^{-1}\|_{\rm F}^2
$$
and
$$
\|\calP_{{U}{V}}^{\perp}{X}_i\calC_{{U}{V}}^{2}\|_{\rm F}^2\leq C_1m_{\star}\|{\Lambda}^{-2}\|_{\rm F}^2
$$
where we used the fact $\|{\Lambda}^{-2}\|_{\rm F}\leq \|{\Lambda}^{-1}\|_{\rm F}^2$ and $\|{\Lambda}^{-4}\|\leq \|{\Lambda}^{-2}\|_{\rm F}^2$. 
Meanwhile, with probability at least $1-ne^{-t}$ for all $n+1\leq i\leq 2n$, $\langle\bDelta,{X}_i \rangle^2\leq \|\bDelta\|_{\rm F}^2t$.
Therefore, with probability at least $1-2ne^{-t}-ne^{-m_{\star}}$ for all $t\geq \log n$
\begin{align*}
\sum_{i=n+1}^{2n}\big|\langle\calP_{{U}{V}}^{\perp}{X}_i\calC_{{U}{V}}, {K} \rangle \big|^2\leq& C_1m_{\star}\Big(t^2+\frac{1}{n}\Big)\|\bDelta\|_{\rm F}^2\|{\Lambda}^{-2}\|_{\rm F}^2
+C_2tm_{\star}^2\frac{\|\bDelta\|_{\rm F}^2\|{\Lambda}^{-1}\|_{\rm F}^4}{n}.
\end{align*}
Then, we conclude with probability at least $1-(2n+1)e^{-t}-ne^{-m_{\star}}$ for $t\geq \log n$,
\begin{align*}
\big|\langle\calP_{{U}{V}}^{\perp}&{E}_1\calC_{{U}{V}}, \calP_{{U}{V}}^{\perp}{E}_2\calC_{{U}{V}}\rangle\big|\\
\leq& C_1\sigma_\xi\Big(t^{3/2}+\frac{t^{1/2}}{n^{1/2}}\Big)\frac{m_{\star}^{1/2}}{n}\|\bDelta\|_{\rm F}\|{\Lambda}^{-2}\|_{\rm F}+C_2\sigma_\xi t^{1/2}\frac{m_{\star}}{n^{3/2}}\|\bDelta\|_{\rm F}\|{\Lambda}^{-1}\|_{\rm F}^2\\
\leq&C_1\sigma_\xi\|\bDelta\|_{\rm F}\|{\Lambda}^{-2}\|_{\rm F}\cdot \frac{t^{3/2}m_{\star}^{1/2}}{n}+C_2\sigma_\xi t^{1/2}\frac{m_{\star}}{n^{3/2}}\|\bDelta\|_{\rm F}\|{\Lambda}^{-1}\|_{\rm F}^2
\end{align*}
where we used the fact $\|{\Lambda}^{-2}\|_{\rm F}\leq \|{\Lambda}^{-1}\|_{\rm F}^2$. Together with Proposition~\ref{prop:nuclear_pen}, we conclude that with probability at least $1-(2n+1)e^{-t}-ne^{-m_{\star}}-c_1e^{-c_2\bar{m}}$ for $t\geq \log n$,
\begin{align*}
\big|\langle\calP_{{U}{V}}^{\perp}&{E}_1\calC_{{U}{V}}, \calP_{{U}{V}}^{\perp}{E}_2\calC_{{U}{V}}\rangle\big|
\leq C_1t^{3/2}\sigma_\xi^2\|{\Lambda}^{-2}\|_{\rm F}\cdot \frac{r^{1/2}\bar{m}}{n^{3/2}}+C_2t\sigma_\xi^2 \|\Lambda^{-1}\|_{\rm F}^2\cdot \frac{r^{1/2}\bar{m}^{3/2}}{n^{2}}
\end{align*}
for absolute constants $C_1,C_2>0$.

\subsubsection{Proof of Theorem~\ref{thm:hatP-normal_approx}}
Denote by 
$$
\hat T=\frac{\|\calP_{\hat{U}\hat{V}}-\calP_{{U}{V}}\|_{\rm F}^2-\EE\|\calP_{\hat{U}\hat{V}}-\calP_{{U}{V}}\|_{\rm F}^2}{\sqrt{8}\sigma_\xi^2\|{\Lambda}^{-2}\|_{\rm F}\cdot \frac{m_{\star}^{1/2}}{n}}
$$
and write $\hat T=\hat T_0+\hat T_1$ where
$$
\hat T_0=\frac{\|\calP^{\perp}_{{U}{V}}{E}_1\calC_{{U}{V}}\|_{\rm F}^2-\EE\|\calP^{\perp}_{{U}{V}}{E}_1\calC_{{U}{V}}\|_{\rm F}^2 }{\sqrt{2}\sigma_\xi^2\|{\Lambda}^{-2}\|_{\rm F}\cdot \frac{m_{\star}^{1/2}}{n}}
$$
and
\begin{align*}
\hat T_1=&\frac{\big(\|\calP_{\hat{U}\hat{V}}-\calP_{{U}{V}}\|_{\rm F}^2-\|\calL_{N}({E})\|_{\rm F}^2\big)-\EE\big(\|\calP_{\hat{U}\hat{V}}-\calP_{{U}{V}}\|_{\rm F}^2-\|\calL_{N}({E})\|_{\rm F}^2\big)}{\sqrt{8}\sigma_\xi^2\|{\Lambda}^{-2}\|_{\rm F}\cdot \frac{m_{\star}^{1/2}}{n}}\\
&+\frac{\|\calP_{{U}{V}}^{\perp}{E}_2\calC_{{U}{V}}\|_{\rm F}^2-\EE\|\calP_{{U}{V}}^{\perp}{E}_2\calC_{{U}{V}}\|_{\rm F}^2 }{\sqrt{2}\sigma_\xi^2\|{\Lambda}^{-2}\|_{\rm F}\cdot \frac{m_{\star}^{1/2}}{n}}
+\frac{\sqrt{2}\big<\calP_{{U}{V}}^{\perp}{E}_1\calC_{{U}{V}},\calP_{{U}{V}}^{\perp}{E}_2\calC_{{U}{V}}\big>}{\sigma_\xi^2\|{\Lambda}^{-2}\|_{\rm F}\cdot \frac{m_{\star}^{1/2}}{n}}.
\end{align*}
By Lemma~\ref{HSlemma}, Lemma~\ref{lemma:PN-E2-CN-frob} and Lemma~\ref{lem:PN-E1-CN-PN-E2-CN} with $t= 2\log n$, we get that with probability at least $1-2e^{-n}-ne^{-c_1m_{\star}}-\frac{2n+6}{n^2}-c_1e^{-c_2\bar{m}}$,
\begin{align*}
\|\hat T_1\|\leq& C_1\frac{\sigma_\xi}{\|{\Lambda}^{-2}\|_{\rm F}\lambda_r^3}\cdot \frac{r\bar{m}^{1/2}\log^{1/2}n}{n^{1/2}}
+C_2\frac{\|{\Lambda}^{-1}\|_{\rm F}^2}{\|{\Lambda}^{-2}\|_{\rm F}}\cdot \frac{r^{1/2}\bar{m}^{1/2}\log n}{n^{1/2}}+C_3\frac{r^{1/2}\bar{m}^{1/2}\log^{3/2}n}{n^{1/2}}\\
\leq&C_1\frac{\sigma_\xi}{\lambda_r}\cdot \frac{r\bar{m}^{1/2}\log^{1/2}n}{n^{1/2}}+C_2\frac{r\bar{m}^{1/2}\log^{3/2}n}{n^{1/2}}\leq C_1(\beta\vee 1)\cdot \frac{r\bar{m}^{1/2}\log^{3/2}n}{n^{1/2}}
\end{align*}
where we used the fact $\|{\Lambda}^{-2}\|_{\rm F}\geq\lambda_r^{-2}$ and $\frac{\|{\Lambda}^{-1}\|_{\rm F}^2}{\|{\Lambda}^{-2}\|_{\rm F}}\leq r^{1/2}$ and $n\geq C_1r^2\bar{m}$.
By Lemma~\ref{lem:PN-E-CN}, we have
$$
\hat T_0\stackrel{{\rm d}}{=}\underbrace{\frac{\sum_{k=1}^r\lambda_k^{-2}\sum_{j_k=1}^{m_{\star}}(z_{k,j_k}^2-1)}{\sqrt{2}m_{\star}^{1/2}\|{\Lambda}^{-2}\|_{\rm F}}}_{\hat T_{00}}+\underbrace{\frac{\frac{\sum_{i=n+1}^{2n}(\xi_i^2-\sigma_\xi^2)}{n^2}\cdot \sum_{k=1}^r\lambda_k^{-2}\sum_{j_k=1}^{m_{\star}}z_{k,j_k}^2}{\sqrt{2}\sigma_{\xi}^2\|{\Lambda}^{-2}\|_{\rm F}\cdot \frac{m_{\star}^{1/2}}{n}}}_{\hat T_{01}}
$$
where $\{z_{k,j_k}\}_{k\in[r]}^{j_k\in[m_{\star}]}$ are i.i.d. standard normal random variables. Observe that $\{\xi_i\}_{i=n+1}^{2n}$ are independent with $\{z_{k,j_k}\}$. With probability at least $1-\frac{1}{n}-re^{-m_{\star}}$, we have
$$
|\hat T_{01}|\leq C_1\frac{\|{\Lambda}^{-1}\|_{\rm F}^2}{\|{\Lambda}^{-2}\|_{\rm F}}\cdot \frac{\bar{m}^{1/2}\log n}{n^{1/2}}\leq C_1\cdot \frac{r^{1/2}\bar{m}^{1/2}\log n}{n^{1/2}}.
$$
By Berry-Esseen theorem, for any $x\in \RR$,
\begin{align*}
\Big|\PP\big\{\hat T_{00}\leq x\big\}-\Phi(x)\Big|\leq \frac{\|{\Lambda}^{-3}\|_{\rm F}^2}{\|{\Lambda}^{-2}\|_{\rm F}^3}\cdot \frac{C_2}{\bar{m}^{1/2}}\leq \frac{C_2}{\bar{m}^{1/2}}
\end{align*}
where we used the facts $\|{\Lambda}^{-3}\|_{\rm F}^2\leq \|{\Lambda}^{-2}\|_{\rm F}^3$ and
$$
\EE \sum_{k=1}^r\lambda_k^{-4}\sum_{j_k=1}^{m_{\star}}(z_{k,j_k}^2-1)^2=2m_{\star}\|{\Lambda}^{-2}\|_{\rm F}^2
$$
and
$$
\EE \sum_{k=1}^r\lambda_k^{-6}\sum_{j_k=1}^{m_{\star}}|z_{k,j_k}^2-1|^3\leq C_1m_{\star}\|{\Lambda}^{-3}\|_{\rm F}^2.
$$
Now, recall that $\hat T=\hat T_{00}+ \hat T_{01}+\hat T_1$. Then, we get
\begin{align*}
\PP\Big(\hat T\leq x\Big)\leq& \PP\bigg(\hat T_{00}\leq x+C_1(\beta\vee 1)\cdot \frac{r\bar{m}^{1/2}\log^{3/2}n}{n^{1/2}}\bigg)\\
&\quad+2e^{-n}+c_1(n+r)e^{-c_2m_{\star}}+\frac{3n+6}{n^2}\\
\leq&\Phi\bigg(x+C_1(\beta\vee 1)\cdot \frac{r\bar{m}^{1/2}\log^{3/2}n}{n^{1/2}}\bigg)+c_1ne^{-c_2m_{\star}}+\frac{3n+6}{n^2}+\frac{C_3}{\bar{m}^{1/2}}\\
\leq&\Phi(x)+C_1(\beta\vee 1)\cdot \frac{r\bar{m}^{1/2}\log^{3/2}n}{n^{1/2}}+c_1ne^{-c_2m_{\star}}+\frac{3n+6}{n^2}+\frac{C_3}{\bar{m}^{1/2}}
\end{align*}
where the last inequality is due to the Lipschitz property of function $\Phi(\cdot)$. Similarly, we have
\begin{align*}
\PP\Big(\hat T\leq x\Big)\geq& \PP\bigg(\hat T_{00}\leq x-C_1(\beta\vee 1)\cdot \frac{r\bar{m}^{1/2}\log^{3/2}n}{n^{1/2}}\bigg)\\
&\quad-2e^{-n}-c_1(n+r)e^{-c_2m_{\star}}-\frac{3n+6}{n^2}\\
\geq&\Phi\bigg(x-C_1(\beta\vee 1)\cdot \frac{r\bar{m}^{1/2}\log^{3/2}n}{n^{1/2}}\bigg)-c_1ne^{-c_2m_{\star}}-\frac{3n+6}{n^2}- \frac{C_3}{\bar{m}^{1/2}}\\
\geq&\Phi(x)-C_1(\beta\vee 1)\cdot \frac{r\bar{m}^{1/2}\log^{3/2}n}{n^{1/2}}-c_1ne^{-c_2m_{\star}}-\frac{3n+6}{n^2}- \frac{C_3}{\bar{m}^{1/2}}
\end{align*}
By combining the above two inequalities, we obtain the claimed bound.

\subsubsection{Proof of Lemma~\ref{lem:V_n-B_n}}
By the definition of $\hat\sigma_\xi^2$, we can write
\begin{align*}
\hat\sigma_\xi^2:=& \frac{1}{n}\sum_{i=n+1}^{2n}\xi_i^2+\frac{1}{n}\sum_{i=n+1}^{2n}\langle\bDelta,{X}_i \rangle^2+\frac{2}{n}\sum_{i=n+1}^{2n}\xi_i\langle\bDelta, {X}_i\rangle.
\end{align*}
By the concentration inequality of the sum of sub-exponential random variables (see \cite{vershynin2010introduction}), we conclude that with probability at least $1-\frac{2}{n^2}$, 
\begin{align*}
\big|\hat\sigma_\xi^2-(\sigma_\xi^2+\|\bDelta\|_{\rm F}^2) \big|\leq C_1\sigma_\xi^2 \cdot \frac{\log n}{n^{1/2}}
\end{align*}
for some absolute constant $C_1>0$, where we also used the fact $\|\Delta\|_{\rm F}^2=O_P\Big(\sigma_\xi\cdot \frac{r\bar{m} }{n}\Big)$ and $n\gg r\bar{m}$. To prove the the concentration bound for $\hat B_n$ and $\hat V_n$, we apply the results from random matrix theory \cite{ding2017high}.   Then, we can immediate show that the following bounds hold with probability at least $1-\frac{1}{\bar{m}^2}$ for all $1\leq j\leq r$,  
\begin{align*}
\Big|\hat\lambda_j^2-\lambda_j^2-\sigma_\xi^2\cdot \frac{m_1+m_2}{n} \Big|\leq &C_2\beta^2\cdot \frac{\sigma_\xi^2\bar{m}^2}{n^2}+C_2\sigma_\xi^{3/2}\lambda_j^{1/2}\cdot \frac{\bar{m}^{1/4}}{n^{3/4}}
+C_3\lambda_j\cdot \|{Z}_2\|
\end{align*}
, where $C_2, C_3>0$ are absolute constants. 
Together with Lemma~\ref{lem:E_op}, we conclude that with probability at least $1-\frac{1}{\bar{m}^2}-3e^{-\bar{m}}-e^{-n}$, 
\begin{align*}
\big|\hat\lambda_j^{-2}-\lambda_j^{-2} \big|=\frac{|\lambda_j^2-\hat\lambda_j^2|}{\hat\lambda_j^2\lambda_j^2}\leq& C_2\frac{|\hat\lambda_j^2-\lambda_j^2|}{\lambda_j^4}\\
\leq&C_2\beta^2\cdot \frac{\sigma_\xi^2\bar{m}^2}{\lambda_j^4n^2}+C_2\frac{\sigma_\xi^{3/2}}{\lambda_j^{7/2}}\cdot \frac{\bar{m}^{1/4}}{n^{3/4}}+C_2\frac{\sigma_\xi}{\lambda_j^3}\cdot \frac{r^{1/2}\bar{m}\log^{1/2}\bar{m}}{n}
\end{align*}
for all $1\leq j\leq r$. Therefore, with the same probability, we get
\begin{align*}
\big|\hat B_n-\|\Lambda^{-1}\|_{\rm F}^2 \big|\leq C_2\|\Lambda^{-1}\|_{\rm F}^2\cdot& \Big(\frac{\beta^2\bar{m}^2}{n^2}+\frac{\beta^{3/2}\bar{m}^{1/4}}{n^{3/4}}+\frac{\beta r^{1/2}\bar{m}\log^{1/2}\bar{m}}{n}\Big)\\
\leq&C_2\|\Lambda^{-1}\|_{\rm F}^2(\beta\vee 1)^2\cdot \frac{r^{1/2}\bar{m}\log^{1/2}\bar{m}}{n}
\end{align*}
for some absolute constant $C_2>0$. Similarly, with the same probability, we get
\begin{align*}
\big|\hat V_n-\|\Lambda^{-2}\|_{\rm F}^2 \big|\leq C_2\|\Lambda^{-2}\|_{\rm F}^2(\beta\vee 1)^2\cdot \frac{r^{1/2}\bar{m}\log^{1/2}\bar{m}}{n}.
\end{align*}

\subsubsection{Proof of Theorem~\ref{thm:hatT_bUbV-normal-approx}}

By definition of $\hat T_{{U}{V}}$, we write
\begin{align*}
\hat T_{{U}{V}}:=&\frac{\|\calP_{\hat{U}\hat{V}}-\calP_{{U}{V}}\|_{\rm F}^2-\sigma_\xi^2\|{\Lambda}^{-1}\|_{\rm F}^2\cdot \frac{2m_{\star}}{n}}{\sqrt{8}\sigma_\xi^2\|{\Lambda}^{-2}\|_{\rm F}\cdot \frac{m_{\star}^{1/2}}{n}}+\underbrace{\frac{\sigma_\xi^2\big(\|{\Lambda}^{-1}\|_{\rm F}^2-\hat B_n\big)\cdot\frac{2m_{\star}}{n}}{\sqrt{8}\hat V_n^{1/2}\sigma_\xi^2\cdot\frac{m_{\star}^{1/2}}{n}}}_{\Xi_1}\\
&\quad+\underbrace{\frac{\|\calP_{\hat{U}\hat{V}}-\calP_{{U}{V}}\|_{\rm F}^2-\sigma_\xi^2\|{\Lambda}^{-1}\|_{\rm F}^2\cdot \frac{2m_{\star}}{n}}{\sqrt{8}\sigma_\xi^2\|{\Lambda}^{-2}\|_{\rm F}\cdot \frac{m_{\star}^{1/2}}{n}}\cdot \bigg[\frac{\|{\Lambda}^{-2}\|_{\rm F}}{\hat V_n^{1/2}}-1\bigg]}_{\Xi_2}\\
&\quad+\underbrace{\frac{\|\calP_{\hat U\hat V}-\calP_{UV}\|_{\rm F}^2-\hat B_n \hat\sigma_\xi^2\cdot \frac{2m_{\star}}{n}}{\sqrt{8}\hat V_n^{1/2}\hat\sigma_\xi^2\frac{\sqrt{m_{\star}}}{n}}-\frac{\|\calP_{\hat U\hat V}-\calP_{UV}\|_{\rm F}^2-\hat B_n \sigma_\xi^2\cdot \frac{2m_{\star}}{n}}{\sqrt{8}\hat V_n^{1/2}\sigma_\xi^2\frac{\sqrt{m_{\star}}}{n}}}_{\Xi_3}.
\end{align*}
By Lemma~\ref{lem:V_n-B_n}, we get that with probability at least $1-\bar{m}^{-2}-c_1e^{-c_2\bar{m}}$, $\hat V_n\geq \frac{\|{\Lambda}^{-2}\|_{\rm F}^2}{2}$ as long as $\frac{n}{(\beta\vee 1)^2}\geq C_1r\bar{m}\log \bar{m}$ for large enough $C_1>0$. Therefore, by Lemma~\ref{lem:V_n-B_n}, with the same probability,
\begin{align*}
\big|\Xi_1 \big|\leq& C_6(\beta\vee 1)^2 \frac{\|\Lambda^{-1}\|_{\rm F}^2}{\|\Lambda^{-2}\|_{\rm F}}\cdot \frac{r^{1/2}\bar{m}^{3/2}\log^{1/2}\bar{m}}{n}\leq C_6(\beta\vee 1)^2\frac{r\bar{m}^{3/2}\log^{1/2}\bar{m}}{n}
\end{align*}
where we used the fact $\|{\Lambda}^{-1}\|_{\rm F}^2\leq r^{1/2}\|{\Lambda}^{-2}\|_{\rm F}$. 
By Lemma~\ref{lem:V_n-B_n}, with the same probability,
\begin{align*}
\bigg|\frac{\|{\Lambda}^{-2}\|_{\rm F}}{\hat V_n^{1/2}}-1\bigg|\leq C_6(\beta\vee 1)^2\cdot \frac{r^{1/2}\bar{m}\log^{1/2}\bar{m}}{n}.
\end{align*}
By Theorem~\ref{thm:hatPN-PN-frob-con}, with probability at least $1-\frac{2n+10}{n^2}-3e^{-n}-c_1ne^{-c_2m_{\star}}-\frac{1}{\bar{m}^2}$,
\begin{align*}
\bigg|&\frac{\|\calP_{\hat{U}\hat{V}}-\calP_{{U}{V}}\|_{\rm F}^2-\sigma_\xi^2\|{\Lambda}^{-1}\|_{\rm F}^2\cdot \frac{2m_{\star}}{n}}{\sqrt{8}\sigma_\xi^2\|{\Lambda}^{-2}\|_{\rm F}\cdot \frac{m_{\star}^{1/2}}{n}}\bigg|\\
\leq& C_1\frac{\sigma_\xi}{\|{\Lambda}^{-2}\|_{\rm F}\lambda_r^3}\cdot\frac{r\bar{m}^{1/2}\log^{1/2}n}{n^{1/2}} +C_2\log^{1/2}n+C_3\frac{r^{3/2}\bar{m}^{3/2}\log^{1/2} n}{n}.
\end{align*}
Therefore, with probability at least $1-\frac{2n+11}{n^2}-4e^{-n}-c_1ne^{-c_2\bar{m}}-\frac{1}{\bar{m}^2}$, we get that
\begin{align*}
\big|\Xi_2\big|\leq&C_7(\beta\vee 1)^4\cdot \bigg(\frac{r^{1/2}\bar{m}}{n}+\frac{r^{3/2}\bar{m}^{3/2}}{n^{3/2}}+\frac{r^2\bar{m}^{5/2}}{n^2}\bigg)\log n.
\end{align*}
Moreover, by Lemma~\ref{lem:V_n-B_n}, we get that with probability at least $1-\frac{2}{n^2}$, 
\begin{align*}
\big|\Xi_3 \big|=\frac{\|\calP_{\hat U\hat V}-\calP_{UV}\|_{\rm F}^2}{\sqrt{8}\hat V_n^{1/2}\hat\sigma_\xi^2\sigma_\xi^2\frac{\sqrt{m_{\star}}}{n}}\cdot \big|\hat\sigma_\xi^2-\sigma_\xi^2\big|\leq \frac{\|\calP_{\hat U\hat V}-\calP_{UV}\|_{\rm F}^2}{\sqrt{8}\hat V_n^{1/2}\hat\sigma_\xi^2\sigma_\xi^2\frac{\sqrt{m_{\star}}}{n}}\cdot \Big(\|\Delta\|_{\rm F}^2+\frac{\sigma_\xi^2\log n}{\sqrt{n}}\Big).
\end{align*}
By the simple fact $\|\calP_{\hat U\hat V}-\calP_{UV}\|_{\rm F}^2\leq C_1\frac{\sigma_\xi^2}{\lambda_r^2}\cdot \frac{r\bar{m}}{n}$ which holds with probability at least $1-C_1e^{-c_2\bar{m}}$ for some absolute constants $C_1,c_1,c_2>0$. Therefore, we conclude with 
$$
\big|\Xi_3\big|\leq C_4(\beta\vee 1)^2\frac{r\bar{m}^{3/2}}{n}+C_5(\beta\vee 1)^2\frac{r\bar{m}^{1/2}\log n}{n^{1/2}}
$$
for constants $C_4,C_5>0$ depending on $C_1,C_2,C_3$.
Together with Corollary~\ref{cor:hatP-normal_approx}, we obtain 
\begin{align*}
\sup_x\Big|\PP\big\{&\hat T_{{U}{V}}\leq x\big\}-\Phi(x)\Big|\\
\leq& C_7(\beta\vee 1)^4\cdot \bigg(\frac{r\bar{m}^{1/2}\log^{3/2}n}{n^{1/2}}+\frac{r^{3/2}\bar{m}^{3/2}\log n}{n}\bigg)+6e^{-n}+(2n+r)e^{-m_{\star}}+\frac{5n+17}{n^2}\\
&\hspace{5cm}+c_1e^{-c_2\bar{m}}+\frac{C_8}{\bar{m}^{1/2}}
\end{align*}
for absolute constants $c_1,c_2,C_7, C_8>0$.

\section{Proof of additional lemmas}\label{sec:proof_additional}
The following lemmas will be frequently used through our proof. Basically, the Gaussian isoperimetric inequality can provide us with tight concentration bounds for Lipschitz functions.
\begin{lemma}\label{lem:gaussian_con}
 Let ${X}_1,\ldots,{X}_n\in\RR^{m}$ be i.i.d. centered Gaussian random vector with $\bSigma=\EE {X}{X}^{\top}$. 
Let $h(\cdot)$  be a function $\mathbb{R}^{nm}\mapsto \mathbb{R}$ satisfying the following 
 Lipschitz condition with some constant $L>0:$
 \begin{align*}
 |h(\{X_i\}_{i=1}^n)-h(\{X_i'\}_{i=1}^n)|\leq L\Big(\sum_{i=1}^n\|X_i-&X_i'\|_{\ell_2}^2\Big)^{1/2}, \\
 &\forall X_1,\cdots,X_n,X_1',\cdots,X_n'\in\RR^m.
 \end{align*}
Then, there exists some constant $C_1>0$ such that for all $t\geq 1$,
\begin{equation*}
 \mathbb{P}\Big\{\big| h(\{{X}_i\}_{i=1}^n)-\mathbb{E}f(\{{X}\}_{i=1}^n)\big|\geq C_1L\|\bSigma\| t^{1/2}\Big\}\leq e^{-t}.
\end{equation*}
\end{lemma}

\subsubsection{Proof of Lemma~\ref{lem:E_op}}
Recall that ${E}_1=\mathfrak{D}({Z}_1)$ with ${Z}_1=\frac{1}{n}\sum_{i=n+1}^{2n}\xi_i{X}_i$. Therefore, $\|{E}_1\|=\|{Z}_1\|$. Meanwhile, conditional on $\{\xi_i\}_{i=n+1}^{2n}$, 
$$
{Z}_1\stackrel{{\rm d}}{=}{X}\cdot \frac{\sqrt{\sum_{i=n+1}^{2n}\xi_i^2}}{n}
$$
where ${X}$ has i.i.d. standard Gaussian entries. By \cite{bandeira2016sharp}, 
$$
\EE_{{X}}\|{Z}_1\|\leq C_1\frac{\sqrt{\sum_{i=n+1}^{2n}\xi_i^2}}{n}\bar{m}^{1/2}.
$$
By Jensen's inequality, we get
\begin{gather*}
\EE\|{E}_1\|=\EE_{\xi}\EE_{X}\|{Z}_1\|\leq C_1\EE_{\xi}\frac{\sqrt{\sum_{i=n+1}^{2n}\xi_i^2}}{n}\bar{m}^{1/2}\\
\leq C_1\frac{\bar{m}^{1/2}}{n}\Big(\EE\sum_{i=n+1}^{2n}\xi_i^2 \Big)^{1/2}=C_1\sigma_{\xi}\frac{\bar{m}^{1/2}}{n^{1/2}}.
\end{gather*}
Conditional on $\{\xi_i\}_{i=n+1}^{2n}$, we view $\|{Z}_1\|$ as a function of $\{{X}_i\}_{i=n+1}^{2n}$, i.e.,
$$
h\big(\{{X}_i\}_{i=n+1}^{2n}\big)=\|{Z}_1\|=\Big\| n^{-1}\sum_{i=n+1}^{2n}\xi_i{X}_i\Big\|.
$$
Clearly, we have
\begin{gather*}
\Big|h\big(\{{X}_i\}_{i=n+1}^{2n}\big)-h\big(\{{X}_i'\}_{i=n+1}^{2n}\big)\Big|\leq \Big\| n^{-1}\sum_{i=n+1}^{2n}\xi_i({X}_i-{X}_i')\Big\|\\
\leq \frac{1}{n}\Big(\sum_{i=n+1}^{2n}\xi_i^2\Big)^{1/2}\Big(\sum_{i=n+1}^{2n}\|{X}_i-{X}_i'\|_{\rm F}^2\Big)^{1/2}
\end{gather*}
implying that $h(\cdot)$ is Lipschitz with constant $n^{-1}\big(\sum_{i=n+1}^{2n}\xi_i^2\big)^{1/2}$. By Lemma~\ref{lem:gaussian_con}, we get that with probability at least $1-e^{-t}$ for $t\geq 1$,
\begin{gather*}
\big|\|{E}_1\| -\EE_{{X}}\|{E}_1\|\big|\leq\frac{C_1}{n}\Big(\sum_{i=n+1}^{2n}\xi_i^2\Big)^{1/2}t^{1/2}
%\leq \frac{C_2}{n}\Big(\sum_{i=n+1}^{2n}\xi_i^2\Big)^{1/2}\bar{m}^{1/2}+\frac{C_1}{n}\Big(\sum_{i=n+1}^{2n}\xi_i^2\Big)^{1/2}t^{1/2}.
\end{gather*}
Since $\{\xi_i\}_i$ are i.i.d. Gaussian random variables, we get $\PP\Big(\sum_{i=n+1}^{2n}\xi_i^2\leq c_1n\sigma_\xi^2\Big)\geq 1-e^{-n}$ for some absolute constant $c_1>0$. We conclude that with probability at least $1-e^{-t}-e^{-n}$,
\begin{gather}\label{eq:E_1-EEXE_1}
\big|\|{E}_1\|-\EE_{X}\|{E}_1\|\big|\leq C_2\sigma_{\xi}\frac{t^{1/2}}{n^{1/2}}.
\end{gather}
We then view $\EE_{{X}}\|{E}_1\|$ as a function $\{\xi_i\}_i$, i.e., 
$$
h_1\big(\{\xi_i\}_{i=n+1}^{2n}\big)=\EE_{X}\Big\|\frac{\big(\sum_{i=n+1}^{2n}\xi_i^2\big)^{1/2}}{n}{X}\Big\|.
$$
Then, by denoting $\bxi=(\xi_{n+1},\cdots,\xi_{2n})^{\top}\in\RR^n$, we get
\begin{align*}
\big| h_1\big(\{\xi_i\}_{i=n+1}^{2n}\big)-h_1\big(\{\xi_i'\}_{i=n+1}^{2n}\big)\big|\leq \|\bxi-\bxi'\|_{\ell_2}\cdot \frac{\EE_{{X}}\|{X}\|}{n}\leq C_1\frac{\bar{m}^{1/2}}{n}\cdot \|\bxi-\bxi'\|_{\ell_2}.
\end{align*}
By Lemma~\ref{lem:gaussian_con}, we get that with probability at least $1-e^{-t}$ for all $t\geq 1$,
\begin{align}\label{eq:EEE_1-EEXE_1}
\big|\EE\|{E}_1\|-\EE_{X}\|{E}_1\|\big|\leq C_1\sigma_\xi\cdot \frac{\bar{m}^{1/2}t^{1/2}}{n}.
\end{align}
By (\ref{eq:E_1-EEXE_1}) and (\ref{eq:EEE_1-EEXE_1}), we conclude that with probability at least $1-2e^{-t}-e^{-n}$,
\begin{align*}
\big|\|{E}_1\|-\EE\|{E}_1\| \big|\leq C_1\sigma_\xi\cdot\bigg[ \frac{t^{1/2}}{n^{1/2}}+\frac{\bar{m}^{1/2}t^{1/2}}{n}\bigg].
\end{align*}
Now, we turn to the proof of $\EE\|{E}_2\|$.
Recall that $\|{E}_2\|=\|{Z}_2\|$ where 
$$
{Z}_2=n^{-1}\sum_{i=n+1}^{2n}\big(\langle\bDelta,{X}_i \rangle{X}_i-\bDelta\big).
$$
The following bounds are standard
$$
\|\langle\bDelta,{X} \rangle{X}\|_{\psi_1}\lesssim \|\langle\bDelta,{X} \rangle\|_{\psi_2}\cdot\big\|\|{X}\|\big\|_{\psi_2}\lesssim \|\bDelta\|_{\rm F}\bar{m}^{1/2}
$$
where the Orlicz $\psi_{\alpha}$-norm, for $\alpha\in[1,2]$, of a random variable $X$ is defined as
$$
\|X\|_{\psi_{\alpha}}:=\inf\big\{u>0: \EE\exp(|X|^{\alpha}/u^{\alpha})\leq 2\big\}.
$$
By matrix Bernstein inequality \cite{koltchinskii2011neumann}, with probability at least $1-e^{-t}$ for $t\geq 0$, we have
$$
\|{Z}_2\|\leq C_1\|\bDelta\|_{\rm F}\sqrt{\frac{\bar{m}(t+\log\bar{m})}{n}}+C_2\|\bDelta\|_{\rm F}\frac{\bar{m}^{1/2}(t+\log \bar{m})}{n}.
$$
By integrating over $t$, as long as $n\geq \log \bar{m}$, we end up with
$$
\EE\|{E}_2\|=\EE\|{Z}_2\|\leq C_1\|\bDelta\|_{\rm F}\frac{m^{1/2}\log^{1/2}\bar{m}}{n^{1/2}}.
$$
We denote by ${\rm vec}(\bDelta)$ the vectorization of $\bDelta$ and $\calM({v})$ the matricization of a vector ${v}\in\RR^{m_1m_2}$ such that $\calM({\rm vec}(\bDelta))=\bDelta$. We write ${Z}_2={Z}_{21}+{Z}_{22}$ with
\begin{gather*}
{Z}_{21}:=\frac{1}{n}\sum_{i=n+1}^{2n}\big(\langle\bDelta,{X}_i \rangle\calM(\calP_{{\rm vec}(\bDelta)}{\rm vec}({X}_i))-\bDelta\big)\\
{Z}_{22}:=\frac{1}{n}\sum_{i=n+1}^{2n}\langle\bDelta,{X}_i \rangle\calM\big(\calP_{{\rm vec}(\bDelta)}^{\perp}{\rm vec}({X}_i)\big)
\end{gather*}
where $\calP_{{v}}$ denotes the orthogonal projection onto ${v}$, i.e., $\calP_{{v}}({u})=\frac{{v}\cdot({v}^{\top}{u})}{\|v\|_{\ell_2}^2}$. More explicitly, we have
$$
\calP_{{\rm vec}(\bDelta)}{\rm vec}({X}_i)={\rm vec}(\bDelta)\cdot\frac{\langle\bDelta,{X}_i \rangle}{\|\bDelta\|_{\rm F}^2}.
$$
Since $\langle\bDelta,{X}_i \rangle$ and $\calM\big(\calP_{{\rm vec}(\bDelta)}^{\perp}{\rm vec}({X}_i)\big)$ both have Gaussian distributions, we claim that $\langle\bDelta,{X}_i \rangle$ is independent with $\calM\big(\calP_{{\rm vec}(\bDelta)}^{\perp}{\rm vec}({X}_i)\big)$ in view of  their uncorrelation. %Now,
%$$
%\|{Z}_1\|\leq \|{Z}_{11}\|+\|{Z}_{12}\|.
%$$

We view $\|{Z}_{22}\|$ as a function of $\big\{\calM\big(\calP_{{\rm vec}(\bDelta)}^{\perp}{\rm vec}({X}_i)\big)\big\}_{i=n+1}^{2n}$, conditioned on $\{\langle\bDelta,{X}_i \rangle\}_{i=n+1}^{2n}$. More exactly, we define
$$
h_2\Big(\big\{\calM\big(\calP_{{\rm vec}(\bDelta)}^{\perp}{\rm vec}({X}_i)\big)\big\}_{i=n+1}^{2n}\Big)=\Big\|\frac{1}{n}\sum_{i=n+1}^{2n}\langle\bDelta,{X}_i \rangle\calM\big(\calP_{{\rm vec}(\bDelta)}^{\perp}{\rm vec}({X}_i)\big)\Big\|.
$$
Observe that $h_2\big(\cdot\big)$ is a Lipschitz function with constant $n^{-1}\Big(\sum_{i=1}^n\langle\bDelta,{X}_i \rangle^2\Big)^{1/2}$. By Lemma~\ref{lem:gaussian_con}, conditioned on $\{\langle\Delta,X_i \rangle\}_{i=n+1}^{2n}$, we get that with probability at least $1-e^{-t}$ for all $t\geq 1$,
\begin{gather*}
\|{Z}_{22}\|-\EE_{\{\calP_{{\rm vec}(\bDelta){\rm vec}({X_i})}^{\perp}\}_i}\|{Z}_{22}\|\leq \frac{C_1t^{1/2}}{n}\Big(\sum_{i=n+1}^{2n}\langle\bDelta,{X}_i \rangle\Big)^{1/2}.
\end{gather*}
Similarly, we conclude that with probability at least $1-e^{-t}-e^{-n}$,
\begin{gather}\label{eq:Z22-EEperpZ22}
\Big|\|{Z}_{22}\|-\EE_{\{\calP_{{\rm vec}(\bDelta){\rm vec}({X_i})}^{\perp}\}_i}\|{Z}_{22}\|\Big|\leq C_1\|\bDelta\|_{\rm F}\cdot\frac{t^{1/2}}{n^{1/2}}.
\end{gather}
Following the same fashion, we view $\EE_{\{\calP_{{\rm vec}(\bDelta){\rm vec}({X_i})}^{\perp}\}_i}\|{Z}_{22}\|$ as a function of $\{\langle \bDelta,{X}_i\rangle\}_i$ and define
\begin{align*}
h_3\big(\big\{\langle\bDelta,{X}_i \rangle\big\}_{i=n+1}^{2n}\big)=&\EE_{\{\calP_{{\rm vec}(\bDelta){\rm vec}({X}_i)}^{\perp}\}_{i}}\Big\|\frac{1}{n}\sum_{i=n+1}^{2n}\langle\bDelta,{X}_i \rangle\calM\big(\calP_{{\rm vec}(\bDelta)}^{\perp}{\rm vec}({X}_i)\big)\Big\|\\
\stackrel{{\rm d}}{=}&\EE_{\calP_{{\rm vec}(\bDelta){\rm vec}({\tilde X})}^{\perp}}\Big\|\frac{1}{n}\Big(\sum_{i=n+1}^{2n}\langle\bDelta,{X}_i \rangle^2\Big)^{1/2}\calM\big(\calP_{{\rm vec}(\bDelta)}^{\perp}{\rm vec}({\tilde X})\big)\Big\|
\end{align*}
where ${\tilde X}$ is an independent copy of ${X}_i$. Denote the vector $\bx_{\bDelta}=\big(\langle\bDelta, {X}_i\rangle\big)_{i=n+1}^{2n}\in\RR^n$. Then,
\begin{align*}
\Big|h_3\big(\big\{\langle\bDelta,{X}_i \rangle\big\}_{i=n+1}^{2n}\big)-&h_3\big(\big\{\langle\bDelta,{X}_i' \rangle\big\}_{i=n+1}^{2n}\big)\Big|\\
\leq&\|\bx_{\bDelta}-\bx_{\bDelta}'\|_{\ell_2}\cdot \frac{1}{n}\EE_{\calP_{{\rm vec}(\bDelta){\rm vec}({\tilde X})}^{\perp}}\Big\|\calM\big(\calP_{{\rm vec}(\bDelta)}^{\perp}{\rm vec}({\tilde X})\big)\Big\|\\
\leq&C_1\|\bx_{\bDelta}-\bx_{\bDelta}'\|_{\ell_2}\cdot \frac{\bar{m}^{1/2}}{n}.
\end{align*}
Therefore, by Lemma~\ref{lem:gaussian_con}, we get that with probability at least $1-e^{-t}$,
\begin{align}\label{eq:EZ22-EEperpZ22}
\Big|\EE\|{Z}_{22}\|-\EE_{\{\calP_{{\rm vec}(\bDelta){\rm vec}({X_i})}^{\perp}\}_i}\|{Z}_{22}\|\Big|\leq C_1\|\bDelta\|_{\rm F}\cdot \frac{\bar{m}^{1/2}t^{1/2}}{n}.
\end{align}
By (\ref{eq:Z22-EEperpZ22}) and (\ref{eq:EZ22-EEperpZ22}), we conclude that with probability at least $1-2e^{-t}-e^{-n}$,
\begin{align*}
\Big|\|{Z}_{22}\|-\EE\|{Z}_{22}\| \Big|\leq C_1\|\bDelta\|_{\rm F}\cdot \bigg[\frac{t^{1/2}}{n^{1/2}}+\frac{\bar{m}^{1/2}t^{1/2}}{n}\bigg].
\end{align*}
Similarly, by matrix Bernstein inequality (\cite{koltchinskii2011neumann}), we conclude that, with probability at least $1-e^{-t}$,
\begin{gather*}
\|{Z}_{21}\|\leq C_1\|\bDelta\|_{\rm F}\bigg(\frac{(t+\log\bar{m})^{1/2}}{n^{1/2}}+\frac{t+\log\bar{m}}{n}\bigg)
\end{gather*}
and thus
$$
\EE\|{Z}_{21}\|\leq C_1\|\bDelta\|_{\rm F}\leq C_1\|\bDelta\|_{\rm F}\cdot \frac{\log^{1/2}\bar{m}}{n^{1/2}}.
$$
By putting the above three bounds together and adjusting the constants, we obtain
\begin{align*}
\big|\|{E}_2\|-\EE\|{E}_2\|\big|\leq C_1\|\bDelta\|_{\rm F}\cdot\bigg[\frac{t^{1/2}+\log^{1/2}\bar{m}}{n^{1/2}}+\frac{\bar{m}^{1/2}t^{1/2}+t+\log\bar{m}}{n}\bigg]
\end{align*}
%$$
%\|{E}_1\|=\|{Z}_1\|\leq C_1\|\bDelta\|_{\rm F}\frac{(\bar{m}\log\bar{m})^{1/2}}{n^{1/2}}+C_2\|\bDelta\|_{\rm F}\bigg(\frac{t^{1/2}}{n^{1/2}}+\frac{t}{n}\bigg)
%$$
with probability at least $1-3e^{-t}-e^{-n}$ for all $t\geq 1$.

\subsubsection{Proof of Lemma~\ref{lemma:varphi-delta_n}}
Recall from eq. (\ref{eq:SN_k}) that
\begin{equation*}
\begin{split}
 \varphi_{k,\bar{\delta}}({E})&=-2\big<\calS_{{N},k}({E}),\calP_{{U}{V}}\big>\phi\Big(\frac{\|{E}\|}{\bar{\delta}}\Big)
%=&\sum_{1\leq |k|\leq r}(-1)^k\frac{1}{(\mu_i-\eta)^2}\tr\big(P_iE_1\big[\calR_{N}(\eta)E_1\big]^{k-1}P_i\big)\phi\big(\frac{\|E_1\|_{\rm op}}{\bar{\delta}}\big).
\end{split}
\end{equation*}
where the matrix $\calS_{{N},k}({E})$ is defined with non-negative integers $\{s_i\}_{i=1}^{k+1}$ so that 
$$
\calS_{{N},k}({E})=\sum_{\bs: s_1+\cdots+s_{k+1}=k}(-1)^{1+\tau(\bs)}\cdot \calP_{{U}{V}}^{-s_1}{E}\calP_{{U}{V}}^{-s_{2}}{E}\cdots{E}\calP_{{U}{V}}^{-s_{k+1}}.
$$
\paragraph*{Case 1} If $\|{E}\|, \|{E}'\|\geq \frac{9}{8}\bar{\delta}$, then $\phi\big(\|{E}\|/\bar{\delta}\big)=\phi\big(\|{E}'\|/\bar{\delta}\big)=0$. The first claim bound trivially holds. 

\paragraph*{Case 2} If $\|{E}\|, \|{E}'\|\leq \frac{9}{8}\cdot\bar{\delta}$, then for $s_1, s_{k+1}\geq 1$ and $s_2,\cdots, s_k\geq 0$, we get
\begin{align*}
\Big|\tr\big(\calP_{{U}{V}}^{-s_1}{E}&\calP_{{U}{V}}^{-s_{2}}{E}\cdots{E}\calP_{{U}{V}}^{-s_{k+1}}\big)\phi\Big(\frac{\|{E}\|}{\bar{\delta}}\Big)-\tr\big(\calP_{{U}{V}}^{-s_1}{E}'\calP_{{U}{V}}^{-s_{2}}{E}'\cdots{E}'\calP_{{U}{V}}^{-s_{k+1}}\big)\phi\Big(\frac{\|{E}'\|}{\bar{\delta}}\Big)\Big|\\
\leq &2kr\|{E}-{E}'\|\cdot \frac{\bar\delta^{k-1}}{\lambda_r^{k}}\Big(\frac{9}{8}\Big)^{k-1}
+\big|\tr\big(\calP_{{U}{V}}^{-s_1}{E}'\calP_{{U}{V}}^{-s_{2}}{E}'\cdots{E}'\calP_{{U}{V}}^{-s_{k+1}}\big)\big|\bigg|\phi\Big(\frac{\|{E}\|}{\bar{\delta}}\Big)-\phi\Big(\frac{\|{E}'\|}{\bar{\delta}}\Big)\bigg|\\
\leq& k(2r)\big(\frac{1}{\lambda_r}\big)^{k}\Big(\frac{9\bar{\delta}}{8}\Big)^{k-1}\|{E}-{E}'\|
+2r\big(\frac{1}{\lambda_r}\big)^{k}\Big(\frac{9\bar\delta}{8}\Big)^{k}\frac{8}{\bar{\delta}}\|{E}-{E}'\|
\end{align*}
, where the last inequality is due to the Lipschitz property of function $\phi(\cdot)$. Therefore,
$$
\big|\varphi_{k,\bar{\delta}}({E})-\varphi_{k,\bar{\delta}}({E}')\big|\leq 8r\cdot\frac{(k+9)}{\lambda_r}\Big(\frac{9\bar{\delta}}{2\lambda_r}\Big)^{k-1}\|{E}-{E}'\|
$$
which proves the first claim.

\paragraph*{Case 3} If $\|{E}\|\leq \frac{9}{8}\cdot\bar{\delta}$ and $\|{E}'\|\geq \frac{9}{8}\cdot\bar{\delta}$, then $\phi\big(\|{E}'\|/\bar{\delta}\big)=0$. For $s_1, s_{k+1}\geq 1$ and $s_2,\cdots,s_k\geq0$, we write
\begin{align*}
\Big|\tr&\big(\calP_{{U}{V}}^{-s_1}{E}\calP_{{U}{V}}^{-s_{2}}{E}\cdots{E}\calP_{{U}{V}}^{-s_{k+1}}\big)\phi\Big(\frac{\|{E}\|}{\bar{\delta}}\Big)\Big|\\
=&\Big|\tr\big(\calP_{{U}{V}}^{-s_1}{E}\calP_{{U}{V}}^{-s_{2}}{E}\cdots{E}\calP_{{U}{V}}^{-s_{k+1}}\big)\phi\Big(\frac{\|{E}\|}{\bar{\delta}}\Big)-\tr\big(\calP_{{U}{V}}^{-s_1}{E}\calP_{{U}{V}}^{-s_{2}}{E}\cdots{E}\calP_{{U}{V}}^{-s_{k+1}}\big)\phi\Big(\frac{\|{E}\|}{\bar{\delta}}\Big)\Big|\\
\leq& \Big|\tr\big(\calP_{{U}{V}}^{-s_1}{E}\calP_{{U}{V}}^{-s_{2}}{E}\cdots{E}\calP_{{U}{V}}^{-s_{k+1}}\big)\phi\Big(\frac{\|{E}\|}{\bar{\delta}}\Big)\Big|\cdot\Big|\phi\big(\frac{\|{E}\|}{\bar{\delta}}\big)-\phi\big(\frac{\|{E}'\|}{\bar{\delta}}\big)\Big|\\
\leq& 16r\Big(\frac{1}{\lambda_r}\Big)^{k}\bar{\delta}^{-1}\Big(\frac{9\bar\delta}{8}\Big)^{k}\|{E}-{E}'\|.
\end{align*}
Therefore, we get 
\begin{gather*}
\big|\varphi_{k,\bar{\delta}}({E})-\varphi_{k,\bar{\delta}}({E}')\big|\leq 72\frac{r}{\lambda_r}\Big(\frac{9\bar{\delta}}{2\lambda_r}\Big)^{k-1}\|{E}-{E}'\|
\end{gather*}
which also proves the first claim. 

\paragraph*{Case 4} If $\|{E}'\|\leq \frac{9}{8}\cdot\bar{\delta}$ and $\|{E}\|\geq \frac{9}{8}\cdot\bar{\delta}$. The proof is identical to {\it Case 3}.

\paragraph*{Proof of second claim} By first claim,
\begin{align*}
 \big|\varphi_{\bar{\delta}}&({E})-\varphi_{\bar{\delta}}({E}')\big|\leq \sum_{k\geq 3}\big|\varphi_{k,\bar{\delta}}({E})-\varphi_{k,\bar{\delta}}({E}')\big|\\
\leq& \frac{8r}{\lambda_r}\sum_{k\geq 3}(k+9)\Big(\frac{9\bar{\delta}}{2\lambda_r}\Big)^{k-1}\|{E}-{E}'\|
\leq \frac{C_5r}{\lambda_r}\Big(\frac{9\bar{\delta}}{2\lambda_r}\Big)^2\|{E}-{E}'\|
\end{align*}
for an absolute constant $C_5>0$.

\subsubsection{Proof of Lemma~\ref{lem:Eperp-E}}
To this end, we define
\begin{gather*}
{A}=\frac{1}{n}\sum_{i=n+1}^{2n}\big(\xi_i+\langle\bDelta,{X}_i \rangle\big)\mathfrak{D}\circ\calM\big(\calP_{{\rm vec}(\bDelta)}{\rm vec}({X}_i)\big)-\bDelta
\end{gather*}
and
\begin{gather*}
{A}^{\perp}=\frac{1}{n}\sum_{i=n+1}^{2n}\big(\xi_i+\langle\bDelta,{X}_i \rangle\big)\mathfrak{D}\circ\calM\big(\calP_{{\rm vec}(\bDelta)}^{\perp}{\rm vec}({X}_i)\big)
\end{gather*}
such that ${E}={A}+{A}^{\perp}$. Then, we write
\begin{align*}
 \big|\mathbb{E}_{\{\calP^{\perp}_{{\rm vec}(\bDelta)}{X_i}\}_i}&\big[\varphi_{\bar{\delta}}({A}+{A}^{\perp})\big]-\EE\big[\varphi_{\bar{\delta}}({A}+{A}^{\perp})\big]\big|\\
 \leq & \big|\mathbb{E}_{\{\calP^{\perp}_{{\rm vec}(\bDelta)}{X_i}\}_i}\big[\varphi_{\bar{\delta}}({A}+{A}^{\perp})\big]-\mathbb{E}_{\{\calP^{\perp}_{{\rm vec}(\bDelta)}{X_i}\}_i}\big[\varphi_{\bar{\delta}}({A}^{\perp})\big]\big|\\
 &+\big|\mathbb{E}_{\{\calP^{\perp}_{{\rm vec}(\bDelta)}{X_i}\}_i}\big[\varphi_{\bar{\delta}}({A}^{\perp})\big]-\EE\big[\varphi_{\bar{\delta}}({A}^{\perp})\big]\big|+\big|\EE\big[\varphi_{\bar{\delta}}({A}^{\perp})\big]-\EE\big[\varphi_{\bar{\delta}}({A}+{A}^{\perp})\big]\big|.
\end{align*}
By Lemma~\ref{lemma:varphi-delta_n}, 
\begin{align*}
 \big|\mathbb{E}_{\{\calP^{\perp}_{{\rm vec}(\bDelta)}{X_i}\}_i}&\big[\varphi_{\bar{\delta}}({A}+{A}^{\perp})\big]-\mathbb{E}_{\{\calP^{\perp}_{{\rm vec}(\bDelta)}{X_i}\}_i}\big[\varphi_{\bar{\delta}}({A}^{\perp})\big]\big|\\
 \leq& \mathbb{E}_{\{\calP^{\perp}_{{\rm vec}(\bDelta)}{X_i}\}_i}\big| \varphi_{\bar{\delta}}({A}+{A}^{\perp})-\varphi_{\bar{\delta}}({A}^{\perp})\big|
 \leq\frac{C_6r}{\lambda_r}\Big(\frac{9\bar\delta}{2\lambda_r}\Big)^2\|{A}\|.
\end{align*}
By matrix Bernstein inequality (\cite{koltchinskii2011neumann}), we conclude that (see also Proposition~\ref{prop:nuclear_pen}) with probability at least $1-e^{-t}-c_1e^{-c_2\bar{m}}$,
\begin{align*}
\|{A}\|\leq C_5\sigma_{\xi}\Big(\frac{(t+\log\bar{m})^{1/2}}{n^{1/2}}+\frac{t+\log\bar{m}}{n}\Big)
\end{align*}
where the first term dominate if $t\leq n$ and $n\geq \log\bar{m}$.
Therefore, with probability at least $1-e^{-t}-c_1e^{-c_2\bar{m}}$ for $1\leq t\leq n$, we have
\begin{align*}
\big|&\mathbb{E}_{\{\calP^{\perp}_{{\rm vec}(\bDelta)}{X_i}\}_i}\big[\varphi_{\bar{\delta}}({A}+{A}^{\perp})\big]-\mathbb{E}_{\{\calP^{\perp}_{{\rm vec}(\bDelta)}{X_i}\}_i}\big[\varphi_{\bar{\delta}}({A}^{\perp})\big]\big|
\leq \frac{C_5r}{\lambda_r}\Big(\frac{9\bar{\delta}}{2\lambda_r}\Big)^2\sigma_{\xi}\cdot \frac{(t+\log\bar{m})^{1/2}}{n^{1/2}}.
\end{align*}
Similarly, by integrating out $t$, we get
\begin{align*}
\big|\EE\big[&\varphi_{\bar{\delta}}({A}^{\perp})\big]-\EE\big[\varphi_{\bar{\delta}}({A}+{A}^{\perp})\big]\big|\leq C_6\frac{r}{\lambda_r}\Big(\frac{9\bar\delta}{2\lambda_r}\Big)^2\cdot \EE\|A\|^2
\leq\frac{C_6r}{\lambda_r}\Big(\frac{9\bar{\delta}}{2\lambda_r}\Big)^2\sigma_{\xi}\cdot \frac{\log^{1/2}\bar{m}}{n^{1/2}}.
\end{align*}
It remains to bound $\big|\mathbb{E}_{\{\calP^{\perp}_{{\rm vec}(\bDelta)}{X_i}\}_i}\big[\varphi_{\bar{\delta}}({A}^{\perp})\big]-\EE\big[\varphi_{\bar{\delta}}({A}^{\perp})\big]\big|$. Recall that $\{\xi_i\}_{i=n+1}^{2n}, \{\langle\bDelta,{X}_i \rangle\}_{i=n+1}^{2n}$ and $\{\calP_{{
\rm vec}(\bDelta)}^{\perp}{\rm vec}({X}_i)\}_i$ are mutually independent.
 Therefore, conditional on $\{\xi_i\}_i$ and $\{\langle\bDelta,{X}_i \rangle\}_i$, we have
\begin{gather*}
{A}^{\perp}\stackrel{{\rm d}}{=}\frac{1}{n}\Big(\sum_{i=n+1}^{2n}(\xi_i+\langle\bDelta,{X}_i \rangle)^2\Big)^{1/2}\mathfrak{D}\circ\calM\big(\calP_{{\rm vec}(\bDelta)}^{\perp}{\rm vec}(F)\big)
\end{gather*}
where $F$ is a copy of ${X}_i$ being independent with $\{\xi_i\}_i$ and $\{\langle\bDelta,{X}_i \rangle\}_i$. We define the function
\begin{align*}
h_5\big(\{\xi_i\}_i, &\{\langle\bDelta,{X}_i \rangle\}_i\big)=\mathbb{E}_{\{\calP^{\perp}_{{\rm vec}(\bDelta)}{X_i}\}_i}\big[\varphi_{\bar{\delta}}({A}^{\perp})\big]\\
=&\EE_{F}\varphi_{\bar{\delta}}\bigg(\frac{1}{n}\Big(\sum_{i=n+1}^{2n}(\xi_i+\langle\bDelta,{X}_i \rangle)^2\Big)^{1/2}\mathfrak{D}\circ\calM\big(\calP_{{\rm vec}(\bDelta)}^{\perp}{\rm vec}(F)\big)\bigg).
\end{align*}
By Lemma~\ref{lemma:varphi-delta_n}, we get
\begin{align*}
\big|h_5\big(\{\xi_i\}_i, &\{\langle\bDelta,{X}_i \rangle\}_i\big)-h_5\big(\{\xi_i'\}_i, \langle\bDelta,{X}_i' \rangle\}_i\big)\big|\\
\leq&\frac{C_6r}{\lambda_r}\Big(\frac{9\bar{\delta}}{2\lambda_r}\Big)^2\EE_F\Big\| \mathfrak{D}\circ\calM\big(\calP_{{\rm vec}(\bDelta)}^{\perp}{\rm vec}(F)\big)\Big\|\\
&\times\bigg| \frac{1}{n}\Big(\sum_{i=n+1}^{2n}(\xi_i+\langle\bDelta,{X}_i \rangle)^2\Big)^{1/2}-\frac{1}{n}\Big(\sum_{i=n+1}^{2n}(\xi_i'+\langle\bDelta,{X}_i' \rangle)^2\Big)^{1/2}\bigg|.
\end{align*}
Define the vectors $\bxi=(\xi_i)_{i=1}^n$ and $\bx_{\bDelta}=\big(\langle\bDelta,{X}_i \rangle\big)_{i=1}^n$. Then, we have
\begin{align*}
\bigg| \frac{1}{n}\Big(\sum_{i=n+1}^{2n}&(\xi_i+\langle\bDelta,{X}_i \rangle)^2\Big)^{1/2}-\frac{1}{n}\Big(\sum_{i=n+1}^{2n}(\xi_i'+\langle\bDelta,{X}_i' \rangle)^2\Big)^{1/2}\bigg|\\
=&\frac{1}{n}\big|\|\bxi+\bx_{\bDelta}\|_{\ell_2}-\|\bxi'+\bx'_{\bDelta}\|_{\ell_2} \big|
\leq\frac{1}{n}\big(\|\bxi-\bxi'\|_{\ell_2}+\|\bx_{\bDelta}-\bx'_{\bDelta}\|_{\ell_2}\big)\\
\leq&\frac{1}{n}\big(\|\bxi-\bxi'\|_{\ell_2}^2+\|\bx_{\bDelta}-\bx'_{\bDelta}\|_{\ell_2}^2\big)^{1/2}.
\end{align*}
Meanwhile, by operator norm of random matrix (\cite{vershynin2010introduction}), we can easily get
\begin{align*}
\EE_F\Big\| \mathfrak{D}\circ\calM&\big(\calP_{{\rm vec}(\bDelta)}^{\perp}{\rm vec}(F)\big)\Big\|
\leq \EE_{F}\|F\|+\EE_{F}\Big\| \mathfrak{D}\circ\calM\big(\calP_{{\rm vec}(\bDelta)}{\rm vec}(F)\big)\Big\|\leq C_1\bar{m}^{1/2}.
\end{align*}
Therefore, we conclude that $h_5\big(\cdot\big)$ is Lipschitz with respect to $\{\xi_i\}_i$ and $\{\langle\bDelta,{X}_i \rangle\}_i$ with constant $\frac{C_6r}{\lambda_r}\big(\frac{9\bar{\delta}}{2\lambda_r}\big)^2\frac{\bar{m}^{1/2}}{n}$. Since $\xi\sim\calN(0,\sigma_\xi^2)$ and $\langle\bDelta,{X} \rangle\sim\calN(0,\|\bDelta\|_{\rm F}^2)$, we apply Lemma~\ref{lem:gaussian_con} and conclude that with probability at least $1-e^{-t}-c_1e^{-c_2\bar{m}}$,
\begin{align*}
\big|\mathbb{E}_{\{\calP^{\perp}_{{\rm vec}(\bDelta)}{X_i}\}_i}&\big[\varphi_{\bar{\delta}}({A}^{\perp})\big]-\EE\big[\varphi_{\bar{\delta}}({A}^{\perp})\big]\big|\leq C_7(\sigma_\xi+\|\bDelta\|_{\rm F})\frac{rt^{1/2}}{\lambda_r}\Big(\frac{9\bar{\delta}}{2\lambda_r}\Big)^2\frac{\bar{m}^{1/2}}{n}\\
\leq &C_8r\frac{\sigma_\xi t^{1/2}}{\lambda_r}\Big(\frac{9\bar{\delta}}{2\lambda_r}\Big)^2\frac{\bar{m}^{1/2}}{n}.
\end{align*}
Therefore, with probability at least $1-2e^{-t}-2c_1e^{-c_2\bar{m}}$ for all $t\in[1,n]$,
\begin{align*}
 \big|\mathbb{E}_{\{\calP^{\perp}_{{\rm vec}(\bDelta)}{X_i}\}_i}&\big[\varphi_{\bar{\delta}}({A}+{A}^{\perp})\big]-\EE\big[\varphi_{\bar{\delta}}({A}+{A}^{\perp})\big]\big|
 \leq C_7r\frac{\sigma_\xi}{\lambda_r}\Big(\frac{9\bar{\delta}}{2\lambda_r}\Big)^2\cdot\frac{(t+\log\bar{m})^{1/2}}{n^{1/2}}
\end{align*}
where we used the fact $n\geq r\bar{m}$.

\bibliographystyle{abbrv}
\bibliography{refer}
\end{document}